\pdfoutput=1
\newif\ifOneColumn
\OneColumnfalse

\newif\iftcnsVersion
\tcnsVersiontrue

\ifOneColumn
\documentclass[journal,usenames,dvipsnames,onecolumn]{IEEEtran}
\else
\documentclass[journal,usenames,dvipsnames]{IEEEtran}
\fi 

\newif\ifArxivVersion
\ArxivVersiontrue

\newif\ifNotArxivVersion

\makeatletter\providecommand{\mdseries@tt}{}\makeatother

\usepackage{etex}
\reserveinserts{28}

\usepackage{makeidx}         
\usepackage{graphicx}        
\usepackage{multicol}        

\makeindex             

\usepackage{enumitem, fancyhdr, amsthm, amsmath, amsfonts,  amssymb, csquotes, mathtools, url, graphicx, balance, comment, caption, mathrsfs, upgreek, wrapfig, cases, pgfplots, subfig,  fancybox,  multirow, hhline, mathtools, chngcntr, booktabs, algpseudocode, mathabx, floatrow, algorithm, algorithmicx, mathrsfs, etex, etoolbox, xfrac,  xstring, xparse, ifmtarg, xifthen, bm, array, multicol, layouts, printlen,  soul, cancel, tabularx, ragged2e, relsize, scalefnt, varwidth, setspace, dsfont, bbm}

\algnewcommand\algorithmicforeach{\textbf{for each}}
\algdef{S}[FOR]{ForEach}[1]{\algorithmicforeach\ #1\ \algorithmicdo}

\newcolumntype{x}[1]{>{\centering\let\newline\\\arraybackslash\hspace{0pt}\vspace{1pt}}m{#1}}

\usepackage[]{todonotes}
\usepackage[numbers,sort&compress]{natbib}
\usepackage[normalem]{ulem}
\pgfplotsset{compat=1.13}

\usepackage{tikz,lmodern}
\usetikzlibrary{patterns,arrows.meta,shapes,snakes,automata,backgrounds,petri,fit,calc,positioning,external,intersections,trees}
\usepackage{color,xcolor}
\usepackage[colorlinks=true,
            raiselinks=true,
            linkcolor=MidnightBlue,
            citecolor=Mahogany,
            urlcolor=ForestGreen,
            pdftitle={Evaluating Resilience of Electricity Distribution Networks via A Modification of  Generalized Benders Decomposition Method},
            pdfkeywords={Cyber-physical security, network resilience, smart grids, Generalized Benders Decomposition},
            pdfsubject={Revised Manuscript},
            plainpages=false
            ]{hyperref}

\newcommand{\xfnm}[1][]{\ifx!#1!\else\unskip,\space#1\fi} 
\usepackage[bottom,multiple]{footmisc} 

\newcounter{algsubstate}


\newcommand{\subparagraph}{}

\usepackage[compact]{titlesec}

\usepackage[nameinlink]{cleveref}

\newif\ifLongVersion
\LongVersionfalse

\ifArxivVersion

	\NotArxivVersionfalse
\else
	
	\renewcommand{\baselinestretch}{1}

	\NotArxivVersiontrue
\fi

\crefname{section}{Sec.}{§§}
\Crefname{section}{Sec.}{§§}
\crefname{subsection}{Sec.}{§§}
\crefname{equation}{Eq.}{}
\Crefname{equation}{Equation}{Equations}
\crefname{appsec}{Appendix}{Appendices}

\newtheoremstyle{theoremdd}
{5pt}
{5pt}
{\itshape}
{0pt}
{\bfseries}
{.}
{ }
{\thmname{#1}\thmnumber{ #2}\thmnote{ (#3)}}

\theoremstyle{theoremdd}

\newtheorem{lemma}{Lemma}

\newtheorem{proposition}{Proposition}
\newtheorem{definition}{Definition}

\newtheorem*{formulation*}{}

\theoremstyle{remark}

\crefname{table}{Table}{Tables}
\crefname{figure}{Fig.}{Figures}
\crefname{theorem}{Theorem}{Theorems}
\crefname{proposition}{Prop.}{Propositions}
\crefname{lemma}{Lemma}{Lemmas}
\crefname{algorithm}{Algorithm}{Algorithms}
\crefname{myclaim}{Claim}{Claims}

\makeatletter
\def\thm@space@setup{%
  \thm@preskip=0.2cm plus 0.2cm minus 0.2cm
  \thm@postskip=\thm@preskip 
}
\renewenvironment{proof}[1][\proofname]{\par
  \pushQED{\qed}%
  \normalfont
  \topsep0pt \partopsep5pt 
  \trivlist
  \item[\hskip\labelsep
        \itshape
    #1\@addpunct{.}]\ignorespaces
}{%
  \popQED\endtrivlist\@endpefalse
  \addvspace{6pt plus 6pt} 
}
\makeatother

\makeatletter

\newtheoremstyle{named}{}{}{\itshape}{}{\bfseries}{.}{.5em}{\thmnote{#3's }#1}
\theoremstyle{named}

\newtheoremstyle{mynamed}{}{}{\itshape}{}{\bfseries}{.}{.5em}{#1 \thmnote{A }}
\theoremstyle{named}

\newcommand{\abs}[1]{\left\lvert{#1}\right\rvert}

\newcommand{\normSquared}[1]{\left\lVert{#1}\right\rVert_2}

\allowdisplaybreaks

\counterwithout{equation}{section} 
\counterwithout{figure}{section} 

\algrenewcommand{\algorithmiccomment}[1]{\hskip3em$\slash\slash$ #1}

\newcommand{\linfinityNorm}[1]{\left\lVert#1\right\rVert_\infty}

\newcommand{\pathFromRoot}[1]{\mathcal{P}_{#1}}
\newcommand{\subtree}[1]{\mathcal{N}_{#1}}
\newcommand{\subtreeEdges}[1]{\mathcal{W}_{#1}}

\newcommand{\isep}{\mathrel{{.}\,{.}}\nobreak}

\def \emult {\small\odot}
\def \unity {\mathbf{1}}
\def \zero {\mathbf{0}}

\newcommand{\N}[1][]{\mathcal{N}_{#1}}
\newcommand{\E}[1][]{\mathcal{E}_{#1}}
\newcommand{\G}[1][]{\mathcal{G}_{#1}}

\def \P {\mathcal{P}}
\def \O {\mathcal{O}}

\def \arcm {{\mathrm{k}}}

\def \R {\mathbb{R}}
\def \Z {\mathbb{Z}}

\def \x {{x}}

\def \j {\mathbf{j}}

\def \state {\eta}
\def \opr {0}
\def \con {1}
\def \attack {a}
\def \optimalAttack {oa}
\def \defend {d}
\def \hardMin {hmin}
\def \hardMax {hmax}

\def \NN {{\mathrm{N}}}

\def \C {{\text{W}}}

\def \lineloss {\text{LL}}
\def \load {\text{LC}}
\def \shed {\text{LS}}

\def \Clineloss {\C^{\lineloss}}
\def \Cshed {\C^{\shed}}
\def \Cload {\C^{\load}}

\def \Csupply {\C^{\text{AC}}}

\def \Clovr {\C^{\slovr}}

\def \setConfigurations {\mathcal{K}}

\def \intermediateState {\text{\fontsize{7}{5}\selectfont in}}

\def \noResponse {\text{\fontsize{5}{5}\selectfont AD}}

\def \maxMin {\text{\fontsize{6}{5}\selectfont Mm}}

\def \noresponseSmall {{\text{\fontsize{7}{5}\selectfont nr}}}

\def \Second {\mathcal{D}_\arcm}

\def \ploss {\mathcal{P}}
\def \loss {\mathcal{L}}
\def \losshat {\widehat{\mathcal{L}}}

\def \lossMaxmin {\mathcal{L}_{\maxMin}}

\def \lossMaxminhat {\hat{\mathcal{L}}_{\maxMin}}

\newcommand{\lossNoResponse}{\loss_{\noResponse}}

\def \ltarget {\loss_\text{target}}
\def \resilience {\mathcal{R}}
\def \resiliencehat {\widehat{\mathcal{R}}}
\def \resilienceTarget {\resilience_\text{target}}
\def \resilienceMaxmin {\resilience_{\maxMin}}
\def \resilienceNoResponse {\resilience_{\noResponse}}
\def \resilienceMaxminhat {\resiliencehat_{\maxMin}}

\def \lcompleteShed {\loss_\mathrm{max}}

\def \pg { pg}
\def \qg { qg}

\def \y { y}

\def \i  {I}

\def \slovr {\text{VR}}

\def \xx {X}

\def \l {\mathcal{L}}

\def \a {\delta}

\def \Xnpf {\mathcal{X}}

\def \ref {\text{ref}}

\def \second {{d}}
\def \dualCoefficient {C}

\def \third {u}

\newcommand{\Third}[1][]{{ \mathcal{U}^{#1}}}

\def \crossTermLineSet {\mathcal{H}}
\def \crossTermNodeSet {\mathcal{I}}
\def \setConsumptionGenerationTuple {\mathcal{J}}
\def \consumptionSet {\mathcal{M}}

\def \x {x}

\def \masterProblem {\text{A-MILP}}
\def \mcp {\text{MCP}}
\def \slaveProblem {\text{O-MISOCP}}

\newcommand{\resistance}[1]{\mathbf{r}_{#1}} 
\newcommand{\reactance}[1]{\mathbf{x}_{#1}} 
\newcommand{\impedance}[1]{\mathbf{z}_{#1}} 

\newcommand{\commonResistance}[1]{\mathbf{R}_{#1}} 
\newcommand{\commonReactance}[1]{\mathbf{X}_{#1}}

\def \f {\omega}

\def \ssum {\textstyle\sum}

\def \mmax {\textstyle\max}
\def \mmin {\textstyle\min}

\usepackage{chngcntr}
\counterwithout{subsubsection}{subsection}

\titleformat{\paragraph}[runin]
{\normalfont\normalsize\itshape}{\thesubsubsection)}{2pt}{}[:]

\def \problemMaxminHat {$(\widehat{\text{Mm}})$}
\newcommand{\tcbtext}[1]{{\color{black}#1}\color{black}}

\newcommand{\tcr}[1]{\color{Mahogany} #1  \color{black}}

\usepackage{tcolorbox}

\usepackage{mdframed}

\DeclareDocumentCommand{\mycommand}{ O{mydefault} m o o o }{%
	p:#2%
	\IfNoValueTF{#3}%
	{}%
	{p:#3}%
	\IfNoValueTF{#4}%
	{}%
	{\@ifmtarg{#4}{}{ p:#4}}%
	\IfNoValueTF{#5}%
	{}%
	{\@ifmtarg{#5}{}{ p:#5}}
	p:#1
}

\newcommand{\test}[3][o o o]{%
	\ifthenelse{\equal{#1}{}}{omitted}{given}%
}

\def \mincardinalityProblem {\text{{\emph{Min-cardinality disruption}}}}
\def \maxlossProblem {\text{{\emph{Budget-k-max-loss}}}}

\def \pre {o}
\def \post {c}

\def \edge {ij}
\def \period {t}

\def \scenarioIdx {s}
\def \node {i}
\def \ts {\period}
\def \nperiod {\mathrm{T}}

\newcommand{\myc}[3]{
	\def \firstString {{#1}}	
	\IfEqCase{#2}{
		{}{\firstString_{#3}}
		{\opr}{\firstString_{#3}^{\opr}}
		{\con}{\firstString_{#3}^{\con}}
		{\star}{\firstString_{#3}^{\star}}
		{\ts}{\firstString_{#3}^{\ts}}
		{\scenarioIdx}{\firstString_{#3}^{\scenarioIdx}}
		{\ts-1}{\firstString_{#3}^{\ts-1}}
		{0}{\firstString_{#3}^{0}}
		{1}{\firstString_{#3}^{1}}
		{2}{\firstString_{#3}^{2}}
		{0\star}{\firstString_{#3}^{0\star}}
		{1\star}{\firstString_{#3}^{1\star}}
		{u\ts}{\widecheck{\firstString}_{#3}^{\ts}}
		{u\ts-1}{\widecheck{\firstString}_{#3}^{\ts-1}}
		{u0}{\widecheck{\firstString}_{#3}^{0}}
		{u1}{\widecheck{\firstString}_{#3}^{1}}
		{u2}{\widecheck{\firstString}_{#3}^{2}}
		{u\scenarioIdx}{\firstString_{#3}^{\scenarioIdx}}
		{\ts-1}{\firstString_{#3}^{\ts-1}}
		{\period}{\firstString_{#3}^{\period}}
		{\period-1}{\firstString_{#3}^{\period-1}}
		{\nperiod}{\firstString_{#3}^{\nperiod}}
		{n}{\firstString_{#3}}
		{l}{\widehat{\firstString}_{#3}}
		{u}{\widecheck{\firstString}_{#3}}
		{nr}{\firstString^{\noresponseSmall}_{#3}}
		{in}{\firstString^{\intermediateState}_{#3}}
		{in\star}{\firstString^{\intermediateState\star}_{#3}}
		{max}{\mathbf{\overline{\firstString}}_{#3}}
		{min}{\mathbf{\underline{\firstString}}_{#3}}
		{\hardMax}{\bm{\overline{\overline{\firstString}}_{#3}}}
		{\hardMin}{\bm{\underline{\underline{\firstString}}_{#3}}}	
		{constant}{\mathbf{\firstString}_{#3}}
		{\pre}{\firstString_{#3}^{\pre}}
		{\post}{{\firstString}_{#3}^{c}}
		{act}{\firstString_{#3}^{act}}
		{set}{\firstString_{#3}^{set}}
		{ref}{\mathbf{\firstString}_{#3}^{\text{ref}}}
		{r}{\firstString_{#3}^{r}}
		{nom}{\mathbf{\bm{\firstString}_{#3}^{nom}}}
		{stab}{\bm{\firstString_{#3}^{stab}}}
		{devmax}{\firstString_{#3}^{dev,max}}
		{reg}{\bm{\firstString_{#3}^{reg}}}
		{ev}{\firstString_{#3}^{ev}}
		{nev}{\firstString_{#3}^{nev}}
		{maxev}{\bm{\overline{\firstString}_{#3}^{ev}}}
		{maxnev}{\bm{\overline{\firstString}_{#3}^{nev}}}
		{\state}{\firstString_{#3}^{\state}}
		{\attack}{\firstString_{#3}^{\attack}}
		{\optimalAttack}{\firstString_{#3}^{\attack\star}}
		{\defend}{\firstString_{#3}^{\defend}}
	}
}

\newcommand{\kcc}[2]{\myc{kc}{#1}{#2}}	
\newcommand{\kgc}[2]{\myc{kg}{#1}{#2}}

\newcommand{\configurationc}[2]{\myc{\kappa}{#1}{#2}}	

\newcommand{\Lossc}[1]{\myc{L}{#1}{}}

\newcommand{\Costc}[2]{\myc{\mathcal{C}}{#1}{#2}}

\newcommand{\costMaxmin}[1]{\Costc{#1}{\maxMin}}

\newcommand{\ellc}[2]{\myc{\ell}{#1}{#2}}
\newcommand{\etac}[2]{\myc{\eta}{#1}{#2}}
\newcommand{\nuc}[2]{\myc{\mathrm{v}}{#1}{#2}}

\newcommand{\Pc}[2]{\myc{P}{#1}{#2}}
\newcommand{\Qc}[2]{\myc{Q}{#1}{#2}}

\newcommand{\nucc}[2]{\myc{vc}{#1}{#2}}
\newcommand{\nugc}[2]{\myc{vg}{#1}{#2}}

\newcommand{\pcc}[2]{\myc{pc}{#1}{#2}}
\newcommand{\qcc}[2]{\myc{qc}{#1}{#2}}

\newcommand{\ptc}[2]{\myc{p}{#1}{#2}}
\newcommand{\qtc}[2]{\myc{q}{#1}{#2}}

\newcommand{\pgc}[2]{\myc{pg}{#1}{#2}}
\newcommand{\qgc}[2]{\myc{qg}{#1}{#2}}

\newcommand{\vdc}[2]{\myc{{\Delta \mathrm{v}}}{#1}{0}} 
 
\newcommand{\xc}[2]{\myc{\x}{#1}{#2}}
\newcommand{\Xc}[2]{\myc{\Xnpf}{#1}{#2}}

\newcommand{\Zc}[2]{\myc{\mathcal{Z}}{#1}{#2}}

\newcommand{\lcc}[2]{\myc{\beta}{#1}{#2}}

\newcommand{\Flowc}[2]{\myc{\mathcal{F}}{#1}{#2}}
\newcommand{\Voltc}[2]{\myc{\mathcal{V}}{#1}{#2}}

\newcommand{\transpose}[1]{{#1}^{\top}}

\newcommand{\mycc}[4]{
	\def \firstString {{#1}}	
	\IfEqCase{#2}{
		{n}{\firstString_{#3}^{#4}}
		{max}{\overline{\firstString}_{#3}^{#4}}
		{min}{\underline{\firstString}_{#3}^{#4}}
	}
}

\def \nperm {Z}
\def \YY {Y}
\def \VV {V}

\def \startParameter {\mathrm{m}}
\def \seqEnd {\mathrm{e}^j}
\def \seqStart {\mathrm{s}^j}
\def \permutation {\sigma^j}

\title{\huge Evaluating Resilience of Electricity Distribution Networks via A Modification of  Generalized Benders Decomposition Method}
\author{Devendra Shelar, Saurabh Amin, and Ian Hiskens
	\thanks{Manuscript resubmitted on November 30, 2020. This work was supported by awards: AFOSR ``Building attack resilience into complex networks", NSF CAREER (CNS-1453126), and \enquote{Modeling \& Analysis of Load Ensembles} (ECCS-1810144).}
	\thanks{D. Shelar and S. Amin are with the Laboratory for Information and Decision Systems, Massachusetts Institute of Technology (MIT), Cambridge, MA 02139 USA (e-mail: \{shelard, amins\}@mit.edu, phone: 857-253-8964).}
	\thanks{I. A. Hiskens is with the Department of Electrical Engineering and Computer	Science, University of Michigan, Ann Arbor, MI 48109 USA (e-mail: hiskens@umich.edu).}
}

\ifArxivVersion
\renewcommand{\baselinestretch}{1}
\fi 
\ifOneColumn
\renewcommand{\baselinestretch}{2}
\else

\iftcnsVersion
\usepackage[margin=0.75in]{geometry}

\ifArxivVersion
\renewcommand{\baselinestretch}{1}
\begin{document}
	
	\maketitle

	\begin{abstract}

\tcbtext{This paper presents a computational approach to evaluate the resilience of electricity Distribution Networks (DNs) to cyber-physical failures. In our model, we consider an attacker who targets multiple DN components to maximize the loss of the DN operator. We consider two types of operator response: (i) Coordinated emergency response; (ii) Uncoordinated autonomous disconnects, which may lead to cascading failures. To evaluate resilience under response (i), we solve a Bilevel Mixed-Integer Second-Order Cone Program which is computationally challenging due to mixed-integer variables in the inner problem and non-convex constraints. Our solution approach is based on the Generalized Benders Decomposition method, which achieves a reasonable tradeoff between computational time and solution accuracy. Our approach involves modifying the Benders cut based on structural insights on power flow over radial DNs. We evaluate DN resilience under response (ii) by sequentially computing autonomous component disconnects due to operating bound violations resulting from the initial attack and the potential cascading failures. Our approach helps estimate the gain in resilience under response (i), relative to (ii).} 
		
	\end{abstract}
	
	
	\section{Introduction}\label{sec:introduction} 
	~
	
	Despite the ongoing modernization of electricity distribution networks (DNs), many distribution system operators face both strategic and operational challenges in ensuring a reliable and secure service to their customers. The integration of distributed generators (DGs) and new monitoring and control capabilities has enabled flexible operations, which can be utilized to respond to typical failure events such as sudden voltage drop~\cite{henrikSandberg0,dorfler,shelarAminHiskens}. However, these capabilities also expose the vulnerabilities of DNs to adversaries~\cite{lalitaShankar,bulloDetectionAndIdentification}. Particularly, cyber-physical failures in  DNs can result in contingencies that cause cascading network outages~\cite{abhinavVerma,salmeron,karlJohansson}. In this article, we argue that the flexibility of modern DNs can be leveraged to generate a timely and effective response to cyber-physical failures. We present a computational approach for evaluating the DN resilience under realistic response capabilities.  
	
	For a given cyber-physical failure model and an operational response capability of the operator, our optimization-based approach can be used to evaluate the \emph{worst-case} post-contingency loss due to various factors, such as the impact of voltage degradation and costs of load control, load shedding, and line losses. By evaluating this loss for different response (or control) operations, we can compute their relative \emph{value} in maintaining the DN resilience against the given class of failures. From a strategic viewpoint, this computation is useful for knowing which control capabilities, if deployed in the DN, will be most effective in response to contingencies arising from such failures. From an operational viewpoint, it can help the operator to implement the response in a timely manner to limit the uncontrolled outages resulting from the triggering of automatic protection mechanisms. 
		
	Indeed, defining the appropriate operational response capability, the cyber-physical failure model, and the nature of attacker-operator interaction are all crucial aspects of our problem, which we introduce next. 
	
\begin{figure}[htbp!]\setstretch{1}
		\centering
		\tikzset{every node/.append ={font size=tiny}} 
		\tikzstyle{dernode}=[circle, fill=blue!70, minimum size=0.5pt, inner sep = 4]
		\tikzstyle{load}=[draw,rectangle, minimum height=4pt]
		\def \drawgrid {\draw[step=1,gray, ultra thin, draw opacity = 0.5] (0,0) grid (3,4);}
		\def \drawSubstation {\draw[-, line width = 2pt] (0.8,4.05) -- (2.2,4.05)  node [midway,above] {};}
		\def \drawZero {\node (0) at (1.5,3.85) {};}
		\def \drawBus at (#1,#2,#3); {\draw [line width=2.5pt,-] (#1,#2-0.8) 
			-- (#1,#2) node (#3) {}
			-- (#1,#2+0.8);}
		\def \drawBuss at (#1); {\draw [line width=2.5pt,-] ([yshift=-0.7]#1) -- ([yshift=0.7]#1);}
		\def \drawLines at (#1,#2,#3,#4); {\draw [line width=1.5pt,-] (#1,#2) -- (#3,#4);}
		\def \connectBuses (#1,#2); {\draw [line width=1.5pt,-] (#1.center) -- (#2.center);}
		\def \drawArrows at (#1,#2,#3,#4); {\draw [red,line width=1.5pt,->, >=latex] (#1,#2) -- (#3,#4);}
		\def \drawPointers (#1,#2); {\draw [red,line width=1.5pt,->, >=latex] (#1) -- (#2);}
		\def \drawDG at (#1,#2,#3); {\node[circle,fill=blue,minimum size=5pt] (#3) at (#1,#2) {}; }
		\def \connectComponent (#1,#2); {\path[draw] let \p1=(#1), \p2=(#2) in (\x1,\y2) -- (#2.west); }
		\tikzstyle{alterConn} = [draw,line width=1pt,-,dashed]
		\def \switch (#1,#2); {
			\path (#1) ++(40:1) node (#1a) {};
			\path (#1) ++(0:0.8) node (#1b) {};
			\path (#1) ++(40:0.8) node (#1c) {};
			\draw[#2] (#1.center) -- (#1a.center); 
			\draw[->, line width=0.5pt,>=latex,#2] (#1b.center) to [out=90,in=-60] (#1c.center);
		}
		\scalefont{0.7}
		\resizebox{0.9\textwidth}{!}{
		\begin{tikzpicture}[scale=0.6]
		
		
		\node (substation) at (0.5,0) {};
		
		\drawBus at (0.5,0,tnbus); \drawBus at (5,0,ssbus); \drawBus at (8.5,0,jbus); \drawBus at (11,1.4,kbus); \drawBus at (11,-1.4,lbus); 

		\connectBuses (ssbus,jbus); \connectBuses (kbus,jbus); \connectBuses (lbus,jbus); 
		
		\draw [] (tnbus.center) -- (ssbus.center);
		
		
		\node[align=left,above right=0.65 and -0.55 of kbus,color=blue](lcd) {Load control, \\ Disconnections};
		\node[align=left,left= -0.1  of lcd,color=blue](){\normalsize (c)};
		
		\node[align=left,below right=0.5 and -0.75 of lbus](ald) {- Local/Autonomous\\ \hspace{0.15cm} Disconnects};
		\node[align=left,left= -0.1  of ald](){\normalsize (b)};
		
		\node[align=right, above=2 of tnbus,xshift=0.2cm,red] (attacker) {Targeted\\attack
		}; 	
		\node[align=left, above=0.5 of tnbus] (tnname) {TN}; 
		\node[draw,align=left, above right=0.0 and 0 of ssbus] (SA) {SA}; 
		\node[align=left, below right=-0.1 and 0 of ssbus] (sbname) {substation}; 
		\node[draw,align=left, above right=-0.4 and 0.3 of tnname] (dgms) {Control center\\--Dispatch\\--fault/outage \\\hspace{0.08cm} management
		};
		\draw[->,>=latex,red] (attacker.east) to [out=0,in=90]  (dgms.north);
		\draw[->,>=latex,red,dashed] (dgms.east) to [out=0,in=180]  (11.2,2.4);
		\draw[->,>=latex,blue,dashed] (SA.east) to [out=0,in=180]  (11.4,1.6);
		
		\node[align=left, below=0.5 of ssbus,red] (vd) {$-\vdc{}{0}$}; 
		\node[align=center, below=0 of vd,red] (tsd) {TN-side disturbance}; 
		
		\node[align=left,above right= -0.5 and 0.5 of dgms,red]  {DN-side\\disturbance}; 
		\node[align=left,above right= -0.8 and -0.05  of dgms](){\normalsize(a)};
		\node[align=left,below right= -0.3 and 0  of SA](){\normalsize (c)};
		\node[align=left,above right= -0 and 0  of SA,blue]  {emergency\\control}; 
		
		\node[dernode,below right= 0.15 and 0.75 of kbus] (kder1) {};
		\node[dernode,above right= 0.15 and 0.75 of kbus,black] (kder2) {};
		\node[align=left,right= 0  of kder2](){DG};
		\node[dernode,below right= 0 and 0.75 of lbus,black] (lder) {};
		
		\node[load, right=0.75 of kbus, minimum width=5pt,fill=blue](kload1) {};
		\node[align=left,right= 0.2  of kload1](){Load};
		\node[load, right=0of kload1, minimum width=5pt](kload2) {};
		\node[load, above right=0pt and 0.75 of lbus, minimum width=10pt,fill=blue,black](lload) {};
		
		\node[align=left,right= 0  of lload, black](){Load};
		\node[align=left,right= 0  of lder,black](){DG};
		
		\connectComponent(kbus,kder1); \connectComponent(lbus,lder); 	\connectComponent(kbus,kder2); 
		\connectComponent(kbus,kload1);	\connectComponent(lbus,lload);
		
		\node[left= 0.65 of kload1](kloadsw) {}; \switch(kloadsw,blue);	
		\node[left= 0.60 of kder1](kder1sw) {}; \switch(kder1sw,blue);	
		\node[left= 0.60 of kder2](kder2sw) {}; \switch(kder2sw,red);	
		\node[left= 0.60 of lder](ldersw) {}; \switch(ldersw,black);	
		\node[left= 0.60 of lload](lloadsw) {}; \switch(lloadsw,black);
		\end{tikzpicture}
	}\caption[]{An illustration of modern DN operations: (b) and (c) are secure  (\textcolor{blue}{blue}), and (a) is compromised  (\textcolor{red}{red}).}\label{fig:attackerDefenderInteractions}
	\end{figure}
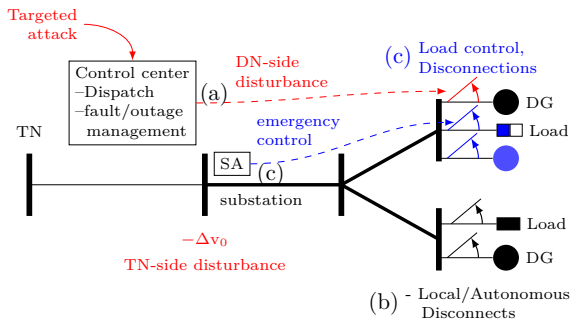

	\tcbtext{Firstly, we consider three different operations supported by modern DNs; see \cref{fig:attackerDefenderInteractions}. Operation (a) refers to remote control by the control center; operation (b) - autonomous disconnects of components due to activation of local protection systems; and operation (c) -  emergency control by the Substation Automation (SA) systems. Typically, operation (a) may include coordination of one or more DNs, dispatch of generators, fault/outage management, etc. Furthermore, operation (a) is typically exercised during normal operating conditions over relatively longer time scales (every 15 minutes or more) and more regularly than operation (b) or (c)~\cite{controlCenter}. 
	}
		
	\tcbtext{In contrast, operation (b) or (c) are executed in emergency situations when certain components of a DN experience operating bound violations. In our model, operation (b) is an \emph{uncoordinated} tripping of components based on local checks of operating bounds at the DN nodes. However, operation (c) is a  \emph{coordinated} action comprising of  DG dispatch, load control, and preemptive tripping of components. Operation (c) utilizes information from DN meters that includes node-level consumption, distributed generation, and nodal voltages. When either operation (b) or  (c) is executed, it happens at a faster timescale (few seconds). We assume that the operation (c) subsumes (b) by making decisions that are coordinated across the SA system. Hence, (b) and (c) are never simultaneously active.} 
	
	\tcbtext{
	Secondly, we study an attack model in which operation (a) (see~\cref{fig:attackerDefenderInteractions}) is compromised by an attacker~\cite{nescor}. Control center operations are prone to cyber attacks by remote entities as indicated by recent real-world incidents~\cite{ukraine}. 
	In fact, the impact of such security failures may be aggravated due to a failure in the adjoining transmission network (TN). We model the impact of such TN-side failure as a voltage sag (drop in the substation voltage), and that of the security attack on operation (a) as disturbances resulting due to the tripping of DGs, at DN nodes. Thus, an important feature of our attack model is that it can capture the effect of contingencies resulting from simultaneous TN and DN failures.} 

\tcbtext{ 	
	The operator can respond to abovementioned cyber-physical failures by operation (b) or (c). We refer to these operations as \emph{response} (b) and \emph{response} (c) in order to clearly distinguish them from operation (a) which is prone to attack. We argue that both response (b) and (c) can be considered secure against remote cyberattacks. Response (b) relies on local checks of operating parameters, and is implemented by actuators which are geographically distributed. The SA systems (response (c)) were recently required to secure both physical and cyber-security perimeters by NERC~\cite{nerccip}.
Therefore, even response (c) can be considered secure from a remote attack. }

	\tcbtext{Thirdly, we develop a computationally tractable approach to address the problem of determining the \emph{worst-case} post-contingency loss when the operator optimally  responds to the attacker actions with response (b) or (c)}. For the case of response (c), we formulate a bilevel mixed-integer second-order cone program (BiMISOCP), which captures the  sequential nature of attacker-operator interaction (\cref{sec:modeling1}). The  inner (operator) problem consists of mixed-binary variables which model response (c), and the second-order cone constraints model the non-linear power flows (NPF) over a radial DN. The operator's objective is to minimize the post-contingency loss. The outer (attacker) problem is to determine an attack that will maximize the operator's loss, assuming the operator responds optimally. The worst case post-contingency loss for response (c) is given by the maximin value of the BiMISOCP. 
	
	To compute the worst-case post-contingency loss under response (b), we present a two-step approach~(\cref{sec:autonomousDisconnects}). For a given attack, the first step evaluates the impact of cascading failures in the DN by determining DG disconnections. The second step determines the effect of DG disconnects on the nodal voltages and evaluates its impact on the load control/shedding. We propose a randomized algorithm to estimate the worst-case post-contingency loss under a maximally disruptive attacker strategy.
	
	Several papers have used bilevel formulations for security assessment of power systems~\cite{karlJohansson,abhinavVerma,bulloDetectionAndIdentification,salmeron,shelarAminTCNS,shelarAminHiskens}. 
	 However, BiMISOCP formulations with mixed-integer variables in the inner problem are extremely challenging to solve, and have received limited attention in the literature. Even under linear constraints, the resulting bilevel mixed-integer linear program (BiMILP), is still hard to solve due to integer variables in the inner problem. Previous works have utilized a 
	 relaxation technique to reformulate this BiMILP as a single-level MILP, which can be solved using an advanced branch-and-bound algorithm~\cite{bard1990,xuWangBnBBiMIP}).  Recent papers have proposed intersection cuts~\cite{fischetti2016new,fischetti2018new}, and disjunction cuts~\cite{disjuctionCutsWaterMelon,disjunctionCutsValueFunction} to introduce stronger cuts. However, these approaches only solve a weak relaxation of the original BiMILP~\cite{kevinwood,baldickBowenHua}. Other methods for solving BiMILPs include the Generalized Benders Decomposition (GBD) method~\cite{baldickBowenHua} and column constraint generation~\cite{zeng2014solving}. However, the presence of integer variables in the inner problem with nonlinear constraints precludes the application of these methods to solve our BiMISOCP problem.  
	
	\tcbtext{
	We address this challenge by making the following contributions: 
	\begin{enumerate}[leftmargin=*]
		\item We derive structural results on non-linear power flows on a radial DN topology~(\cref{sec:technicalResults}), and use these insights to generate more effective cuts for the inner (operator) problem. 
		
		\item We describe a computational approach for solving a BiMISOCP based on its reformulation as an equivalent min-cardinality disruption problem, and then applying a new algorithm that can be viewed as a \emph{modified} GBD method. Typically, in each iteration of the classical GBD method, a generalized Benders cut is added to the master (attacker) subproblem, in which a linear expression in the attack variables is constrained to be greater than zero~\cite{generalizedBenders}. We modify this method by changing the right hand side of the cut to be a positive value $\epsilon$. Also, we introduce a heuristic to select appropriate values of $\epsilon$ for each generalized Benders cut by exploiting the structure of our problem (\cref{sec:solutionApproach}).
	
	\item Additionally, we provide a novel two-step approach to estimate the worst-case post-contingency loss under autonomous disconnects. The difference between this loss and the operator's loss computed by the abovementioned GBD method provides an estimate of the value of timely operator response. 
	\end{enumerate} }

 \tcbtext{The computational results demonstrate that our modified GBD method achieves a good tradeoff between the  computational speed and gap from the optimal (maximin) attack}~(\cref{sec:computationalStudies}). 
	
\tcbtext{	We refer the reader to our online technical report~\cite{technicalReport}, which includes further details on the following:
(i) Further justification of our attacker and operator model, and its technological feasibility (some discussion on extensions to other models is also provided); (ii) Full expressions of the constraints of the attacker and operator subproblems, and the GBD cuts; and (iii) Details of some technical results in \cref{sec:technicalResults}. }

\section{Modeling and Problem formulation}
\label{sec:modeling1}
In this section, we describe our approach to evaluate DN resilience, and then present our sequential attacker-operator bilevel problem formulation.
\subsection{Evaluating resilience of DNs}\label{subsec:resilienceFramework}
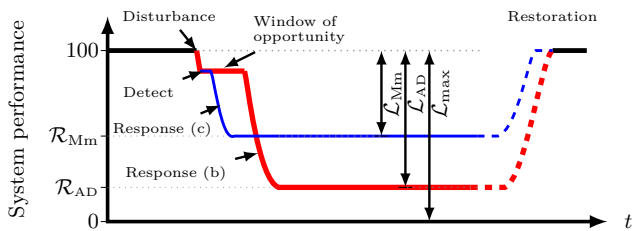
\begin{figure}[htbp!]
	\tikzstyle{mfs} = [font=\fontsize{9}{8}\selectfont,align=left]
	\tikzstyle{mrfs} = [mfs,rotate=90]
	\tikzstyle{fs} = [font=\fontsize{6}{6}\selectfont,align=left]
	\tikzstyle{ft} = [draw,line width=2pt]
	\tikzstyle{ftr} = [draw,line width=2.2pt,red]
	\tikzstyle{ftb} = [draw,line width=1.2pt, blue]
	\tikzstyle{ftg} = [draw,line width=1.2pt, green]
	\tikzstyle{ftrd} = [draw,line width=2.2pt,red,dashed]
	\tikzstyle{ftbd} = [draw,line width=1.2pt, blue,dashed]
	\tikzstyle{ptr} = [-latex,line width=1pt]
	\tikzstyle{mk} = [gray,line width=0.6pt, dotted, opacity=0.8]
	\def \maxp {2.5}
	\def \dip {1}
	\def \rebound {2}
	\resizebox{8.5cm}{!}{\begin{tikzpicture}
			
			\node[](ytip) at (-0.5,3) {}; \node[](xtip) at (6.8,0) {};
			\draw[<->,>=latex,line width=1.5pt] (ytip.center) |- (xtip.center); 
			\node[mfs,right=0 of xtip]() {$t$}; 
			
			\draw[ft] (-0.5,2.5) -- (0.8,2.5) (6,2.5) -- (6.5,2.5);
			\draw[ftr] (0.8,2.5) -- (0.85,2.2) -- (1.5,2.2);
			\path[ftr] (1.5,2.2) .. controls  (1.55,2) and (1.7,0.6) .. (2,0.5) -- (4.8,0.5);
			\path[ftrd] (4.8,0.5) -- (5.3,0.5) .. controls (5.6,0.6)  and (5.85,2.45) .. (6.0,2.5);
			\path[ftb] (0.85,2.2) -- (1.00,2.2) .. controls (1.05,2.2) and (1.15,1.15) .. (1.35,1.25) -- (4.9,1.25);
			\path[ftbd] (4.9,1.25) -- (5.2,1.25) .. controls (5.5,1.35) and (5.65,2.45) .. (5.75,2.5) --(6.0,2.5);
			\draw [mk] (0.8,2.5) -- (5.0,2.5);
			
			\node[fs] (dn) at (0.4,3) {Disturbance}; 
			\node[fs] (re) at (6.0,3) {Restoration}; 
			\node[fs] (den) at (0.1,1.9) {Detect}; 
			\node[fs] (trn) at (0.3,1.35) {Response (c)}; 
			\node[fs] (woo) at (2.3,2.8) {Window of\\opportunity}; 
			\node[fs] (cas) at (0.5,0.7) {Response (b)};
			\draw[ptr] (dn) -- (0.8,2.5);
			\draw[ptr] 	(trn) -- (1.15,1.7);
			\draw[ptr] 	(den) -- (0.9,2.2);
			\draw[ptr] 	(cas) -- (1.7,1);
			\draw[ptr] 	(woo) -- (1.2,2.25);
			\def \a {-1.2}; \def \b {-0.5}; \def \c {2.5}; \def \d {1.25}; \def \e {0.5}; \def \f {0}; \def \g {-0.8}; \def \h {-1.3}; \def \i {-1.8};
			\def \yredlow {0.5}; 
			\def \ybluelow {1.25}; 
			\def \ygreenlow {1.9};
			\def \axmicro {3.4}; 
			\def \axmaxmin {2}; 
			\def \axnoresponse {4.9};
			\def \axxmicro {4.1}; 
			\def \axxmaxmin {3.5}; 
			\def \axxnoresponse {3.85};  
			\def \axxmax {4.2}; 
			\def \ymax {2.5};
			\def \start {-0.5};
			\foreach \y/\l in {0/0,\yredlow/\resilienceNoResponse, \ybluelow/\resilienceMaxmin,  \ymax/100}{
				\draw (-0.6,\y) -- (-0.5,\y);
				\node[mfs,anchor=east] () at (-0.5,\y) {$\l$}; 
			}
			\foreach \x/\y/\z in {\start/\yredlow/2,\start/\ybluelow/1.2}{
				\draw[gray,opacity=0.8,dotted] (\x,\y) -- (\z,\y);
			}
			
			\node[mrfs] (ylabel) at (-1.8,1.5) {System performance};
			
			\tikzstyle{ptre} = [ptr,<->,>=latex]

			\def \w {0.2}
			\foreach \x/\y/\l in {\axxmaxmin/\ybluelow/\lossMaxmin, \axxnoresponse/\yredlow/\lossNoResponse,\axxmax/0/\lcompleteShed} {
				\pgfmathsetmacro\xx{\x - \w/2}
				\pgfmathsetmacro\xxx{\x + \w/2}
				\draw [ptre] (\x,\ymax) -- (\x,\y);
				\draw [] (\xx,\y) -- (\xxx,\y);
				
				\pgfmathsetmacro\lx{\x + 0.165}
				\pgfmathsetmacro\ly{(\ybluelow + \ymax)/2}
				\node[mrfs] (lnr) at (\lx,\ly)  {$\l$};
			}
			\foreach \xx/\y/\xxx in { \axmaxmin/\ybluelow/\axxmaxmin, \axnoresponse/\yredlow/\axxnoresponse} {
				\draw [mk] (\xx,\y) -- (\xxx,\y);
			}
			
	\end{tikzpicture}}
	\caption{System performance under various response capabilities. (The dashed lines indicate the restoration aspect of DN resilience, which is not  the focus of this paper.)} \label{fig:resilienceDefinition}
\end{figure}

A system's resilience is broadly defined as \enquote{its ability to prepare and plan for, absorb, recover from, and more successfully adapt to adverse events}~\cite{resilienceDefnNIAC}. To systematically evaluate a DN's resilience, we need to select both a class of adverse events and the DN's ability to respond to those events. In our setup, the class of disruptions is  denoted by $\Second$, where $\arcm$ is attacker's resource constraint. \tcbtext{We also consider a set of feasible operator strategies (denoted by $\Third$) which model the response (c).}\  We denote by $\lossMaxmin$ the post-contingency loss  which is a measure of the maximum reduction in system performance under $\Second$; see \cref{fig:resilienceDefinition}. Let $\lcompleteShed$ denote the loss incurred when all loads and DGs are disconnected. Then,  $\resilienceMaxmin \coloneqq 100(1-\lossMaxmin/\lcompleteShed)$ i.e., the percentage drop in system performance, can be viewed as a  metric of the DN resilience under the response capabilities $\Third$ against attacks in $\Second$.  

Furthermore, to compare the DN resilience under response (c) with the case of autonomous disconnections (b), we need to estimate the maximum loss corresponding to the response (b) that would be induced by an attack in $\Second$. Let the automatic disconnect actions be denoted by $\third_\noresponseSmall$, the resulting network state by $\xc{}{\noresponseSmall}$, and the corresponding loss by $\lossNoResponse = \Lossc{}(\third_\noresponseSmall,\xc{}{\noresponseSmall}$). Then, the resilience metric of the DN under autonomous disconnections (AD) can be written as $\resilienceNoResponse = 100(1-\lossNoResponse/\lcompleteShed)$. 

\Cref{fig:resilienceDefinition}  qualitatively illustrates the evolution of system performance over time. Initially, the  DN is operating under the nominal conditions. Then, due to the TN/DN-side disturbances, the system performance degrades. If the operator fails to respond in a timely manner (in less than a few seconds), then an uncontrolled cascade can occur due to response (b), resulting in a loss $\lossNoResponse$.


By evaluating the post-contingency loss due to a timely operator response, and comparing it with the loss under the autonomous disconnections, we can estimate the \emph{value} of the timely response toward improving the DN's resilience. We assume that all the devices within the DN are networked in response (c). Hence, $\resilienceMaxmin \ge \resilienceNoResponse$, and we can evaluate the relative value of operational response (or equivalently, the improvement in DN resilience) as  $(\resilienceMaxmin - \resilienceNoResponse)$; see~\cref{sec:computationalStudies}.

\subsection{BiMISOCP formulation for $\lossMaxmin$}
We now describe our bilevel program to evaluate $\lossMaxmin$ over a radial DN for specific attacker and operator models. 

We model the DN as  a tree network of node set $\N\bigcup\{0\}$ and line set $\E$. Without loss of generality, we assume that each  node of the DN has a load and a DG.
Furthermore, we consider the constant power model for both loads and DGs. 
We refer the reader to \cref{tab:notationsTable1} in the Appendix for the table of notations. 

\ifOneColumn
\def \www {0.7\textwidth}
\else
\def \www {\textwidth}
\fi
We formulate a bilevel problem to model the  sequential interaction between the strategic attacker (leader) and the operator (follower). First, we model the effect of a TN-side disruption on the DN as a drop in the substation voltage $\vdc{}{0}$. Next, we  consider a specific attack model where the attacker attacks operation (a) by disrupting a subset of DGs in the DN. We denote an attacker-induced failure by $\second\in\{0,1\}^{\N}$, where $\second_i = 1$ indicates that the DG at node $i$ is disrupted, $\second_i = 0$ otherwise. Let $\arcm$ denote the attacker's resource budget. 
For a given TN-side disruption $\vdc{}{0}$ and attacker action $\second$, we consider that the operator can exercise response (c) by  exercising load control, changing the DG output, and disconnecting the loads and DGs, if necessary. We denote the response (c) by $\third=(\lcc{}{},\pgc{}{},\qgc{}{},\kcc{}{},\kgc{}{})$,  
and use $\xc{}{} = (\pcc{}{},\qcc{}{},\ptc{}{},\qtc{}{},\Pc{}{},\Qc{}{},\nuc{}{},\ellc{}{})$ to denote the post-contingency network state. Finally, we denote by  $\Lossc{}\left(\third,\xc{}{}\right)$ the loss function for a given operator response and network state. We state our problem as follows: 
\begin{alignat}{14}
 \nonumber \lossMaxmin &\coloneqq \mmax_{\second\in\{0,1\}^{\N}} \  \costMaxmin{}(\second)&&&&	\hspace{4cm}\\ 
\label{eq:attackerResourceConstraint}	 \text{s.t.} & \qquad \ssum_{i\in\N} \second_i \le \arcm, \\
	\nonumber \costMaxmin{}(\second) & \coloneqq \mmin_{\third,\xc{}{}} \Lossc{}\left(\third,\xc{}{}\right) &\text{s.t.} 
 	\end{alignat}
 	\begin{alignat}{10}
 		\label{eq:postContingencyVoltage1} \nuc{}{0} &= \nuc{nom}{} - \vdc{}{0} \\
	\label{eq:dgConnectivityPostContingency}
    \kgc{}{i} &\ge \second_i  &&  &&\forall\  i\in \N\\
    \label{eq:dgoutputAttackerAction} \pgc{}{i} &\le \pgc{max}{i}(1-\second_i) &&  &&\forall\ i \in \N\\ 
	\label{eq:dgoutputOperatorAction} \pgc{}{i} &\le \pgc{max}{i}(1-\kgc{}{i}) &&  &&\forall\ i \in \N\\
	\label{eq:dgoutputActivePositive} \pgc{}{i} &\ge 0 &&  &&\forall\ i \in \N\\ 
\label{eq:dgoutputReactiveActiveConstraint}  \qgc{}{i} &\ge-\etac{constant}{i}\pgc{}{i}, &\qgc{}{i} &\le \etac{constant}{i}\pgc{}{i}  && \forall \ i\in\N\\
	\label{eq:integralityConstraints}  \kgc{}{i}&\in\{0,1\}, & \kcc{}{i}&\in\{0,1\}  \quad&& \forall\ i\in \N\\
\label{eq:loadControlSheddingConstraint}   \lcc{}{i} &\ge \left(1-\kcc{}{i}\right)\lcc{min}{i}, & \lcc{}{i}&\le \left(1-\kcc{}{i}\right) && \forall \ i\in\N\\
	\label{eq:loadControlParameterConsumptionConstraint}\pcc{}{i} &=  \lcc{}{i}\pcc{max}{i}, & \qcc{}{i} &=  \lcc{}{i}\qcc{max}{i}  &&\forall\ i \in \N\\ 
\label{eq:voltageDisconnectDG}  \kgc{}{i}&\ge \nugc{min}{i} - \nuc{}{i}, & \kgc{}{i}&\ge\nuc{}{i} - \nugc{max}{i}  \quad&& \forall\ i\in \N\\
\label{eq:voltageDisconnectLoad}  \kcc{}{i}&\ge \nucc{min}{i} - \nuc{}{i},  &\kcc{}{i}&\ge \nuc{}{i} - \nucc{max}{i}  && \forall\ i\in \N\\
	\label{eq:totalPowerConsumption} \ptc{}{i} &= \pcc{}{i} -  \pgc{}{i}, &\qtc{}{i} &= \qcc{}{i} -  \qgc{}{i} \quad &&\forall\ i\in\N\\
	\label{eq:conserveRealTrue} \Pc{}{ij} &= \ssum_{k:(j,k) \in\E} \Pc{}{jk} &+\ptc{}{j} &+ \resistance{ij}\ellc{}{ij}  && {\small \forall\ (i,j)\in\E} 	\\
	\label{eq:conserveReacTrue} \Qc{}{ij} &= \ssum_{k:(j,k) \in\E} \Qc{}{jk} &+\qtc{}{j} &+ \reactance{ij}\ellc{}{ij}  && {\small \forall\ (i,j)\in\E} 	\\
	\label{eq:voltageTrue} \nuc{}{j} &= \nuc{}{i}-2\resistance{ij}\Pc{}{ij}-&2\reactance{ij}&\Qc{}{ij} + \abs{\impedance{ij}^2}\ellc{}{} && {\small \forall\ (i,j)\in\E} 	\\
	\label{eq:currentTrue} \ellc{}{ij}\nuc{}{i} &=  \Pc{}{ij}^2 +  \Qc{}{ij}^2 && && {\small \forall\ (i,j)\in\E} 	
	 		\end{alignat}
where $\lossMaxmin$ denotes the Max-min (Mm) post-contingency loss under response (c).
For a fixed attack, the operator's objective is to minimize the post-contingency loss  $\Lossc{}\left(\third,\xc{}{}\right)$. The attacker's objective is to choose an attack that maximizes  the minimum post-contingency loss. 

We define $\Lossc{}(\third,\xc{}{})$ as the sum of: (i) cost due to loss of voltage regulation, (ii) cost of load control, (iii) cost of load shedding, and (iv) cost of line losses:
\begin{equation}\label{eq:costGCregime}\fontsize{9}{8}\selectfont
\begin{split}
\hspace{-0.25cm}\Lossc{}(\third,\xc{}{}) = \Clovr\linfinityNorm{\nuc{nom}{} - \nuc{}{}} 
+  \ssum_{i\in\N}\ \Cload _i\left(\unity-\lcc{}{i}\right)  \pcc{max}{i}\\ + \ssum_{i\in\N}\ \left(\Cshed_i-\Cload_i\right)\kcc{}{i}\pcc{max}{i} + \Clineloss\ssum_{ij\in\E}\resistance{ij}\ellc{}{ij},
\end{split} 
\end{equation}
where 
for load at node $i$, $\Cload_i \in \R_{+}$ denotes the cost of per unit load controlled, $\Cshed_i \in \R_{+}$ and $\Cshed_i \ge \Cload_i$ is the cost in dollars of per unit load shed;  $\Clineloss\in \R_{+}$ is the cost of unit power lost in line losses; and $\Clovr \in \R_{+}$ is the cost of unit absolute deviation of nodal voltage from the nominal value $\nuc{nom}{}$. 	The weight $\Cshed_i-\Cload_i$ is chosen to enable proper counting of the cost of load control when the load is disconnected. 


\paragraph*{Explanation of constraints} \  \cref{eq:attackerResourceConstraint} states that the attacker can disrupt at most $\arcm$ nodes. \cref{eq:postContingencyVoltage1} models the impact of a TN-side disruption in terms of sudden drop in substation voltage; 
\eqref{eq:dgConnectivityPostContingency} states that if the attacker disrupts a DG at node~$i$, then that DG becomes non-operational, and is \emph{effectively disconnected} from the DN; 
\eqref{eq:dgoutputAttackerAction}-\eqref{eq:dgoutputReactiveActiveConstraint} determine the feasible space for a DG's power output; \eqref{eq:dgoutputActivePositive} states that the active power output of DG is always non-negative; \eqref{eq:dgoutputReactiveActiveConstraint} states that the magnitude of a DG's reactive power output can atmost be $\etac{constant}{i}\ge0$ times its active power output;
and \eqref{eq:dgoutputAttackerAction} (resp. \eqref{eq:dgoutputOperatorAction}) combined with \eqref{eq:dgoutputActivePositive} and \eqref{eq:dgoutputReactiveActiveConstraint} state that the active and reactive power output of DG is zero if it is disconnected due to attacker (resp. operator) action. 

\cref{eq:integralityConstraints} captures the binary constraints of the connectivity variables;  \eqref{eq:loadControlSheddingConstraint} and \eqref{eq:loadControlParameterConsumptionConstraint} together model that if a load is connected to a DN, the operator may change the actual consumption of the load to a fraction of its nominal demand via direct load control; and  \eqref{eq:voltageDisconnectDG} models that a DG is  disconnected if the nodal voltage violates either of the DG's operating voltage bounds, as required by the IEEE standard rules for interconnection of DGs~\cite{ieee1547}.
Similarly, a load at node $i\in\N$ will disconnect if either of its operating voltage bounds is violated \eqref{eq:voltageDisconnectLoad}. 

\cref{eq:totalPowerConsumption} models the net nodal power consumption;  \eqref{eq:conserveRealTrue} (resp. \eqref{eq:conserveReacTrue}) is the active (resp. reactive) power conservation equation;  \eqref{eq:voltageTrue} is the voltage drop equation; and \eqref{eq:currentTrue} models the current-voltage-power relationship~\cite{baran}.


\cref{eq:currentTrue} is a non-convex equation due to which the operator subproblem becomes challenging to solve. \tcbtext{Furthermore, for a fixed operator response, the network state computed using NPF constraints (\eqref{eq:postContingencyVoltage1} and \eqref{eq:conserveRealTrue}-\eqref{eq:currentTrue}) may not even be unique. Considering a linear power flow (LPF) model instead would resolve the uniqueness issue of the network state, and allow for a  straightforward application of Benders cut. However, as stated in \cref{sec:introduction}, an analogous application of the Benders cut does not work for BiMISOCPs with binary variables in the inner problem. Nevertheless, we can address the issue of uniqueness of the network state by considering the convex relaxation of \eqref{eq:currentTrue} as follows~\cite{exactConvexRelaxation}:}
\begin{equation}\label{eq:currentApproxConvex}
	\ellc{}{ij}\nuc{}{i} \ge \Pc{}{ij}^2+\Qc{}{ij}^2 \qquad \forall \ (i,j) \in \E.
\end{equation}

Let $\Second  \coloneqq \{\second \in \{0,1\}^{\NN}  \ | \  \sum_{i\in\N }\second_i \le \arcm\}$ denote the set of feasible attacker strategies. 
Next, we can denote an operator response strategy as $\third\in \Third$, where $\Third\coloneqq \{(\lcc{}{},\pgc{}{},\qgc{}{},\kcc{}{},\kgc{}{})\in  \R^{5\NN} \ | \ \eqref{eq:dgoutputOperatorAction}-\eqref{eq:loadControlSheddingConstraint} \text{ hold}\}$. Finally, we denote the set of response strategies feasible after an attack $\second$ by $\Third(\second) \coloneqq \{\third\in\Third \ | \ \text{such that } \crefrange{eq:dgConnectivityPostContingency}{eq:dgoutputAttackerAction} \text{ hold}\}$. 

For a $\third\in\Third$, let $\Xc{}{}(\third)=\{\xc{}{}\in\R^{5\abs{\N}+3\abs{\E}} \ | \ \eqref{eq:postContingencyVoltage1},  \eqref{eq:loadControlParameterConsumptionConstraint}-\eqref{eq:voltageTrue},  \eqref{eq:currentApproxConvex} \text{ hold}\}$ denote the set of feasible post-contingency states. 
 Then, we can succinctly express the attacker-operator interaction in the presence of TN-side disturbance as follows:
\begin{align}\tag{Mm}\label{eq:Mm-cascade}
\begin{aligned}
\hspace{-0.5cm}\lossMaxmin\; \coloneqq\;  \max_{\second \in\Second}&  \quad \costMaxmin{}(\second) \\
\text{s.t.} & \quad \costMaxmin{}(\second) \coloneqq \min_{\third\in\Third(\second), \xc{}{} \in \Xc{}{}\left(\third\right)} \;  \Lossc{}\left(\third,\xc{}{}\right).
\end{aligned}
\end{align}
Here, the attacker's (resp. operator's) objective is to maximize (resp. minimize) the loss $\Lossc{}$ subject to DG and load models, nonlinear power flows, TN-side disruption, and the impact of failure captured by  $\third\in\Third(\second)$.  We refer the problem \eqref{eq:Mm-cascade}  as the $\maxlossProblem$ problem.  

One can indeed compare the solution of \eqref{eq:Mm-cascade} with the  analogous Bilevel Mixed-Integer Linear Problem (BiMILP). \tcbtext{The BiMILP is different from the BiMISOCP in two main aspects: i) The constraints involve the LPF model instead of the NPF model,  and ii) the objective function does not contain the line loss term. We distinguish the variables and the quantities computed using the LPF by the hat symbol. Thus, we denote the BiMILP problem as \problemMaxminHat, the max-min value as $\lossMaxminhat$, the  minimum post-contingency loss for a given attack $\second$ as $\costMaxmin{l}(\second)$, the network state as $\xc{l}{}$, and the set of feasible network states  $\Xc{l}{}$, and so on. }

To summarize, our problem is to determine the maximin optimal attacker-operator strategies to compute the worst-case post-contingency loss $\lossMaxmin$  for NPF model.

\subsection{Assumptions}
We assume that DN lines have positive, but small impedances, i.e., $
0 < \resistance{ij} \ll 1, 0 < \reactance{ij} \ll 1 \quad \forall \ (i,j)\in\E, $
and that voltage lower bounds are positive, i.e., $\nucc{min}{i} > 0$, $\nugc{min}{i} > 0 \quad \forall\ i\in\N$. These are rather mild assumptions and hold true for DNs in practice~\cite{exactConvexRelaxation,baran}. 

We also assume the following \emph{no reverse power flow} condition. For $i\in\N$, let $\subtree{i}\subseteq \N$ denote the subset of nodes that belong to the subtree rooted at node $i$; then: 
\begin{definition}
	We say that DN satisfies the \emph{No Reverse Power Flow} condition (NRPF) if
	\begin{equation*}
	\ssum_{j\in\subtree{i}}\ \  \ptc{}{j} \ge 0, \qquad \ssum_{j\in\subtree{i}}\ \  \qtc{}{j} \ge 0 \qquad \forall\ i\in\N. 
	\end{equation*} 
\end{definition}
Under the NRPF condition, the flows computed using either linear or nonlinear power flow constraints are non-negative, i.e. on any DN line,  power does not flow towards the substation. Hence, the name \enquote{no reverse power flow}. 

We assume that the NRPF condition holds even when all DGs are producing maximum output, i.e. $\pgc{}{i} = \pgc{max}{i} \text{ and }\qgc{}{i}=\etac{}{i}\pgc{max}{i}\ \forall\ i \in \N$. An important consequence of the NRPF condition is that the convex relaxation of \eqref{eq:currentTrue} is exact~\cite{exactConvexRelaxation}, i.e., for fixed net nodal consumption, there is a unique NPF solution such that inequality \eqref{eq:currentApproxConvex} is tight. 


\tcbtext{Note that this
	property may hold even when NRPF condition is not satisfied. For example, the convex relaxation is still exact under identical resistance-to-reactance ratio~\cite{exactConvexRelaxation}. However, under general conditions when NRPF does not hold, the value of DN's resilience as estimated using convex relaxation provides a non-trivial upper bound on the true resilience value.} 

\section{Theoretical  Results}\label{sec:technicalResults}
In this section, we present novel structural results based on power flows in radial DNs. We will use these results in \cref{sec:solutionApproach} to reduce the computational time required for solving \eqref{eq:Mm-cascade}. 

\tcbtext{
For fixed $\ptc{}{}$ and $\qtc{}{}$, let $(\Pc{}{},\Qc{}{},\nuc{}{},\ellc{}{})$ denote the network state obtained using NPF constraints  \eqref{eq:postContingencyVoltage1} and  \eqref{eq:conserveRealTrue}-\eqref{eq:currentTrue}. Again, for fixed $\ptc{}{}$ and $\qtc{}{}$, let  $(\Pc{l}{}, \Qc{l}{}, \nuc{l}{}$, $\ellc{l}{})$ denote the network state obtained using LPF constraints.
}



\tcbtext{
 Let $\Flowc{}{} = \{\Pc{}{ij},\Qc{}{ij},\ellc{}{ij}\}_{(i,j)\in\E}$ denote the set of flow quantities, and $\Voltc{}{} =  \{\nuc{}{i}\}_{i\in\N}$ the set voltage quantities. Also, let $\Flowc{l}{}$ and $\Voltc{l}{}$ denote the corresponding sets of LPF quantities.  
Let $\consumptionSet = \{\ptc{}{i},\qtc{}{i}\}_{i\in\N}$ denote the set of net nodal consumption quantities. Let  $\crossTermLineSet = \{(\Pc{}{ij},\Pc{l}{ij}), (\Qc{}{ij},\Qc{l}{ij}), (\ellc{}{ij},\ellc{l}{ij})\}_{(i,j)\in\E}$ and $\crossTermNodeSet = \{(\nuc{}{i},\nuc{l}{i})\}_{i\in\N}$ be the sets consisting of tuples each with an entry of the NPF quantity and its corresponding LPF quantity. 
}

%

\tcbtext{
 Our first proposition relates the signs and relative  magnitudes of the partial derivatives of NPF and LPF quantities with respect to net nodal consumption. 
 (We refer the reader to \Cref{app:proofs} for the proofs of the technical results.)
}
\tcbtext{
\begin{proposition}\label{prop:relatingQuantitiesToNetConsumption} 
	Under NRPF, the following hold: 
	\begin{itemize}[leftmargin=*]\renewcommand\labelitemi{--}
		\item flow quantities $f$ (resp. $\hat{f}$) computed using NPF (resp. LPF)  are increasing (resp. non-decreasing) in the net nodal consumption, 
		\item nodal voltages $v$ (resp. $\hat{v}$) for NPF or LPF are strictly decreasing in the net nodal consumption, and 
		\item the impact of a change in consumption is greater on the NPF values than for the LPF values, i.e.
	\end{itemize}
	\begin{align*}\fontsize{9}{9} \selectfont
		\begin{aligned}
		\frac{\partial f}{\partial c}  \ge  \frac{\partial\hat{f}}{\partial  c} \ge 0 > \frac{\partial \hat{v}}{\partial c} > \frac{\partial v}{\partial c} \ \ &&\forall \ (f,\hat{f})\in\crossTermLineSet, (v,\hat{v})\in\crossTermNodeSet, c\in\consumptionSet.
		\end{aligned}
	\end{align*}  
\end{proposition}
Intuitively, \cref{prop:relatingQuantitiesToNetConsumption} holds because increasing net consumption reduces the voltage at all nodes, which in turn, increases the power flows and the currents on all lines. 
}


\tcbtext{
The following proposition relates the optimal DG output to its connectivity under optimal operator response.  
\begin{proposition}\label{prop:optimalDGoutput}
	For a fixed attacker strategy, for $i\in\N$, let  $(\pgc{}{i}^\star,\qgc{}{i}^\star,\kgc{}{i}^\star)$ (resp. $(\pgc{l}{i}^\star,\qgc{l}{i}^\star,\kgc{l}{i}^\star)$) be the  optimal operator response values for the variables of DG $i$ computed using NPF (resp. LPF). Under NRPF, 
	\begin{align*}
	\begin{aligned}
	\pgc{}{i}^\star &= \pgc{max}{i}(1-\kgc{}{i}^\star),  \quad\qquad &&\pgc{l}{i}^\star &&= \pgc{max}{i}(1-\kgc{l}{i}^\star),\\ 
	\qgc{}{i}^\star &=  \etac{}{i}\pgc{max}{i}(1-\kgc{}{i}^\star), \quad &&  \qgc{l}{i}^\star &&= \etac{}{i}\pgc{max}{i}(1-\kgc{l}{i}^\star). 
	\end{aligned}		
	\end{align*} 
\end{proposition}	
}

\Cref{prop:optimalDGoutput} implies that under NRPF, the active and reactive power capability of connected DGs will be fully exhausted leaving no room for response via DG output control. An important consequence of \cref{prop:optimalDGoutput} is that the operator response can be simplified to $\lcc{}{},\kcc{}{},\kgc{}{}$ since the DG output is uniquely determined by whether it is connected or not.\footnote{\tcbtext{In practice, the DGs may not be able to generate output at their maximum capacity. In this case, the operator's loss will be even higher. Thus, the DN's resilience which we compute will be an upper bound on the true DN's resilience. } }

\tcbtext{
 Henceforth, with a slight abuse of notation, we use the notation $\Third$ to denote the projection of the set $\{\third \in \R^{5\NN}  \text{ such that } \qgc{}{i} = \etac{}{i}\pgc{}{i} = \etac{}{i}\pgc{max}{i}(1-\kgc{}{i}) \quad \forall \ i \in \N \text{ and } \eqref{eq:integralityConstraints}-\eqref{eq:loadControlSheddingConstraint} \text{ hold}\}$ onto the space of $(\lcc{}{},\kcc{}{},\kgc{}{})-$variables. Then, an operator response can be denoted by $\third = (\lcc{}{},\kcc{}{},\kgc{}{}) \in\Third$.
}

\tcbtext{
	The next proposition relates the impact of change in net nodal consumption of a  downstream node \emph{versus} an upstream node. } \tcbtext{
\begin{proposition}\label{prop:downstreamImpactIsHigher}
	Consider $k,l\in\N$ such that $k\prec l$, i.e. $k$ is a upstream of $l$.  Let $\mathcal{Z}_{kl} = \{(\ptc{}{k},\ptc{}{l}),(\qtc{}{k},\qtc{}{l})\}$. Under NRPF, the impact of change in net nodal consumption at $l$ on the flow and voltage quantities is larger than that due to an equivalent change in the net consumption at $k$, i.e.  
	\begin{align*}\fontsize{9}{9}\selectfont
	\begin{aligned}
	\frac{\partial f}{\partial c_l} &> \frac{\partial f}{\partial c_k} > 0 > \frac{\partial v}{\partial c_k} > \frac{\partial v}{\partial c_l} \ \ &&\forall \ f\in\Flowc{}{}, v \in \Voltc{}{}, (c_k,c_l)\in \mathcal{Z}_{kl}\\
	\frac{\partial \hat{f}}{\partial  c_l} &\ge \frac{\partial  \hat{f}}{\partial c_k} \ge 0  > \frac{\partial \hat{v}}{\partial  c_k} \ge  \frac{\partial \hat{v}}{\partial c_l} \ \ &&\forall \ \hat{f}\in\Flowc{l}{}, \hat{v} \in \Voltc{l}{}, (c_k,c_l)\in \mathcal{Z}_{kl}.
	\end{aligned}
	\end{align*}  
\end{proposition}}



\tcbtext{
	The next proposition relates the values of  post-contingency loss under optimal operator response (c) and the DN's resilience computed using NPF and LPF. 
\begin{proposition}\label{prop:lossNPFisHigherThanlossLPF} 
	For any attacker strategy, the minimum operator loss as computed using NPF is greater than the corresponding loss computed using LPF, i.e. 
	\begin{equation*}
		\costMaxmin{}(\second)> \costMaxmin{l}(\second)\qquad \forall \ \second\in \Second.		
	\end{equation*}
	Consequently, the resilience computed using LPF upper bounds the corresponding value computed using NPF, i.e. \begin{equation*}
		\resilienceMaxmin \le \resilienceMaxminhat. 
	\end{equation*} 
\end{proposition}}

\tcbtext{The next proposition describes the monotonicity property of DN's resilience with respect to $\arcm$.
\begin{proposition}\label{prop:monotonicity}
	If, for attacks $\second', \second'' \in \{0,1\}^{\N}$, the set of DGs attacked in $\second'$ is a subset of those attacked in  $\second''$, i.e., 
	\begin{equation} \label{eq:subsetAttack}
		\{i\in \N\ |\ \second_i' = 1\} \subseteq \{i\in \N\ |\ \second_i'' = 1\}, 
	\end{equation}
then the operator's loss due to $\second'$ would not be greater than that due $\second''$, i.e., 
	\begin{equation*} \fontsize{9}{9}\selectfont
		  \costMaxmin{}(\second')\le  \costMaxmin{}(\second'').
	\end{equation*}
	Consequently, the DN's resilience is monotonically non-increasing as attack cardinality increases. That is, if  $\resilienceMaxmin^\arcm$ denotes the DN's resilience under attack cardinality  $\arcm$, then 
	\begin{equation*}
		\resilienceMaxmin^{\arcm'} \ge \resilienceMaxmin^{\arcm''} \qquad \forall\ 0\le \arcm' \le \arcm'' \le \NN. 
	\end{equation*}
\end{proposition}}


\tcbtext{
\cref{prop:upstreamLoadConnectivityPreference} (resp. \cref{prop:upstreamDGConnectivityPreference}) relates the connectivity of loads (resp. DGs) on downstream \emph{versus} upstream nodes under optimal operator response (c). 
\begin{proposition}\label{prop:upstreamLoadConnectivityPreference} Consider  $i,j\in\N$ such that $i \prec j$. If (i) the lower voltage bound, the nominal active and reactive power demand, load control parameter, and the cost coefficient of load control at $i$ are at most the corresponding values at $j$, and (ii) the cost coefficient of load shedding at $i$ is at least as much as that at $j$, then, in an optimal operator response, the upstream load being shed implies that the downstream load is also shed. That is,
	\begin{equation}\label{eq:downstreamLoadPreference}\hspace{-0.6cm}
	\left.\begin{array}{r@{\ } l@{\qquad} r@{\ } l}
		\nucc{min}{i} &\le \nucc{min}{j},  \quad   &\pcc{max}{i} &\le \pcc{max}{j}\\	
		\lcc{min}{i} &\le \lcc{min}{j},  &\qcc{max}{i} &\le \qcc{max}{j}\\
		\Cload_i &\le \Cload_j, \qquad  &\Cshed_i  &\ge \Cshed_j	
	\end{array} \right\rbrace \implies 
	\kcc{\star}{i} \le \kcc{\star}{j}.
	\end{equation}
\end{proposition}}

\tcbtext{
\begin{proposition}\label{prop:upstreamDGConnectivityPreference}
	Consider $i,j\in\N$ such that $i\prec j$. If (i) DGs at both $i$ and $j$ are not attacked, (ii) DG at $i$ has a capacity larger than that of the DG at $j$, and (iii) the voltage lower bound at $i$ is smaller than that at $j$, then, in an optimal operator response, the upstream DG being disconnected  implies that the downstream DG is also disconnected, i.e.  
		\begin{equation}\label{eq:downstreamDGPreference}\hspace{-0.6cm}
	\left.\begin{array}{r@{\ } l@{\qquad} r@{\ } l}
	&  \quad   &\nugc{min}{i} &\le \nugc{min}{j}\\	
	\second_i &= 0,  &\second_j &= 0\\
	\pgc{max}{i} &\ge \pgc{max}{j}, \qquad  &\etac{constant}{i}  &\ge \etac{constant}{j}
	\end{array} \right\rbrace \implies 
	\kgc{\star}{i} \le \kgc{\star}{j}.
	\end{equation}
\end{proposition}}

	Propositions \ref{prop:upstreamLoadConnectivityPreference} and \ref{prop:upstreamDGConnectivityPreference} characterize the notion of keeping the more beneficial (\enquote{superior}) DN components connected. That is, if  the operator cannot keep the \enquote{superior} components connected that help reduce the overall loss and provide more flexibility in operation, then the operator must disconnect the \enquote{inferior} components first. \tcbtext{As a special case, if all other  parameters of two components (DGs or loads) are identical, then the component which is located upstream is more beneficial to the  operator than the downstream component. }  In \cref{sec:solutionApproach}, we use these results  to add cuts to the operator subproblem of \eqref{eq:Mm-cascade} and evaluate their effect on the computational speedup in \cref{sec:computationalStudies}. 

\section{Evaluating $\resilienceMaxmin$ - a Modified GBD Method} 
\label{sec:solutionApproach}


Our approach for evaluating $\resilienceMaxmin$ relies on using a modified Generalized Benders Decomposition algorithm~\cite{generalizedBenders} to solve \eqref{eq:Mm-cascade} on a reformulated problem. The overall approach is as follows. First, we argue that $\lossMaxmin$ can be obtained by solving a \emph{Min-cardinality} problem instead. Then, we implement the GBD algorithm, which decomposes the min-cardinality problem into a master (attacker) problem (an integer program) and an operator subproblem (a mixed-integer program). Then, the algorithm solves these two problems in an iterative manner, until either an optimal  min-cardinality attack is obtained or all the attacks are exhausted.

\subsection{Min-cardinality disruption problem}\label{subsec:minCardinalityDisruption}

Recall that in problem \eqref{eq:Mm-cascade}, the attacker's goal is to determine an optimal  attack of size at most $\arcm$ (attack resource). On the other hand, in the min-cardinality  problem, the attacker computes a disruption with as few attacked DN nodes as possible to induce a loss to the operator greater than a pre-specified threshold target post-contingency loss, denoted  $\ltarget$. 
\tcbtext{In fact, the min-cardinality problem and \eqref{eq:Mm-cascade} are duals of each other~\cite{abhinavVerma}. Furthermore, any procedure that can obtain an optimal solution of Min-cardinality problem (resp. \eqref{eq:Mm-cascade}), can be used to obtain an optimal solution of its dual \eqref{eq:Mm-cascade} (resp. min-cardinality problem) using a binary search on the parameter $\ltarget$ (resp. $\arcm$); see~\cite{technicalReport} for additional details. }

Now, we describe the GBD method to solve the min-cardinality problem. For given load and DG connectivity vectors $\kcc{}{}$ and $\kgc{}{}$, we define a  \emph{configuration} vector as $\configurationc{}{} \coloneqq \left(\kcc{}{}, \kgc{}{}\right)$. Given an attack vector $\second$, let $\setConfigurations(\second) \coloneqq \{ \left(\kcc{}{}, \kgc{}{}\right) \in \{0,1\}^{2\NN}   \text{ such that } \eqref{eq:dgConnectivityPostContingency} \text{ holds} \}$, i.e. $\setConfigurations(\second)$ denotes the set of all possible post-disruption configuration vectors that the operator can choose from. Then, for a fixed attack $\second$ and a fixed configuration vector $\configurationc{}{}\in \setConfigurations(\second)$, consider the following second-order cone program:
\begin{align}\label{eq:innerConvexSP}\tag{O-SOCP}
\begin{aligned}
\hspace{-0.4cm}\ploss\left(\second, \configurationc{}{}\right)\; \coloneqq\; & \mmin_{\lcc{}{}\in[0,1]^{\N}} \;  \Lossc{}\left(\third,\xc{}{}\right) \\
& \hspace{0cm} \text{s.t. } \third = \left(\lcc{}{},\configurationc{}{}\right), \third\in\Third, \xc{}{} \in \Xc{}{}\left(\third\right).
\end{aligned}
\end{align}
Note that \eqref{eq:innerConvexSP} may be infeasible as the chosen $\configurationc{}{}$ may violate  \eqref{eq:voltageDisconnectDG} or \eqref{eq:voltageDisconnectLoad} in the set of constraints $\Xc{}{}(\third)$. In this case, the value of $\ploss\left(\second, \configurationc{}{}\right)$ is set to  $\infty$. 

	Suppose that, for a given DN, we are concerned with a TN-side disturbance $\vdc{}{0}$ and a target $\ltarget$ post-contingency loss. Following \cite{abhinavVerma}, we say that  an attack-induced  disruption $\second\in\Second$ \emph{defeats} a configuration $\configurationc{}{}\in\setConfigurations(\second)$ if  $\ploss\left(\second, \configurationc{}{}\right) \ge \ltarget$, and is \emph{successful} if it defeats \emph{every} $\configurationc{}{} \in \setConfigurations(\second)$. 
We can now state the $\mincardinalityProblem$ problem as follows:
\begin{align}\tag{$\mcp$}\label{eq:mincardinalityProblem}
\begin{aligned}
\mmin_{\second\in \{0,1\}^{\N}} & &&\ssum_{i\in\N }\ \second_i\\
\text{s.t. } &&& \ploss\left(\second, \configurationc{}{}\right) \ge \ltarget  \quad \forall \  \configurationc{}{} \in \setConfigurations(\second). 
\end{aligned}
\end{align}	
If there exists an optimal solution of the problem  \eqref{eq:mincardinalityProblem}, say $\second^\star$, then it is a \emph{min-cardinality disruption} corresponding to $\ltarget$ because it is successful and has minimum number of attacked nodes. 

However, problem \eqref{eq:mincardinalityProblem} is not tractable in its current form because the number of constraints is equal to the cardinality of set $\setConfigurations(\second)$ which can be exponential in $\NN$, and verifying each constraint $\left(\ploss\left(\second, \configurationc{}{}\right) \ge \ltarget\right)$ is itself an SOCP. The GBD algorithm addresses this issue. 

\subsection{Modified Generalized Benders Decomposition} \label{subsec:bendersDecomposition}

The GBD algorithm decomposes \eqref{eq:mincardinalityProblem} into two subproblems: attacker \tcbtext{MILP} \ master problem ($\masterProblem$) and operator MISOCP subproblem ($\slaveProblem$), which are then solved in an iterative manner. In fact, in each iteration, one needs to solve ($\masterProblem$), ($\slaveProblem$), and the dual of \eqref{eq:innerConvexSP}, as discussed below. \Cref{fig:bendersDecomposition} summarizes the overall approach. \tcbtext{The numbers in round brackets indicate the order of the steps.} 

\def \hh {\textwidth}
\ifOneColumn
	\def \hh {0.7\textwidth}
\fi
\begin{figure}[htbp!]
	\resizebox{\hh}{!}{
	\begin{tikzpicture}[scale=0.7]
	\tikzset{every node/.append ={font size=tiny}} 
	
	\tikzstyle{fs} = [font=\fontsize{5}{8}\selectfont,align=left];
	\tikzstyle{rr} = [draw,rounded corners=7pt, line width=1.5pt, fs];
	\tikzstyle{rect} = [draw, line width=1.5pt, fs];
	\tikzstyle{numrect} = [ line width=1pt, dotted, fs];
	\tikzstyle{connector} =[->,>=latex,red!80,line width=1.2pt,text=black,fs]
	\tikzstyle{biconnector} =[<->,>=latex,red!80,line width=1.2pt,text=black,fs]

	\node[rect] (obimip) {Original \\BiMISOCP \eqref{eq:Mm-cascade}}; 
	
	\node[rect, below left=-0.6cm and 0.7cm of obimip] (refbimip) {Min-cardinality \\ 
		 problem \eqref{eq:mincardinalityProblem}}; 
	
	\node[rect,below=1.8cm of refbimip,xshift=-0cm] (amip) { Attacker\\
		MILP \eqref{eq:masterProblemBenders}
	}; 
	
	\node[minimum width=0cm,left=0.3cm of amip.north] (amipcoord){};

	\node[rect,right=-0.8cm of amip, yshift=1.45cm](omip) {Operator\\ \eqref{eq:subproblemBenders} \\with Prop. \ref{prop:upstreamLoadConnectivityPreference},  \ref{prop:upstreamDGConnectivityPreference} }; 
	
	\node[minimum width=0cm,below=0.4cm of omip.west] (omipCoord){};
	
	\node[draw,diamond,line width =1.5pt, right=0.6cm of omip,align=left,aspect=2,scale=0.4,font=\fontsize{18}{0}\selectfont] (dec) {
		$\begin{aligned}
		&\Lossc{}\left(\third^\star,\xc{\star}{}\right)\\
		&\ge \ltarget
		\end{aligned}$
	}; 
	\node[right=0.5cm of dec,fs] (exit) {Exit};
	
	\draw[connector] (amip.140) to  [out=90,in=200] node[above,xshift=-0.1cm]{$\second^\star$} (omip.190)[yshift=-.2cm];
	
	\draw[connector] (obimip.210) to  [out=190,in=30] (omip.80);

	\path let \p2=(amip), \p1=(dec) in node[rect] (olp) at (\x1,\y2) {Operator SOCP\\  $\P(\second^\star,\kappa^\star)$ \eqref{eq:innerConvexSP}
	}; 
	
	\node[rect, below left= 0.2cm and -0.6cm of olp,fs] (olpdual) {Operator Dual\\ SOCP \eqref{eq:OSP2}
	}; 
	
	\foreach \x/\y/\z in {obimip/obimipn/1, refbimip/refbimipn/2, amip/amipn/3, omip/omipn/4, olp/olpn/6, olpdual/olpdualn/7}{
		\node[numrect,left=-0.1cm of \x, yshift=0.3cm](\y) {(\z)};
	}
		\node[numrect,left=-1cm of dec, yshift=0.6cm](decn) {(5)};
		

	\node[fs,yshift=0.2cm] () at (dec.east) {yes}; 
	\node[fs,xshift=-0.3cm] () at (dec.south) {no};
	
	\draw[connector] (obimip.west) to[out=180,in=10] node[right,yshift=-0.15cm]{} (refbimip.east);
	\draw[connector] (refbimip.210) to [out = 230,in = 120] (amip.150);
	
	\draw[connector] (omip.east) to node[above,xshift=0.2cm]{$\third^\star,\xc{\star}{}$} (dec.west);
	\draw[connector] (amip.east) to node[above]{$\second^\star$} (olp.west);
	\draw[connector] (dec.south) to node[right,yshift=-0.05cm]{$\third^\star = (\beta^\star,\kappa^\star), \xc{\star}{}$} (olp.north);
	\draw[connector] (dec.east) --++(0.75,0); 
	
	\draw[connector] (olp.south) to  [out=270,in=0] node[above,fs,align=center,xshift=-0.1cm,yshift=0.0cm]{} (olpdual.east);
	\draw[connector] (olpdual.west) to  [out=180,bend left=30] node[below,fs,align=left,xshift=-0.6cm,yshift=0.4cm]{
		Modified\\
		generalized\\Benders cut} (amip.290);
	\end{tikzpicture}}
	\caption{Computational approach to solve \eqref{eq:Mm-cascade}.}
	\label{fig:bendersDecomposition}
\end{figure}
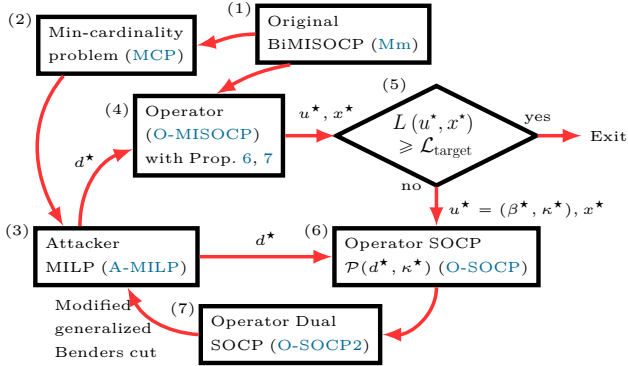

The attacker \tcbtext{MILP }  can be written as follows: 
\vspace{2pt}
\begin{align}\tag{$\masterProblem$}\label{eq:masterProblemBenders}
\begin{aligned}
\mmin_{\second\in \Second} &\quad &&\ssum_{i\in\N }\ \second_i\\
\text{s.t. } &&& \text{set of generalized Benders cuts}. 
\end{aligned}	
\end{align}
The master problem is initialized with only the integrality and budget constraints on the attack variables, and without any generalized Benders cut (to be defined in \eqref{eq:generalizedBendersCut}). In each iteration, solving the master problem \eqref{eq:masterProblemBenders}, which is a bounded MILP, if feasible,  yields an attack $\second^\star$. Then, this attack vector is used as an input parameter for the operator subproblem ($\slaveProblem$). For a fixed attack $\second^\star$, the operator subproblem is the same as the inner  problem of $\eqref{eq:Mm-cascade}$:
\begin{align}\tag{$\slaveProblem$}
\begin{aligned}
\mmin_{\third\in\Third(\second^\star),\xc{}{} \in \Xc{}{}(\third)} &\quad &&\Lossc{}\left(\third,\xc{}{}\right)\\
 \text{s.t.} &&& \eqref{eq:downstreamLoadPreference},\eqref{eq:downstreamDGPreference}.
\end{aligned}\label{eq:subproblemBenders}
\end{align}	
\tcbtext{Note that \eqref{eq:downstreamLoadPreference} and \eqref{eq:downstreamDGPreference} result due to \cref{prop:upstreamLoadConnectivityPreference} and \ref{prop:upstreamDGConnectivityPreference}; see~\cref{sec:technicalResults}}. The problem \eqref{eq:subproblemBenders} is also a bounded MISOCP because the load and DGs have bounded feasible space. If \eqref{eq:subproblemBenders} is feasible, it yields an optimal operator response $\third^\star$ and network state $\xc{\star}{}$ for the disruption $\second^\star$. If the operator's loss $\Lossc{}\left(\third^\star,\xc{\star}{}\right)$  exceeds the target loss $\ltarget$, the algorithm terminates having successfully determined an optimal min-cardinality attack. Otherwise, $\Lossc{}\left(\third^\star,\xc{\star}{}\right) <\ltarget$ which implies that  $\second^\star$ is not a successful disruption. In this case, we need to  generate a generalized Benders cut to eliminate $\second^\star$ from the feasible space of \eqref{eq:masterProblemBenders}.

\def \obj {L}
\def \lhsMatrix {A}
\def \lhssocpMatrix {E}
\def \lhssocpVector {f}
\def \rhssocpCoeff {g}
\def \rhssocpScalar {h}
\def \rhsMatrix {B}

\tcbtext{Note that problem \eqref{eq:innerConvexSP} with parameters ($\second^\star$, $\configurationc{}{}^\star$) can be simplified and rewritten as the following problem:}
\tcbtext{	\begin{align}
	\nonumber\min_{w} \  &\transpose{c}w &&\\
	\text{s.t. }  &\lhsMatrix w &&\ge b + \rhsMatrix\second^\star && : (\lambda) \tag{O-SOCP2}\label{eq:OSP2}\\
	& \nonumber\normSquared{\lhssocpMatrix^jw} &&\le \transpose{\rhssocpCoeff^j}w  && : (\alpha^j,\beta^j) &&\ \forall\ j \in \N,
	\end{align}
where $\normSquared{\cdot}$ is the 2-norm; $w$ is the primal decision vector variable; $\lhsMatrix$, $\rhsMatrix$, and $\lhssocpMatrix^j$s are matrices; and $b$,  $\lhssocpVector^j$s and $\rhssocpCoeff^j$s are vectors of appropriate dimensions.
Also, $\lambda$ and ($\alpha^j, \beta^j$)  for $j \in \N$ are the dual variables corresponding to the linear and SOCP inequalities, respectively. The $\abs{\N}$ second-order cone constraints correspond to \eqref{eq:currentApproxConvex}. }
 
 \tcbtext{Thus, the dual of problem  \eqref{eq:OSP2} can be simply written as:
	\begin{alignat}{8}\fontsize{9}{10}\selectfont
	\hspace{-1cm}
	\max_{\substack{\lambda\ge\zero, \alpha^j, \beta^j}} \quad  &&& \transpose{\left(b+\rhsMatrix\second^\star\right)}\lambda \nonumber\\
	 \text{s.t. } &&& \normSquared{\alpha^j} \le \beta^j \qquad\quad\forall\ j \in \N \tag{D-SOCP2}\label{eq:dspSP2} \\
	 &&& c - \transpose{A}\lambda+\ssum_{j\in\N} \big(\transpose{\lhssocpMatrix^j}\alpha^j - \beta^j \rhssocpCoeff^j\big) = \zero && \nonumber
	\end{alignat}}
\tcbtext{We solve the dual problem (thanks to strong duality, the optimal values are the same) in  \eqref{eq:dspSP2} to compute $\ploss\left({\second}^\star,{\configurationc{}{}}^\star\right)$ and an optimal dual solution $(\lambda^\star,\alpha^{j\star},\beta^{j\star})$. This furnishes a \emph{generalized Benders cut}, which is added to the master problem in the next iteration. 
In particular, if the dual problem in \eqref{eq:dspSP2} has an optimal solution $(\lambda^\star,\alpha^{j\star},\beta^{j\star})$, and its optimal value is $\obj^\star$, then 
\begin{equation}\label{eq:generalizedBendersCut}
	\transpose{\left(b+\rhsMatrix\second\right)}\lambda^\star \ge \epsilon
\end{equation} }
is the desired generalized Benders cut  where $\epsilon$ is a non-negative number. In a classical generalized Benders cut the value of $\epsilon$ is 0. If the inner subproblem of \eqref{eq:Mm-cascade} were convex, such a cut would indeed be useful in eliminating sub-optimal attacker strategies~\cite{generalizedBenders}. However, this cut is not useful in the presence of discrete inner variables, i.e. it does not eliminate any attack vector. 

Hereafter, we refer to the generalized Benders cut in \eqref{eq:generalizedBendersCut} as simply the Benders cut. An exact expression for \eqref{eq:generalizedBendersCut} is provided in~\cite{technicalReport}. Note that $\second^\star$ does not satisfy \eqref{eq:generalizedBendersCut} when $\epsilon>0$ because $\transpose{\left(b+\rhsMatrix\second^\star\right)}\lambda^\star = \ploss\left(\second^\star, {\configurationc{}{}}^\star\right) = \obj^\star < \obj^\star + \epsilon$, where the first equality holds because of strong duality in second-order cone programs. \tcbtext{Thus, choosing $\epsilon >0$ is a \emph{modification} to the Benders cut which helps eliminate $\second^\star$ from attacker's set of feasible strategies. } 
However, due to numerical issues,  an off-the-shelf solver can \enquote{stall} at run-time, and may be unable to generate dual vector values required for the Benders cut. To address this issue, we add the following cut: 
\begin{equation}
	\ssum_{\left(i\in\N:\second^\star_i = 1\right)}\second_i + \ssum_{\left(i\in\N:\second^\star_i = 1\right)} (1-\second_i) \le \NN-1, 
\end{equation} 
which ensures that $\second^\star$ is eliminated. 

Thus, in each iteration, we eliminate suboptimal attacks from the feasible space of \eqref{eq:masterProblemBenders}. Hence, the new master problem obtained by adding a Benders cut is a stronger relaxation of \eqref{eq:mincardinalityProblem}. Consequently, we get a progressively tighter lower bound on the minimum cardinality of the attack as the iteration continues, until we get a successful attack. Since there are a finite number of attacks, whether successful or not, the GBD algorithm is bound to terminate.





Note that the overall algorithm, as depicted in \cref{fig:bendersDecomposition}, is also applicable for solving the BiMILP \problemMaxminHat. In this case, instead of solving an MISOCP and SOCP, the algorithm would simply solve an MILP and an LP. 

\subsection{Choosing $\epsilon$ based on criticality parameter - A heuristic}\label{subsec:choosingEpsilon}
The Benders cut, when simplified, is of the form $\sum_{i\in\N}\dualCoefficient_i\second_i \ge \epsilon^j$, where $\dualCoefficient = \transpose{\lambda^\star}B$ is the coefficient vector, and $\epsilon^j>0$ is a scalar chosen for the $j^{th}$ added Benders cut. The choice of $\epsilon^j$ in the Benders cut is an important issue in our implementation of the GBD algorithm. One way would be to choose a constant value of $\epsilon$ for each Benders cut. However, if we choose too large an $\epsilon$ then many attacks (possibly including the optimal attacks) might be eliminated from the set of feasible attacker strategies in \eqref{eq:masterProblemBenders}. This introduces an approximation  error as a result of which, the  obtained min-cardinality attack may not be optimal. If we choose too small an $\epsilon$, then in each iteration only the current min-cardinality attack vector is eliminated resulting in  performance no better than brute force enumeration over all attacks.

To address this issue, we modify the Benders cut by proposing a novel heuristic to assign varying values for $\epsilon$ in each iteration. 
 Suppose that in iteration $j$, the optimal attack vector obtained is $\second^j$ and the dual coefficient vector is $\dualCoefficient^j = \transpose{\lambda^\star}B$; see~\eqref{eq:generalizedBendersCut}. Let $\arcm^j \coloneqq \sum_i\second^j_i$ be the cardinality of $\second^j$. Let $\permutation$ be a permutation of nodes such that $\dualCoefficient^j_{\permutation(1)} \ge \dualCoefficient^j_{\permutation(2)} \ge \cdots\ge \dualCoefficient^j_{\permutation(\NN)}$, with the ties broken by lexicographical ordering. Here $l= \sigma^j(i)$ indicates that node $l\in\N$ has the $i^{th}$ highest value in the vector $C^j$. Let \tcbtext{$\startParameter\in [0\isep\NN-1]$}\  be a parameter, which we call a \emph{criticality parameter}. We use $\startParameter$ to obtain $\epsilon$ for selecting critical DG nodes to attack. Let $\seqEnd \coloneqq \min(\NN, \startParameter+\arcm^j)$ and $\seqStart = \seqEnd - \arcm^j+1$. Then, one can choose $\epsilon$ for the $(j+1)^{th}$ iteration as follows:  
\begin{equation*}
	\epsilon^{j+1} = \  \underbracket{\dualCoefficient^j_{\permutation(\seqStart)} \ +\  \dualCoefficient^j_{\permutation(\seqStart+1)} \ +\  \cdots \ +\  \dualCoefficient^j_{\permutation(\seqEnd)}}_{\arcm^j \ \text{terms}}.   
\end{equation*}
Essentially, we exclude the top  $\min(\startParameter, \NN - \arcm^j)$ values, and then take the sum of next $\arcm^j$ coefficients. As $\startParameter$ increases, the $\epsilon^j$ value decreases, thereby allowing the GBD algorithm to explore more number of attacks. As a result, one would expect the optimality gap to be lower and the computational time to be higher than the case when $\startParameter$ is small. 

An intuitive reason for why this heuristic works is as follows. By \cref{prop:downstreamImpactIsHigher}, we get the insight that the downstream nodes in a DN are critical. Therefore, the attacker may attack as many downstream nodes as he can subject to his resource constraint. However, in this case the attacker may fail to exploit the cascading nature of the attack. Specifically, the attacker may be better off by not disrupting a few downstream nodes, and instead using his budget on compromising a few upstream nodes. Consequently, the downstream DGs, which are anyway more likely to face voltage bound violations, may be disconnected due to the operator response. That is why choosing a lower value of $\epsilon^j$ allows the GBD algorithm to explore attacks that do not compromise the most critical nodes (as suggested by the dual coefficients in the Benders cut). Essentially, the dual coefficients $C^j$ do not capture the cascading effect due to further disconnection of other DGs and loads because we fix the configuration vector for solving the SOCP. In other words, $C^j$s do not represent the true \enquote{criticality} of the DG nodes because they ignore the cascading effects. Therefore, varying the criticality parameter  $\startParameter$ allows the algorithm to explore attacks on DGs whose criticality as indicated by $C^j$ value is less. As we show in \cref{sec:computationalStudies}, the GBD algorithm with variable value for $\epsilon$ takes significantly fewer iterations compared with brute force or the GBD algorithm with a constant $\epsilon$. 

\section{Evaluating $\resilienceNoResponse$ - A Two-Step Approach}\label{sec:autonomousDisconnects}

\subsection{Autonomous disconnect model - Response (b)}

To model the network state under response (b), i.e. uncoordinated autonomous disconnects, we propose the following two-step approach. In the first step, we compute the subset of DGs which will autonomously disconnect due to the attacker-induced failure as well as due to the resulting voltage bound violations. In the second step, we determine the subset of loads facing voltage bound violations caused by the DG disconnects in the first step. Since voltage bound violations are typically indicative of faults, DGs are disconnected a lot sooner than the loads as a precautionary measure to avoid feeding current to a fault. This is why we focus on only DG disconnections in the first step. \tcbtext{Thus, our approach allows us to compute the worst-case loss due to a cascade. This is the main difference between our approach and the multi-round cascade algorithm described in~\cite{bienstockBook}.} 

Now, we provide the details of our two-step approach. For a fixed operator action $\third\in\Third$, let $\Zc{}{}$ denote the set of network states $\xc{}{}$ that satisfy the constraints \eqref{eq:postContingencyVoltage1},  \eqref{eq:loadControlParameterConsumptionConstraint}, \eqref{eq:voltageDisconnectDG}, \eqref{eq:totalPowerConsumption}-\eqref{eq:voltageTrue} and \eqref{eq:currentApproxConvex}. Note that  $\Xc{}{}(\third)\subseteq\Zc{}{}(\third)$ because  $\Xc{}{}(\third)$ has an additional constraint \eqref{eq:voltageDisconnectLoad}. For a fixed attacker action $\second\in\Second$, let $(\third^{\star}_\intermediateState(\second),\xc{\star}{\intermediateState}(\second))$ denote the intermediate autonomous disconnect action and the corresponding network state. We can extract the information about disconnected DGs and the nodal voltages from this intermediate action and network state  $(\third^{\star}_\intermediateState(\second),\xc{\star}{\intermediateState}(\second))$ to compute the final autonomous disconnect action and the post-contingency state denoted by  $(\third^{\star}_\noresponseSmall,\xc{\star}{\noresponseSmall})$. We formulate a problem to compute $(\third^{\star}_\intermediateState(\second),\xc{\star}{\intermediateState}(\second))$ as follows: 
\begin{align}\tag{P-IN}\label{eq:noresponseIntermediate}
\begin{aligned}
&& \mmin_{\third_\intermediateState,\xc{}{\intermediateState}} &\  &&\Lossc{}\left(\third_\intermediateState,\xc{}{\intermediateState}\right) \\ 
& &\text{s.t. }&&&
\third_\intermediateState \in \Third(\second), \quad && \xc{}{\intermediateState}\in\Zc{}{}\left(\third_\intermediateState\right)\\
&&&&& \lcc{in}{i} = 1 \quad &&  \forall \ i\in\N,
\end{aligned}
\end{align} 
where the intermediate state does not require the loads to satisfy the voltage bound constraint. Note that the load control parameters are set to unity to model the fact that under autonomous disconnections, the operator will not be able to exercise load control. 

Next, to compute $(\third^{\star}_\noresponseSmall,\xc{\star}{\noresponseSmall})$, we extract the value of DG connectivity vector $\kgc{in\star}{}$ and voltage data $\nuc{in\star}{}$ from the intermediate action-state pair  $(\third^{\star}_\intermediateState(\second),\xc{in\star}{}(\second))$. Then, we use this value to parameterize the following problem:  
\begin{align}\tag{P-FN}\label{eq:noresponseFinal} 
\begin{aligned}
\hspace{-0.3cm}
 \min_{\third_\text{nr},\xc{}{\text{nr}}} &&& \Lossc{}\left(\third_\noresponseSmall,\xc{}{\noresponseSmall}\right)\\ 
\text{s.t. } &&& \third_\noresponseSmall \in \Third, \quad && \xc{}{\noresponseSmall}\in\Xc{}{}\left(\third_\noresponseSmall\right) \\
 &&&\lcc{nr}{i} = \kcc{nr}{i} &&\forall \ i\in\N\\
 &&& \kgc{nr}{i} \ge \kgc{in\star}{i}(\second) &&\forall \ i\in\N\\
&&& \kcc{nr}{i} \ge \nuc{}{i}-\nuc{in\star}{i}(\second) \;  &&\forall \ i\in\N\\
 &&& \eqref{eq:downstreamLoadPreference}, \eqref{eq:downstreamDGPreference}.
\end{aligned}
\end{align} 
The optimal solution of the above problem provides us $(\third^{\star}_\noresponseSmall,\xc{\star}{\noresponseSmall})$, i.e the final autonomous disconnect action and the post-contingency state. 

\cref{algo:noDefenderResponse} summarizes the execution of the two-step approach. It takes as input an initial attack-induced contingency $\second$, and generates automatic disconnect actions for one or more components due to the uncontrolled cascade. Note that the load control parameter $\lcc{}{i}=1$ throughout the cascading disconnects of DGs, unless the load becomes fully disconnected, in which case it switches to $\lcc{}{i}=0$. The final connectivity vector $\third^\star_\noresponseSmall$ corresponds to a situation where all the connected components satisfy voltage bounds, and can be used to compute the corresponding post-contingency loss $\Lossc{}\left(\third^\star_\noresponseSmall, \xc{\star}{\noresponseSmall}\right)$. 

\begin{algorithm}[htbp!]
	\fontsize{9}{10}\selectfont
\caption{Uncontrolled cascade under response (b)}
\begin{algorithmic}[1]
\Require attacker action $\second$ (initial contingency)
\State $\third^\star_\noresponseSmall,\xc{\star}{\noresponseSmall}\gets \Call{GetCascadeFinalState($\second$)}{}$
\Function{GetCascadeFinalState}{$\second$}
\State Compute  $\third^\star_\intermediateState(\second),\xc{\star}{\intermediateState}(\second)$ by solving \eqref{eq:noresponseIntermediate}
\State Extract parameters $(\kgc{in\star}{},\nuc{in\star}{})$ from $(\third^\star_\intermediateState,\xc{\star}{\intermediateState})$
\State Instantiate \eqref{eq:noresponseFinal} with parameters $(\kgc{in\star}{},\nuc{in\star}{})$
\State Solve \eqref{eq:noresponseFinal} to compute the final state  $\third^\star_\noresponseSmall,\xc{\star}{\noresponseSmall}$
	\State \Return $\third^\star_\noresponseSmall,\xc{\star}{\noresponseSmall}$ 
		\EndFunction
		\end{algorithmic}
		\label{algo:noDefenderResponse}
		\end{algorithm}

\tcr

\subsection{Randomized algorithm for lower bounding $\lossNoResponse$}\label{subsec:autonomousDisconnectLoss}
			
			For each cardinality $\arcm$, we  can compute the worst case loss under response (b) using brute force. However, that would require evaluating loss over combinatorially many ${\NN \choose \arcm}$ attacks. Therefore, we present a randomized algorithm to compute  worst case loss under the autonomous disconnections; see~\cref{algo:randomAttacksNoResponse}. 
			\begin{algorithm}[htbp!]
				\fontsize{9}{10}\selectfont
				\caption{Random attacks and approximately worst case attack for autonomous disconnections}
				\begin{algorithmic}[1]
					\Require $\nperm$ (number of random permutations)  
					\State Initialize $\YY = \zero_{\NN\times\nperm}$ and $\VV = \zero_{\NN}$
					\For {$t \in [1 \isep \nperm]$}		
					\State Generate a random permutation $\sigma$ of nodes $\N$
					\State Reset $\second = \zero$ 
					\For {$\arcm = 1 \isep \NN$}
					\State Set $\second_{\sigma(\arcm)} = 1$ \Comment{$\arcm$ cardinality attack}
					\State $\left(\third_\noresponseSmall, \xc{}{\noresponseSmall}\right)$ $\gets$ \Call{\text{GetCascadeFinalState}}{$\second$} \\ \Comment{Refer~\cref{algo:noDefenderResponse} for \textsc{GetCascadeFinalState}}
					\State $\YY[\arcm,t] \gets \Lossc{}(\third_\noresponseSmall, \xc{}{\noresponseSmall})$
					\EndFor
					\EndFor
					\For {$\arcm \in [1 \isep \NN]$}
					\State $\VV[\arcm] \gets \max_{t\in [\nperm]} \YY[\arcm,t]$ 
					\EndFor
					\State \Return $\YY, \VV$ 
				\end{algorithmic}
				\label{algo:randomAttacksNoResponse}
			\end{algorithm}
			
			The algorithm performs the following steps: for each random permutation of nodes, for each attack cardinality $\arcm$, it disrupts the first $\arcm$ nodes in that permutation, and computes the loss due to autonomous component disconnects (using~\cref{algo:noDefenderResponse}). Then, for each attack cardinality, it  chooses the maximum among all computed losses.
As shown in \cref{sec:computationalStudies}, for any randomly chosen attack of cardinality $\arcm<\NN$, if we disrupt one more DG, then the loss incurred under autonomous disconnections will increase. This monotonicity of increasing loss for increasing attack cardinality cannot be shown if we simply choose $\NN+1$ random attacks of cardinalities $\arcm \in [0 \isep \NN]$, and plot the loss values \textit{vs.} $\arcm$. This is the main idea behind \cref{algo:randomAttacksNoResponse}. In \cref{sec:computationalStudies}, we implement the modified GBD algorithm and \cref{algo:randomAttacksNoResponse} to compute the \emph{value of timely response}, i.e. $\resilienceMaxmin-\resilienceNoResponse$.

\section{Computational Results}
\label{sec:computationalStudies}
	
	\tcbtext{We refer the reader to the appendix for the setup of our computational study. }
	
	Our computational results are organized to show: (a) the value of timely operator response compared to autonomous disconnections; (b) comparison of the solutions of our GBD approach with the optimal solution (generated for small networks by brute force); and (c) the scalability of our approach to larger networks. 

\ifOneColumn
	\def \www {0.225\textwidth}
	\def \hhh {\www}
\else
	\def \www {0.45\textwidth}
	\def \hhh {\www}
\fi

\paragraph*{Solution accuracy of the modified GBD method}
For a fixed cardinality $\arcm$, we compute the optimal loss $\loss^\star$ using brute force over all disruptions. \tcbtext{For  $\NN=36$ node network, the  brute force method finished after $\approx$24 hours. Thus, under a time limit of 24 hours, the exhaustive search was possible only for $\NN=36$ node network. }  Then, we use $\loss^\star$ as the parameter $\ltarget$ for the problem \eqref{eq:mincardinalityProblem}. If the GBD algorithm applied to \eqref{eq:mincardinalityProblem} computes a successful attack with the same cardinality $\arcm$, then indeed we have obtained the optimal attack of cardinality $\arcm$.  
%

\begin{figure}[htbp!]
	\subfloat[$\NN = 24$]{\label{fig:bruteVsBendersA}
			\includegraphics[width=\www,height=\hhh]{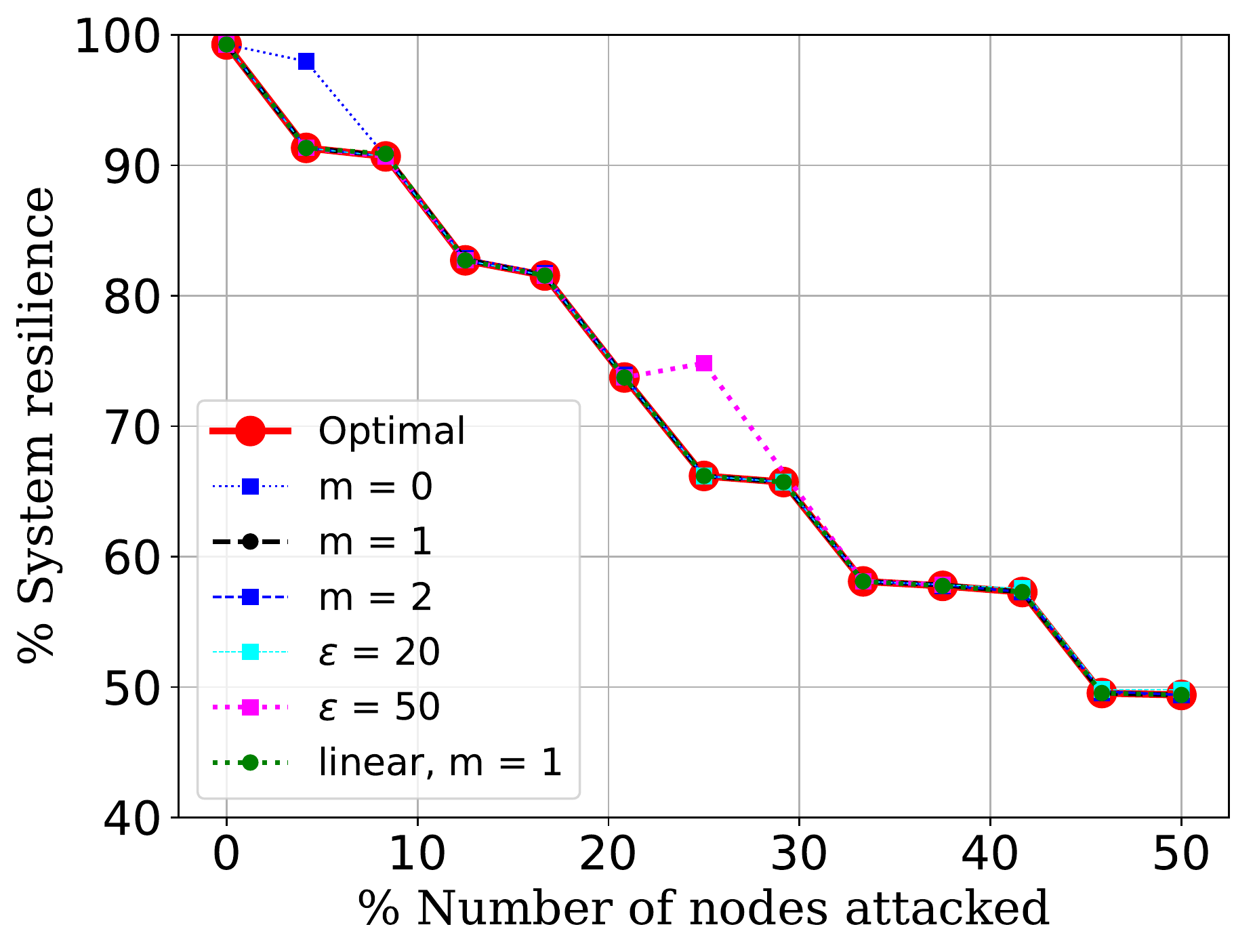}
	}
	\subfloat[$\NN = 36$]{\label{fig:bruteVsBendersB}
			\includegraphics[width=\www,height=\hhh]{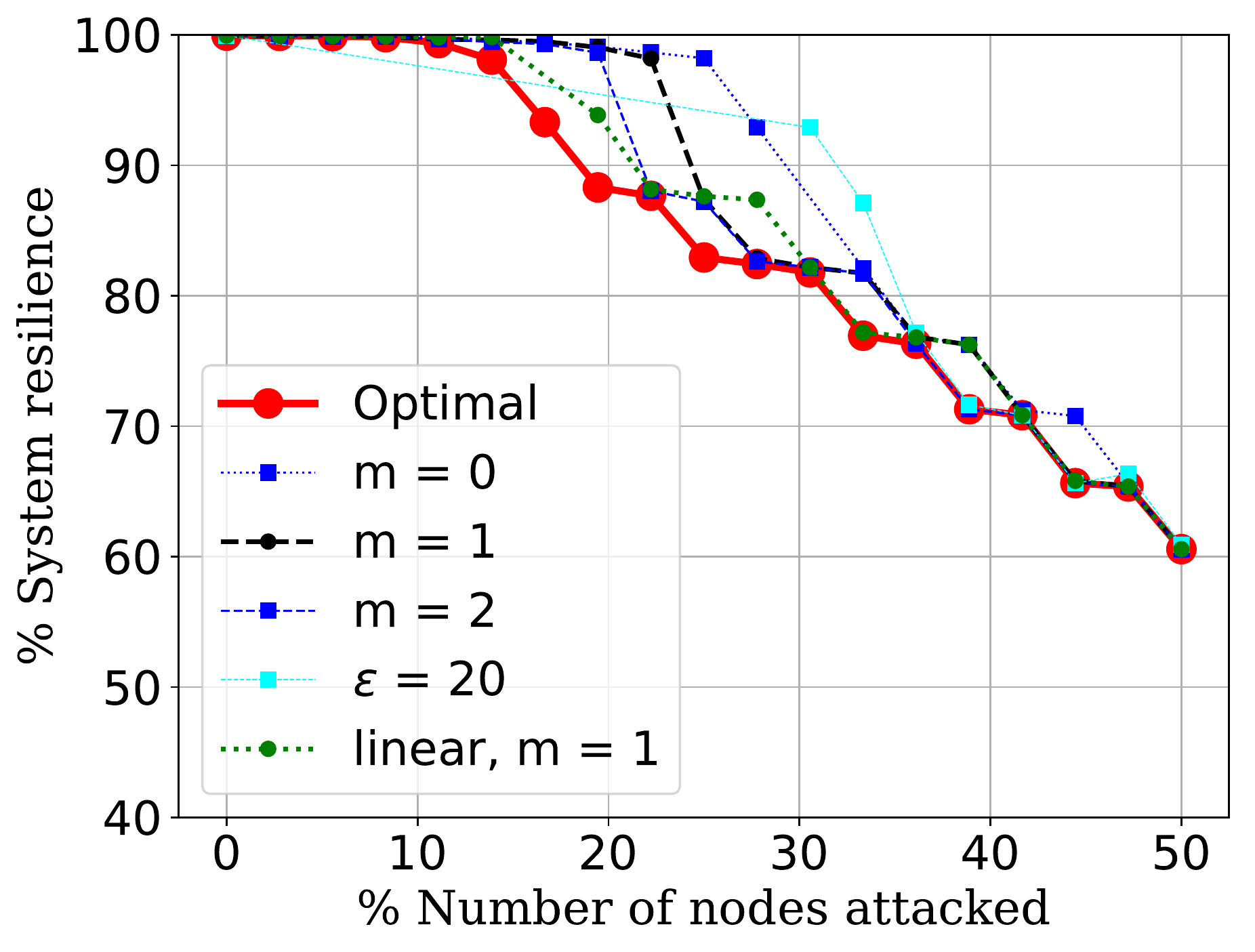}
%
	}
	\caption{\tcbtext{Accuracy of GBD algorithm in computing $\resilienceMaxmin$.
}}
	\label{fig:bendersVsBrute}
\end{figure}
\Cref{fig:bendersVsBrute} shows that our GBD method with variable $\epsilon$ choices performs very well in computing optimal attacks. 
\tcbtext{The accuracy of the modified GBD method decreases as $\epsilon$ increases as shown by curve corresponding to $\epsilon = 50$ in \cref{fig:bruteVsBendersA}, and increases as the criticality parameter $\startParameter$ increases as can be seen in~\cref{fig:bruteVsBendersB};  see~\cref{subsec:choosingEpsilon}. To generate the curve ``linear, $\startParameter = 1$", we first applied the modified BD method to \problemMaxminHat to obtain the optimal attacker strategy, and then computed the operator's post-contingency loss under NPF for that attacker strategy. We explain how we chose the fixed values for $\epsilon$ in the discussion of next experiment. }

\paragraph*{Performance of the modified GBD method}\ 
\Cref{tab:performanceEpsilon} compares the computational time and solution accuracy of the GBD method with constant $\epsilon$ and variable $\epsilon$ choices. We also show the results for our solution approach applied to the BiMILP, where the optimal attacks are then used to evaluate operator's loss using NPF constraints. 
\def \optimalityGap {\text{Gap}}
\tcbtext{In \Cref{tab:performanceEpsilon}, $\optimalityGap$ denotes the percentage gap between the cardinality of  attack obtained by GBD and that of optimal attack; and $n_{iter}$ the number of iterations the algorithm took to reach either convergence or the iteration limit. The first three rows correspond to fixed $\epsilon$ choices.
The next three rows correspond to variable $\epsilon$ choices. In the last two rows, $\text{l}$ indicates that the optimal attacker strategy was computed using the BD method for the \problemMaxminHat, and then  reevaluated using NPF constraints. Results show that GBD method with variable $\epsilon$ provides significant computational speedup, while still retaining solution accuracy.} \tcbtext{${}^\ast$ indicates that algorithm was terminated after it reached the iteration limit of 10000.}

\tcbtext{Consider the case of $\NN = 24$. Since half the nodes have DGs, there are $2^{12} = 4096$ possible attacks. If we choose fixed $\epsilon$ value of 10 or 20, the GBD approach takes explores nearly all 4096 attacks resulting in a performance as bad as the brute force method. For a large $\epsilon = 50$, it takes 1596 iterations (which is more than a third of all the attacks), and provides 8.33\% gap. On the other hand, using the variable $\epsilon$ approach, with $\startParameter = 0$, it converges in 17 iterations while still providing 8.33\% gap. }

\begin{table}[htbp!]
	\centering 
	\def \numberOfRandomAttacks {10}
	\caption[]{\tcbtext{Computational performance vs. $\epsilon$ choices.} }\label{tab:performanceEpsilon}\footnotesize
	\resizebox{1\textwidth}{!}{
		\begin{tabular}{|c|r|r|r|r|r|r|}
			\hline 
			& \multicolumn{3}{|c|}{$\NN=24$} &  \multicolumn{3}{|c|}{$\NN=36$} \\
			\hline 
			& \tcbtext{$\optimalityGap$} & \tcbtext{$n_{iter}$} & Time & \tcbtext{$\optimalityGap$} & \tcbtext{$n_{iter}$} & Time \\
			\hline 
			$\epsilon = 10$ &  0.00\% & 4096 & 3112.4s  &  0.00\% &10000* & 12537s \\
			$\epsilon = 20$ &  0.00\% & 4094 & 3098.9s  & 50.0\% & 4190 & 5045.3s \\
			$\epsilon = 50$ &  8.33\% &  1596 & 815.3s  & -- & -- & -- \\
			\hline
			$\startParameter=0$  & 8.33\% & 17 &  1.49s  & 27.78\% & 22 & 6.44s   \\
			$\startParameter=1$  &  8.33\% & 123 &  13.43s  & 22.22\% & 230 & 46.12s  \\
			$\startParameter=2$  & 5.56\% & 496 &  85.90s  &  16.67\% & 1828 & 825.44s  \\
			\hline
			l, $\startParameter=0$  &  8.33\% & 22 &  2.25s  & 27.78\% & 29 & 4.97s  \\
			l, $\startParameter=1$  &  8.33\% & 161 &  22.65s  & 16.67\% & 198 & 54.42s  \\
			\hline		
		\end{tabular}
	}
\end{table}

\Cref{tab:performanceCuts} shows the benefits of adding cuts \eqref{eq:downstreamLoadPreference} and \eqref{eq:downstreamDGPreference} on the computational time required to solve the operator subproblem. These experiments were carried with variable $\epsilon$ choices for parameter $\startParameter \in \{0,1\}$. Adding the cuts \eqref{eq:downstreamLoadPreference}-\eqref{eq:downstreamDGPreference} become significantly beneficial for large networks, as $\startParameter$ increases.
\begin{table}[htbp!]
	\centering 
	\def \numberOfRandomAttacks {10}
	\caption[]{\tcbtext{Computational speedup due to cuts. }
		}\label{tab:performanceCuts}\footnotesize
	\begin{tabular}{|c|c|c|c|c|}
	\hline 
	&	& $\NN = 24$ & $\NN = 36$ & $\NN = 118$  \\
		\hline
		\multirow{2}{1cm}{$\startParameter = 0$}	&	with cuts & 2.30s  & 4.60s & 87.34s \\
		\cline{2-5}
		&	no cuts  & 1.87s  & 4.55s &  88.83s \\
		\hline 
\multirow{2}{1cm}{$\startParameter = 1$}	&	with cuts & 9.74s  & 24.66s & 613.42s \\
		\cline{2-5}
	&	no cuts  & 10.66s  & 28.29s &  2949.04s \\
		\hline
	\end{tabular}
\end{table}

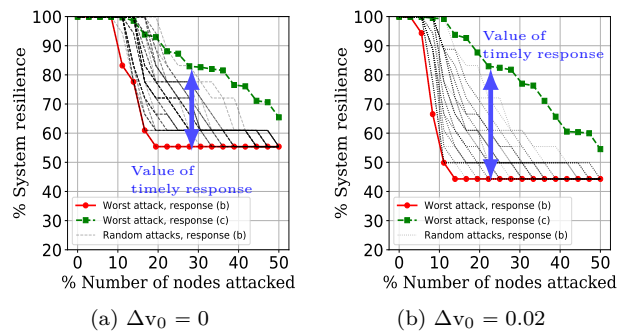
\begin{figure}[htbp!]
			\vspace{-5mm}
			\def \swidth {3}
			\def \drawgrid {\draw[step=0.25,gray,  draw opacity = 0] (0,0) grid (\swidth,\swidth);}
			\def \drawvalue at (#1,#2,#3); {\draw[blue!70, line width = 1.8pt, >= {latex[length=0.5pt,width=1pt]}, <->] (#1,#2) -- (#1,#3);}
			\def \drawresilience at (#1,#2); {\node[blue!90,align=left,font=\fontsize{5}{7}\selectfont] at (#1,#2) {Value of\\timely response};}
			\tikzstyle{nomargin} = [inner sep=0, outer sep = 0]
			\subfloat[$\Delta \mathrm{v}_0 = 0$]{				\label{fig:noResponseVsSequentialN36Dv0}
				\begin{tikzpicture}
				\node[nomargin] at (\swidth/2,\swidth/2) {
					\includegraphics[width=\www,height=\hhh,trim=0cm 0cm 0cm 0cm,clip=true]{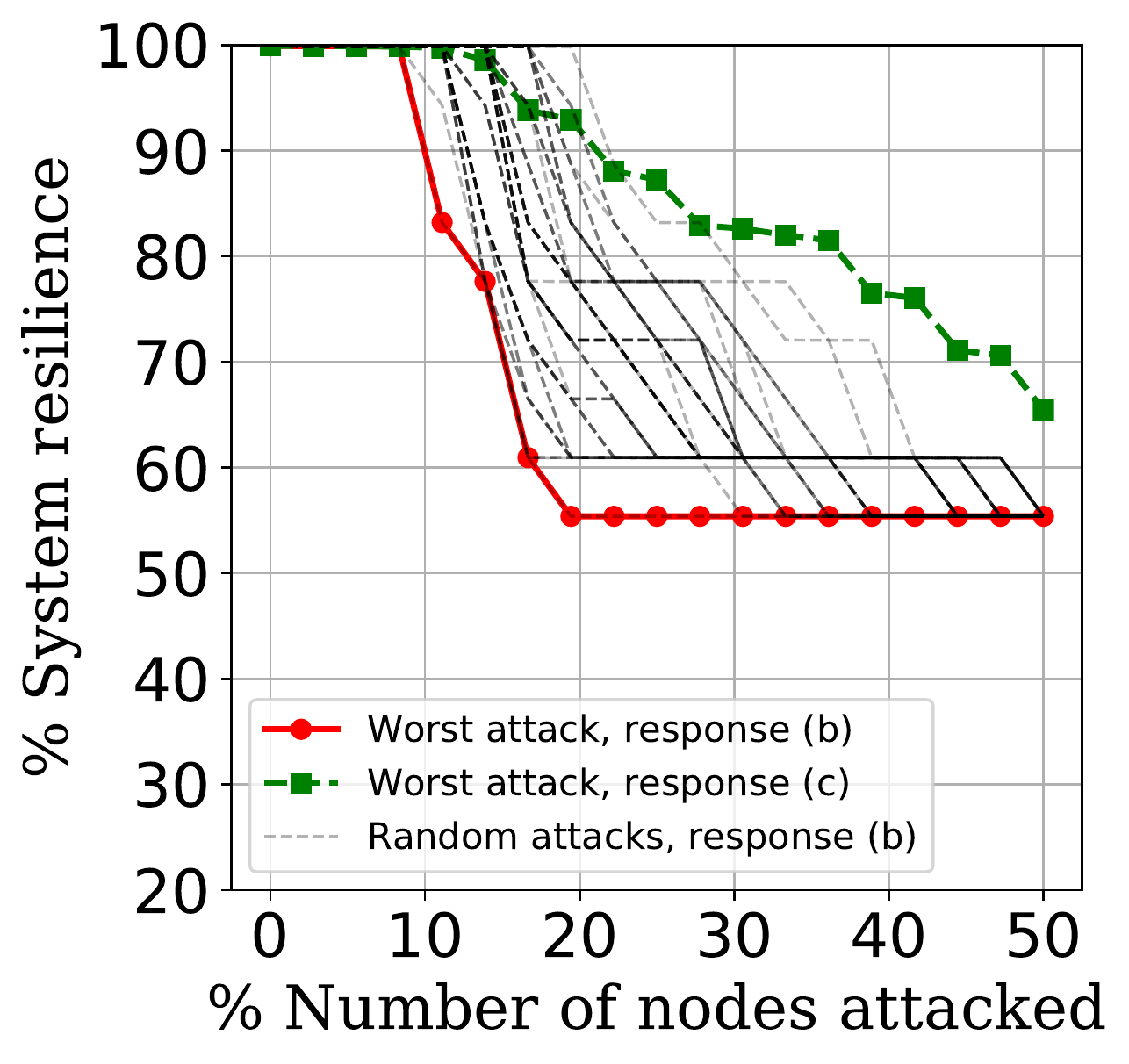}
				};
				
				\drawvalue at (2,1.5,2.6);
				\drawresilience at (2,1.1);
				
				\end{tikzpicture}				
			}
			\subfloat[$\Delta \mathrm{v}_0 = 0.02$]{ 			\label{fig:noResponseVsSequentialN36Dv3}
				\begin{tikzpicture}
				\node[nomargin] at (\swidth/2,\swidth/2)  {
					\includegraphics[width=\www,height=\hhh]{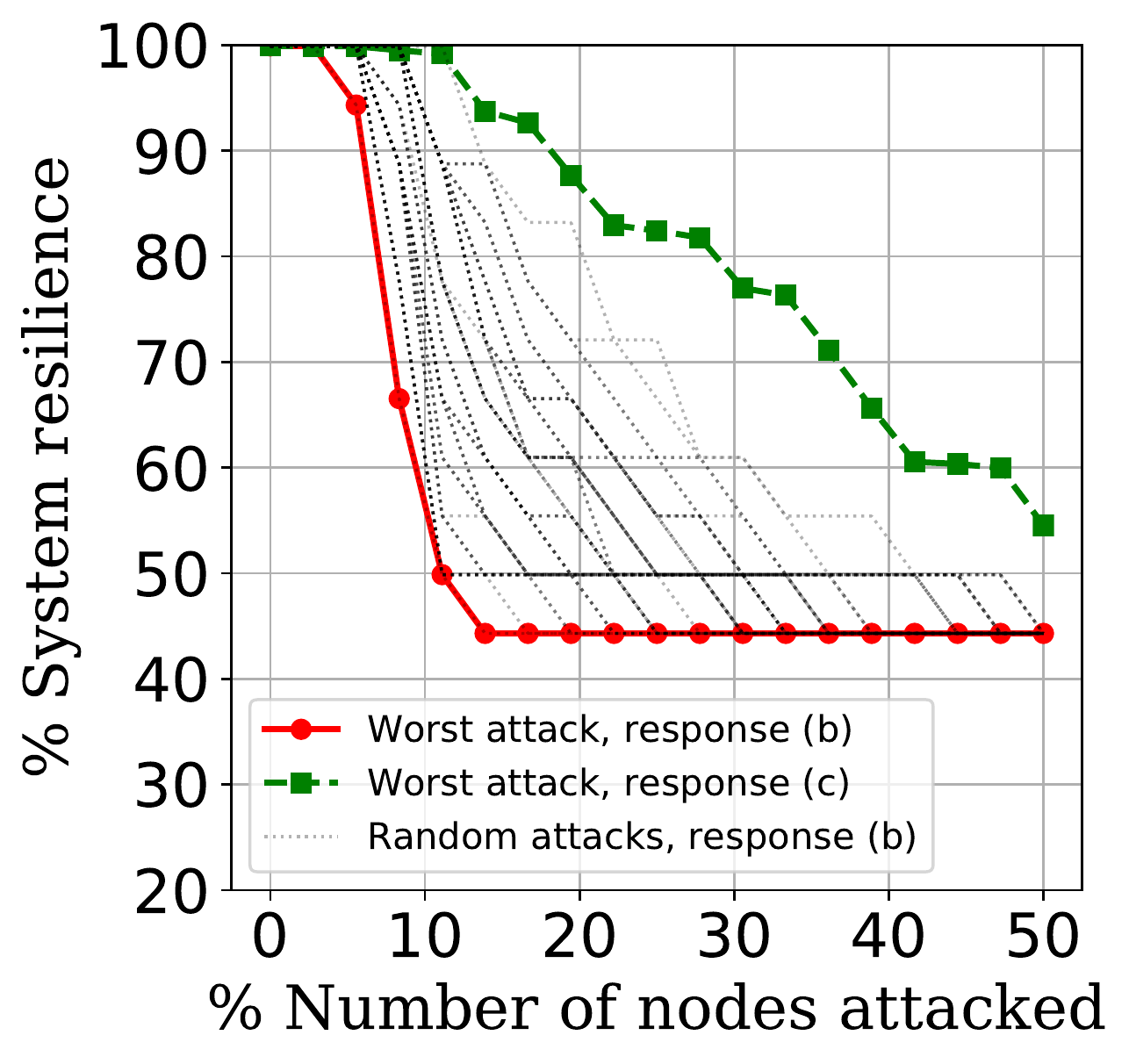}
				};
				
				\drawvalue at (1.7,1.1,2.6);
				\drawresilience at (2.4,2.9);
				
				\end{tikzpicture}
			}
		\vspace{-2mm}
			\caption{\tcbtext{Value of timely response ($\NN = 36$).}}
			\label{fig:noResponseVsSequential}
		\end{figure}

	\paragraph*{Value of timely response} \ 
	Recall that in \cref{sec:introduction}, we used post-contingency loss to define the resilience metric for SA system response ($\resilienceMaxmin$) and autonomous disconnection ($\resilienceNoResponse$) cases; and that $\resilienceMaxmin\ge \resilienceNoResponse$. \cref{fig:noResponseVsSequential} compares the resiliency values for the two cases (response (c) versus autonomous disconnection (b)) for varying number of nodes attacked, where computation of $\resilienceMaxmin$ (resp. $\resilienceNoResponse$) involves using the GBD algorithm (resp. \Cref{algo:noDefenderResponse}). In \Cref{fig:noResponseVsSequential}, the resilience curve due to response (b) under random attacks is obtained by using~\cref{algo:randomAttacksNoResponse} in the~\Cref{subsec:autonomousDisconnectLoss}.

		
		Indeed, under autonomous disconnections, we find that the voltage bound violations cause even the non-disrupted DGs to disconnect resulting in a cascade. However, under operator response, the SA detects these voltage bound violations, and \emph{preemptively} exercises load control and/or disconnects the loads/DGs to reduce the total number of non-disrupted DGs from being disconnected, and minimize the impact of the attack. 
		The difference between the two resiliency curves gives the value of timely response via the SA system. The intermediate curves in \cref{fig:noResponseVsSequential} correspond to the DN resilience under random attacks and autonomous disconnections. 
		Finally, when both a TN-side disturbance and a DN attack are simultaneous, the resilience metric of the DN decreases; see  \Cref{fig:noResponseVsSequentialN36Dv3}. 
		



	\paragraph*{Scalability of GBD algorithm}
		We tabulate the computational time required by the GBD algorithm to compute min-cardinality attacks for different network sizes and varying values of the resilience metric  $\resilienceTarget = 100\left(1-\ltarget/\lcompleteShed\right)$; see~\Cref{tab:bendersPerformance}.  
		Note that even for $\NN = 118$ nodes, which has $2^{118}$ configuration vectors, the GBD algorithm finishes computations in $\approx$10 minutes. In comparison, for $\NN=36$ node network, the brute force method took $\approx$24 hours. 
		The failure cases in \cref{tab:bendersPerformance} correspond to the cases where there does not exist an attack vector that exceeds the target loss values. \tcbtext{The realized resilience metric can significantly fall short of the target resilience metric ($\resilienceTarget = 100\left(1-\ltarget/\lcompleteShed\right)$); for e.g., when the attack cardinality changes from 8 to 9, the resilience for 36-node network decreases sharply from 98.18\% to 87.97\%.
		This means that the 36-node DN is at least 85\% (actual value 87.97\%) resilient to $\arcm=9$ cardinality attacks.}
		\begin{table}[hbtp!]
		\centering
		\caption[]{\tcbtext{Scalability of the modified GBD algorithm. }
			 }\label{tab:bendersPerformance}\footnotesize
		\resizebox{0.98\textwidth}{!}{
			\begin{tabular}{|@{\ }c@{\ }|@{\ }l@{\ }|@{\ }l@{\ }|@{\ }l@{\ }|}
				\hline 
				\multicolumn{4}{|c|}{\parbox{1.1\textwidth}{\textbf{Entries are resilience metric of DN (in percentage), number of iterations (written in brackets), time (in seconds), attack cardinality.}}}\\
				\hline 
				$\resilienceTarget$ 
				& $\NN = 24$ & $\NN = 36$ & $\NN = 118$  \\
				\hline 
				$99$ & 91.33, (3), 1.46, 1  & 98.18, (111), 13.01, 8 & 98.94, (10), 10.6, 6 \\
				\hline
				$95$  & 91.33, (3), 1.46, 1  & 87.97, (112), 13.26, 9 &  94.19, (19), 15.89, 14 \\
				\hline
				$90$ & 82.78, (8), 1.96, 3  & 87.97, (112), 13.26, 9 & 89.89, (29), 23.29, 23 \\
				\hline
				$85$  & 82.78, (8), 1.96, 3  & 82.58, (122), 16.36, 11  &  84.97, (95), 90.75, 39 \\
				\hline $80$ & 74.61, (18), 2.93, 5 & 76.94, (137), 20.69, 13 &  79.71, (86), 613.42, 52 \\
				\hline
				$75$  & 74.61, (18), 2.93, 5  & 71.05, (171), 32.35, 15  &  Failure \\
							\hline $70$  & 66.41, (16), 0.31, 6  & 65.43, (25), 0.67, 18 &   \\
				\hline
				$65$ & 58.17, (54), 8.01, 8  & 60.56, (230), 56.65, 18 &   \\
				\hline		
				$55$  & 49.53, (112), 17.13, 11  & Failure &   \\
				\hline		
				$45$ & Failure & &   \\
				\hline		
			\end{tabular}
		}
	\end{table}

	
	


\section{Concluding remarks}
\label{sec:conclusions}


In this article, we developed a computational approach to evaluate the resilience of DNs under a class of cyberphysical disruptions. We considered an attack model that involves a TN-side voltage disturbance, and DN-side supply-demand disturbance. We formulated the overall problem as a BiMISOCP, and developed a solution approach based on a modification of the GBD method. This modification entails introducing a criticality parameter. Our approach for solving BiMISOCPs with binary variables in the inner problem fills an existing gap in the literature, and can be applied to other resource allocation problems in power systems. We also estimated the value of timely operator response which involves preemptive load control or component disconnections implemented via substation automation.
Future work involves extending our attacker-operator interaction model to DNs with microgrid islanding  capabilities, and using the proposed approach to determine operator strategies for faster system performance recovery after a cyberphysical disruption event.

\section*{Acknowledgements}
We thank the anonymous reviewers for their detailed and useful suggestions. 

\bibliographystyle{IEEEtran}

			\bibliography{IEEEAbbr.bib}

\normalsize

\renewcommand{\baselinestretch}{1}

\begin{IEEEbiography}[{\includegraphics[width=1in,height=1.25in,clip,keepaspectratio]{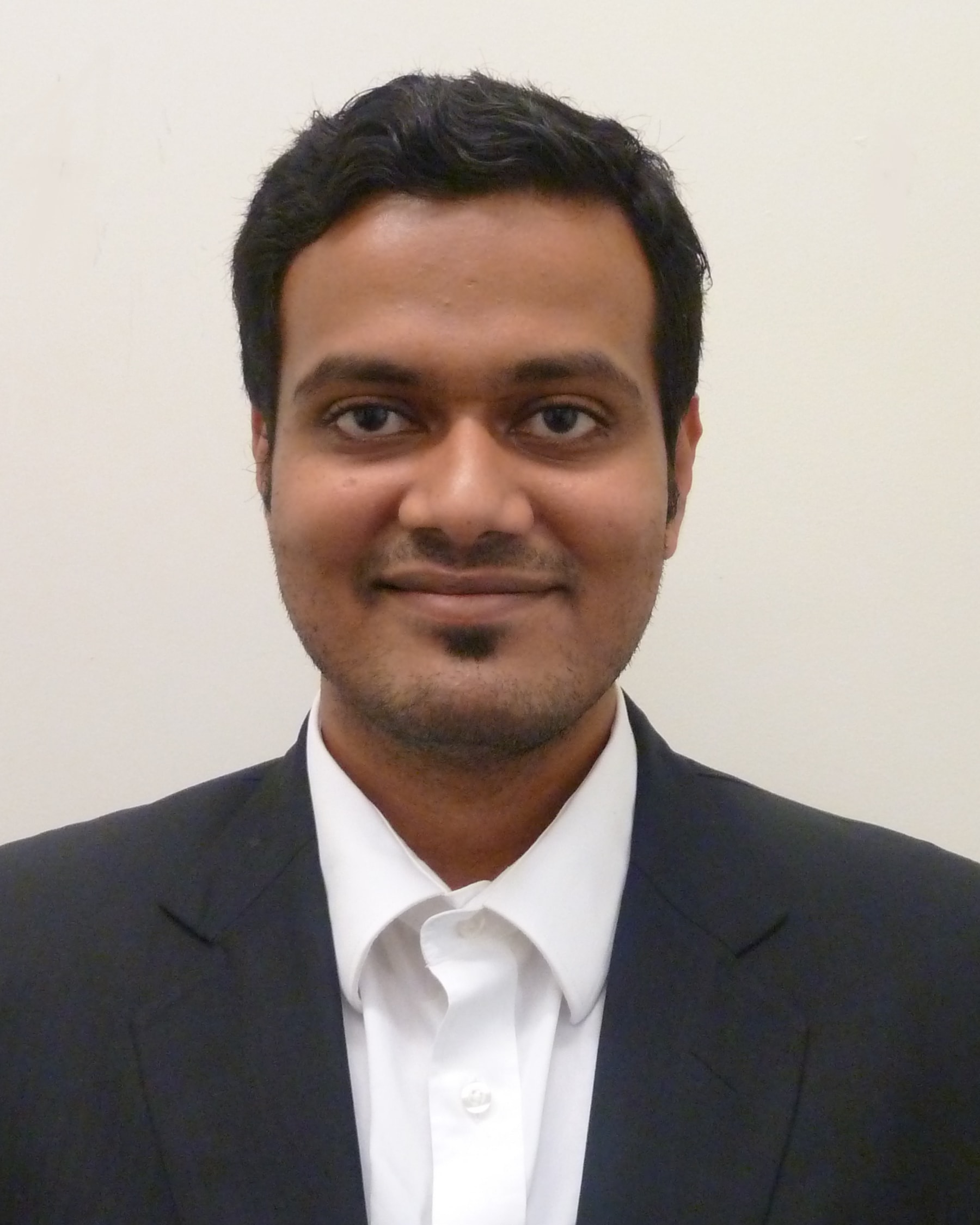}}]{Devendra Shelar}
	is a Postdoctoral Associate in LIDS at MIT. He works on developing resilient control algorithms for cyberphysical systems against extreme weather events and security attacks. His research leverages ideas from large-scale optimization, scheduling, and game theory. Dr. Shelar received his Ph.D. in Computational Science and Engineering from MIT, 2019. 
\end{IEEEbiography}
\begin{IEEEbiography}[{\includegraphics[width=1in,height=1.25in,clip,keepaspectratio]{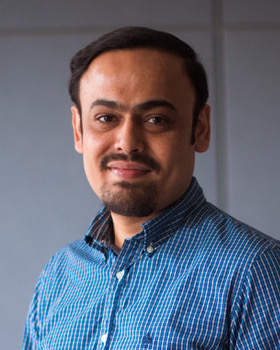}}]{Saurabh Amin}
	is a Associate Professor in the Department of Civil and Environmental Engineering. He is a member of the LIDS at MIT. He received his Ph.D. in Systems Engineering from the UC Berkeley in 2011. His fields of expertise include control and optimization, applied game theory, and networks. His research focuses on the design and implementation of resilient monitoring and control algorithms for networked infrastructures systems.
\end{IEEEbiography}
\begin{IEEEbiography}[{\includegraphics[width=1in,height=1.25in,clip,keepaspectratio]{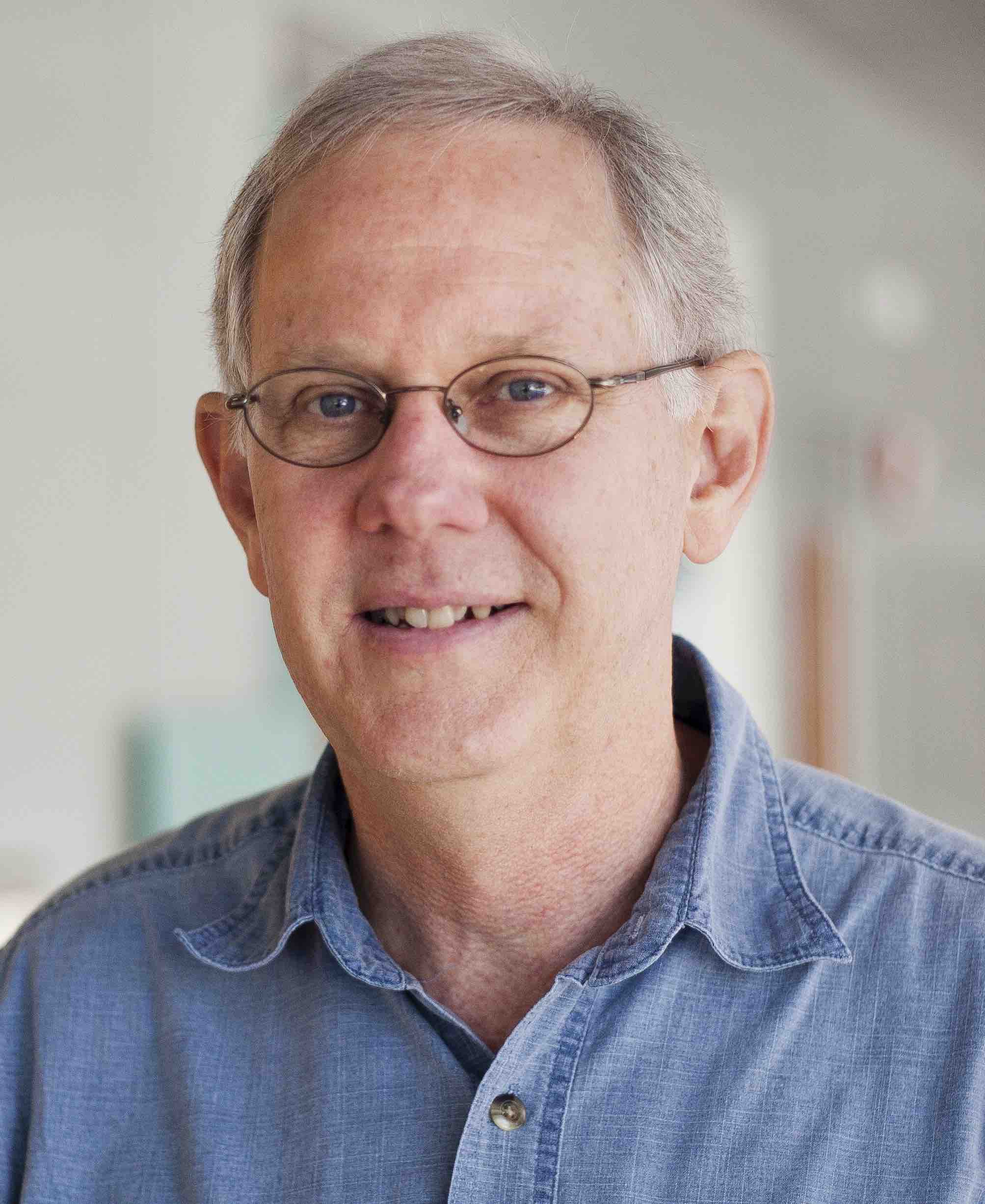}}]{Ian A. Hiskens}
	 is the Vennema Professor of Engineering in the
	Department of Electrical Engineering and Computer Science, University
	of Michigan, Ann Arbor. His research interests lie at the intersection
	of power system analysis and systems theory. He is a Fellow of IEEE
	and a Fellow of Engineers Australia. Dr. Hiskens is a Chartered
	Professional Engineer in Australia and the 2020 recipient of the
	M.A. Sargent Medal from Engineers Australia.
\end{IEEEbiography}
\appendix

\begin{table}[hbtp!]
	\caption{Table of Notations.} 
	\label{tab:notationsTable1}
	\centering
	\def \xx {0.5cm}
	\def \xxx {5.0cm}
	\def \xxxx {7.0cm}
	\def \xxxxx {1cm}
	\def \xxxxxx {10cm}
	\def \xxxxxxx {6.5cm}
	\ifOneColumn
	\def \ww {0.5\textwidth}
	\else
	\def \ww {\textwidth}
	\fi 
	\resizebox{\ww}{!}{
		\begin{tabular}{lp{7.5cm}}
			\multicolumn{2}{l}{\textbf{DN  parameters}} \\
			$0$ & substation node label\\
			\small		$\N$ & set of non-substation nodes in DN \\
			$\E$ & set of edges in DN \\
			$\NN = \abs{\N}$ & number of non-substation nodes in DN\\
			$\j$ &  complex square root of $-1$, $\j = \sqrt{ -1}$ \\
			$\nuc{nom}{}$ & nominal squared voltage magnitude (1 pu)\\
			$\nuc{}{0}$ & squared voltage magnitude at substation node\\
			\multicolumn{2}{l}{\textbf{Nodal quantities of node $i \in \N $}} \\
			$\nuc{}{i}$ & squared voltage magnitude at node $i$\\
			$\pcc{max}{i} +\j\qcc{max}{i}$ & nominal demand at node $i$  \\
			$\pgc{max}{i} +\j\qgc{max}{i}$ & nominal generation at node $i$ \\
			$\etac{constant}{i}$ & $ \max_{\qg_i,\pg_i\ne 0} (\abs{\qg_i}/pg_i)$ maximum ratio of absolute reactive power to active power\\
			$\pcc{}{i} + \j \qcc{}{i}$ & actual power consumed at node $i$  \\
			$\pgc{}{i} + \j\qgc{}{i}$ & actual power generated at node $i$ \\
			$\ptc{}{i} + \j\qtc{}{i}$ & net power consumed at node $i$ \\
			$\nucc{min}{i},\nucc{max}{i}$ & lower, upper voltage bounds for load  at node $i$\\
			$\nugc{min}{i},\nugc{max}{i}$ & lower, upper voltage bounds for DG  at node $i$\\
			$\kgc{}{i}$ & 0 if DG at node $i$ is connected to DN; 1 otherwise \\
			$\kcc{}{i}$ & 0 if load at node  $i$ is connected to DN; 1 otherwise \\		
			$\lcc{}{i}$ & fraction of demand satisfied at node $i$\\
			$\lcc{min}{i}$ & lower bound of  load control parameter $\lcc{}{i}$\\
			$\xc{}{} \in  \R^{6\NN+1}$ &  $\xc{}{} = \left(\ptc{}{}, \qtc{}{}, \Pc{}{}, \Qc{}{}, \nuc{}{}, \ellc{}{}\right)$ the network state\\
			
			\multicolumn{2}{l}{\textbf{Parameters of edge $(i,j) \in \E $}}  \\
			$\Pc{}{\edge} +\j\Qc{}{\edge} $ & power flowing from node $i$ to node $j$ \\
			$\resistance{\edge}, \reactance{\edge}$ & resistance and reactance of line $(i,j) \in \E$  \\
			$\ellc{}{\edge} $ & square of magnitude of current on line $(i,j)$ \\
			\multicolumn{2}{l}{\textbf{Precedence relationship between nodes $i,j \in \N, i\ne j$}}  \\
			$i \prec j$ & Node $i$ precedes node $j$ if $i$ lies on the path connecting $j$ and the substation node $0$\\
			\multicolumn{2}{l}{\textbf{Cyber-physical failure parameter}} \\
			$\vdc{}{0}$ & Drop in substation voltage due to transmission network-side disturbance. \\
			\multicolumn{2}{l}{\textbf{Attack variables}} \\
			$\second \in \{0,1\}^{\N}$ & $\second_i =1$ if DG at node $i$ is disrupted; 0 otherwise. \\
			\multicolumn{2}{l}{\textbf{Operator response variables}} \\
			$\third$ & $\third = \left(\lcc{}{},\pgc{}{},\qgc{}{},\kcc{}{},\kgc{}{}\right)$ an operator  response \\
			\multicolumn{2}{l}{\textbf{Generic math notation}} \\
			$[a \isep b]$ & integer interval set for $a,b\in\Z$ 
	\end{tabular}}
	
\end{table}
\normalsize
\subsection{Setup for Computational Study}
\def \hp {\alpha}

We consider three networks: 24 node, and modified IEEE 36 node and 118 node networks.  Each line has an identical impedance of $\resistance{ij} = 0.01, \reactance{ij} = 0.02$. Half of the nodes have a DG and half have a load. Hence, the maximum cardinality of an attack in our computational study will be half the number of the nodes in the DN.  
Consider a parameter  $\hp = \frac{6}{\NN}$. Before the contingency, each DG has active power output of $\pgc{max}{i} = \hp$, and each load has a demand of $\pcc{max}{i} = 1.25\hp$. Thus, we assume 80\% DG penetration since the total DG output is 80\% of the total demand. The voltage bounds are $\nucc{min}{i} = 0.9$, $\nucc{max}{i} = 1.1$, $\nugc{min}{i} = 0.92$ and $\nugc{max}{i} = 1.08$. The reactive power values are chosen to be exactly one third that of the corresponding active power value, i.e. a 0.95 (lagging) power factor for each load and DG. The values are chosen such that the total net active power demand in the DN is 0.75 pu, and the lowest voltage in the network before any contingency is  close to $\nugc{min}{}$. The maximum load control parameter is $\lcc{min}{i} = 0.8$, i.e. at most 20\% of each load demand can be curtailed. For the sake of simplicity, we assume that all DGs and loads are homogeneous. The values of cost coefficients are chosen to be ${\Cload} = 100/\pcc{max}{i}, \Clovr = 100$ and ${\Cshed} = 1000/\pcc{max}{i}$.

\subsection{Proofs of Technical Results in \cref{sec:technicalResults}}\label{app:proofs}
For $i\in\N$, let $\pathFromRoot{i}\subseteq \E$ denote the subset of DN edges on the path from the substation node $0$ to node $i$. For $i,j\in\N$, let $\commonResistance{ij}$ (resp. $\commonReactance{ij}$) denote the sum of resistances (resp. reactances) of the edges common to $\pathFromRoot{i}$ and $\pathFromRoot{j}$, i.e., 
\begin{align*}
\begin{aligned}
\commonResistance{ij} &\coloneqq \ssum_{\left((k,l)\in\pathFromRoot{i}\cap\pathFromRoot{j}\right)} \resistance{kl}, \qquad \commonReactance{ij} \coloneqq \ssum_{\left((k,l)\in\pathFromRoot{i}\cap\pathFromRoot{j}\right)} \reactance{kl}. 
\end{aligned}
\end{align*}

Let $\subtree{i}\subseteq \N$ be the subset of nodes that form the subtree rooted at node $i$, which includes node $i$, and let $\subtreeEdges{i}\subseteq\E$ be the subset of edges that form the subtree $\subtree{i}$. Then, the following equations can be derived using recursion on the radial tree topology. 
\begin{alignat}{8}
\label{eq:realFlowLinear}	\Pc{l}{ij} &= \ssum_{k\in\subtree{j}}\ptc{}{k}  && \forall\ (i,j)\in\E \\
\label{eq:reacFlowLinear}	\Qc{l}{ij} &= \ssum_{k\in\subtree{j}}\qtc{}{k} && \forall\ (i,j)\in\E\\
\label{eq:voltageLinear}	\nuc{l}{j} &= \nuc{nom}{} - 2\ssum_{k}\left(\commonResistance{jk}\ptc{}{k} +\commonReactance{jk}\qtc{}{k}\right)  &\quad& \forall\ j\in\N \\
\label{eq:realFlowFunction}	\Pc{}{ij} &= \Pc{l}{ij} +  \ssum_{(k,l)\in\subtreeEdges{i}}\resistance{kl}\ellc{}{kl}  && \forall\ (i,j)\in\E \\
\label{eq:reacFlowFunction}	\Qc{}{ij} &= \Qc{l}{ij} +  \ssum_{(k,l)\in\subtreeEdges{i}}\reactance{kl}\ellc{}{kl}  && \forall\ (i,j)\in\E 
\end{alignat}
\begin{equation}\label{eq:voltageFunction}
\begin{split}
\nuc{}{j} = \nuc{l}{j} - 2\ssum_{(k,l)\in\E}(\commonResistance{jl}\resistance{kl} + \commonReactance{jl}\reactance{kl})\ellc{}{kl} + \\ \ssum_{(k,l)\in\pathFromRoot{j}}(\resistance{kl}^2 +\reactance{kl}^2)\ellc{}{kl}
\end{split}
\end{equation}
Thus, we can write $(\Pc{l}{},\Qc{l}{},\nuc{l}{})$ as functions of $(\ptc{}{},\qtc{}{})$ and  $(\Pc{}{},\Qc{}{},\nuc{}{})$ as functions of $(\ptc{}{},\qtc{}{},\ellc{}{})$. Furthermore, we have assumed the NRPF condition. Hence, as shown in \cite{exactConvexRelaxation}, the  NFPF solution is unique. Thus, even $\ellc{}{}$ can be considered a function of $(\ptc{}{},\qtc{}{})$. 


Consider the iterative Backward-Forward Sweep (BFS) algorithm~\cite{backwardForwardSweep} used to compute the NPF values, which we modify to consider the TN-side voltage disturbance. Let  $(\Pc{\ts}{},\Qc{\ts}{},\nuc{\ts}{},\ellc{\ts}{})$ be the values computed in $\ts^{th}$ iteration of the FBS algorithm. 

\textbf{Initialization:} 
\begin{alignat}{8}
\nuc{0}{i} &= \nuc{nom}{}-\vdc{}{0}\quad &&\forall\ i\in\N\backslash\{0\}\\
\nuc{\ts}{0} &= \nuc{nom}{}-\vdc{}{0}\quad &&\forall\ \ts\in[1\isep \nperiod]\\
\ellc{0}{ij} &= 0, \Pc{0}{ij} = \Pc{l}{ij}, \Qc{0}{ij} = \Qc{l}{ij} \qquad &&\forall\ (i,j)\in\E. 
\end{alignat}

\textbf{Backward Sweep:} Starting from the leaf nodes to the substation node, compute:
\begin{alignat}{8}
\label{eq:currentSweep}		\ellc{\ts}{ij} &= \left({(\Pc{\ts-1}{ij})}^2+{(\Qc{\ts-1}{ij})}^2\right)/\nuc{\ts-1}{i} \qquad &&\forall\ (i,j)\in\E\\ 
\label{eq:realFlowSweep}		\Pc{\ts}{ij} &= \ptc{}{j}+\resistance{ij}\ellc{\ts}{ij}+\ssum_{k:(j,k)\in\E}\Pc{\ts}{jk} \qquad &&\forall\ (i,j)\in\E\\
\label{eq:reacFlowSweep}		\Qc{\ts}{ij} &= \qtc{}{j}+\reactance{ij}\ellc{\ts}{ij}+\ssum_{k:(j,k)\in\E}\Qc{\ts}{jk}\qquad &&\forall\ (i,j)\in\E.
\end{alignat}

\textbf{Forward Sweep:} Starting from the children nodes of the substation node to the leaf nodes, compute $\forall\ (i,j)\in\E$:
\begin{equation}\label{eq:voltageSweep}
\nuc{\ts}{j} = \nuc{\ts}{i} - 2\left(\resistance{ij}\Pc{\ts}{ij}+\reactance{ij}\Qc{\ts}{ij}\right)+(\resistance{ij}^2+\reactance{ij}^2)\ellc{\ts}{ij}. 
\end{equation}

The BFS algorithm is bound to converge under mild assumptions of power flows in the DNs, for e.g., small line losses, small line impedances; see~\cite{shelarAminTCNS} for technical definitions of these assumptions. 

\begin{proof}[Proof of \cref{prop:relatingQuantitiesToNetConsumption}]
	Let $\{(\Pc{\ts}{},\Qc{\ts}{},\nuc{\ts}{},\ellc{\ts}{})\}_{\ts=1}^\nperiod$ be the values computed by the BFS algorithm in iteration $\ts=[1 \isep \nperiod]$ where $\nperiod$ is a fixed large number of iterations. Now, suppose that $\ptc{}{k}$ increases  marginally to $\ptc{}{k}+\Delta\ptc{}{k}$, while all other consumption values remain constant. Let $\{(\Pc{u\ts}{},\Qc{u\ts}{},\nuc{u\ts}{},\ellc{u\ts}{})\}_{\ts=1}^\nperiod$ be the new values computed by the BFS algorithm. 
	
	From \eqref{eq:realFlowSweep} and \eqref{eq:reacFlowSweep}, we get:
	%
	%
	%
	%
	%
	\begin{subequations}\label{eq:firstStep}
		\begin{alignat}{12}
			\Pc{u0}{ij} &= \Pc{0}{ij} + \Delta \ptc{}{k}\quad &&\forall\ (i,j) \in \pathFromRoot{k}\\
			\Pc{u0}{ij} &= \Pc{0}{ij} \quad &&\forall\ (i,j) \in \E\backslash\pathFromRoot{k}\\
			\Qc{u0}{ij} &= \Qc{0}{ij} &&\forall\ (i,j) \in \E. 
		\end{alignat}
	\end{subequations}
	By applying \eqref{eq:currentSweep} and \eqref{eq:firstStep}, we get
	\begin{subequations}\label{eq:secondStep}
		\begin{alignat}{12}
			\ellc{u1}{ij} &> \ellc{1}{ij}\qquad &&\forall\ (i,j)\in \pathFromRoot{k}\\
			\ellc{u1}{ij} &= \ellc{1}{ij}\quad &&\forall\ (i,j)\in\E\backslash \pathFromRoot{k}. 
		\end{alignat}
	\end{subequations}
Next, from \eqref{eq:voltageSweep} and \eqref{eq:secondStep}, we get:
	\begin{equation*}
			\quad\nuc{u1}{i} < \nuc{1}{i} \qquad \forall\ i\in\N, 
	\end{equation*}
 which, in turn, implies
	\begin{equation*}
		\quad \ellc{u2}{ij} > \ellc{2}{ij} \qquad \forall\ (i,j)\in\E. 
	\end{equation*} 
	
	Now, by making an inductive argument based on \eqref{eq:currentSweep}-\eqref{eq:voltageSweep}, we can show that 
	\begin{equation*}
		\ellc{u\ts}{ij} > \ellc{\ts}{ij} \quad\forall\ (i,j)\in\E,\ts\ge 2. 
	\end{equation*}
	Furthermore, we can also show that 
	\begin{equation}\label{eq:furtherProof}
	\ellc{u\ts}{ij} - \ellc{u\ts-1}{ij} > \ellc{\ts}{ij} - \ellc{\ts-1}{ij} \quad\forall\ (i,j)\in\E,\ts\ge 2. 
	\end{equation}
	(The proof of \eqref{eq:furtherProof} requires a further detailed analysis which is provided in \cite{technicalReport}.) 
	Thus, $\ellc{u\ts}{ij}$ and $\ellc{\ts}{ij}$ are the $\ts^{th}$ terms of two monotonically increasing and converging sequences such that the difference between consecutive terms of the former sequence are strictly greater than the corresponding difference of the latter. Therefore, the relative ordering also remains true for the converged values in the final iteration, i.e., 
	$\ellc{u}{ij} >\ellc{}{ij} \quad\forall\ (i,j)\in\E.$ 
	Then, by applying   \eqref{eq:realFlowFunction}-\eqref{eq:voltageFunction}, we can show that 
	\begin{alignat*}{8}
		\Pc{u}{ij} -\Pc{}{ij} &> \Pc{u0}{ij} - \Pc{0}{ij} &&\ge 0\qquad &&\forall\ (i,j)\in\E \\
		\Qc{u}{ij} -\Qc{}{ij} &> \Qc{u0}{ij} - \Qc{0}{ij} &&= 0\qquad &&\forall\ (i,j)\in\E \\
		\nuc{u}{i} -\nuc{}{i} &< \nuc{u0}{i} - \nuc{0}{i} &&< 0 \qquad&&\forall\ i\in\N.
	\end{alignat*}
	Then, taking the limit ${\Delta \ptc{}{k} \to 0}$, 
	\begin{equation*}
	\frac{\partial\Pc{}{ij}}{\partial\ptc{}{k}} > \frac{\partial\Pc{l}{ij}}{\partial\ptc{}{k}} \ge 0 > \frac{\partial\nuc{l}{l}}{\partial\ptc{}{k}} > \frac{\partial\nuc{}{l}}{\partial\ptc{}{k}}\quad  \forall\ (i,j)\in\E, l\in\N.
	\end{equation*}
	We conclude the proof by noting that a similar argument can be made had $\qtc{}{k}$ been increased instead of $\ptc{}{k}$.  
	%
	%
	%
	%
	%
	%
	%
	%
\end{proof}

\tcbtext{
\begin{proof}[Proof of \cref{prop:optimalDGoutput}]
	The proof follows from the application of   \cref{prop:relatingQuantitiesToNetConsumption}. Suppose that an optimal response  $(\lcc{\star}{},\kcc{\star}{},\kgc{\star}{},\pgc{\star}{},\qgc{\star}{})$ results in voltages $\nuc{\star}{}$ and currents $\ellc{\star}{}$. 
	Also, for the sake of contradiction, suppose that $\exists\ i\in\N$, $\pgc{\star}{i} < \pgc{max}{i}$. Thus,  increasing $\pgc{}{i}$ will increase the voltages and  reduce line losses. 
	Suppose, keeping everything else a constant, the operator changes his response to $\pgc{}{i} = \pgc{max}{i}$, which results in voltages $\nuc{}{}$ and currents $\ellc{}{}$.  Due to NRPF condition, the new voltage values will satisfy $\nuc{\star}{j}\le \nuc{}{j}\le \nuc{nom}{} \le \nuc{max}{j}\  \forall\ j \in\N$ and $\ellc{}{jk}\le \ellc{\star}{jk}$. Thus, the new response is feasible. Furthermore, the second and third terms in the objective function remain the same, whereas the first and last terms are strictly smaller for the newer response. This contradicts the optimality of $(\lcc{\star}{},\kcc{\star}{},\kgc{\star}{},\pgc{\star}{},\qgc{\star}{})$. 
\end{proof}
}

\begin{proof}[Proof of \cref{prop:downstreamImpactIsHigher}]
	Let  $(\Pc{\ts}{},\Qc{\ts}{},\nuc{\ts}{},\ellc{\ts}{})$ (resp. $(\Pc{u\ts}{},\Qc{u\ts}{},\nuc{u\ts}{},\ellc{u\ts}{}))$ be the values computed in $\ts^{th}$ iteration of the BFS algorithm when $\ptc{}{k}$ (resp. $\ptc{}{l}$) is increased by $\Delta \ptc{}{}$.  Applying \eqref{eq:realFlowLinear} and \eqref{eq:reacFlowLinear}, we get
		\begin{subequations}\label{eq:firstStepProp3}
		\begin{alignat}{12}
		\Pc{u0}{ij} &= \Pc{0}{ij} + \Delta \ptc{}{k}\quad &&\forall\ (i,j) \in \pathFromRoot{l}\backslash\pathFromRoot{k}\\
		\Pc{u0}{ij} &= \Pc{0}{ij} \quad &&\forall\ (i,j) \in \E\backslash(\pathFromRoot{l}\backslash\pathFromRoot{k})\\
		\Qc{u0}{ij} &= \Qc{0}{ij} &&\forall\ (i,j) \in \E. 
		\end{alignat}
	\end{subequations}
This is because when the consumption at $l$ increases, the additional power has to travel a path $\pathFromRoot{l}$ that subsumes the path $\pathFromRoot{k}$. 
	The rest of the proof is similar to that of  \cref{prop:relatingQuantitiesToNetConsumption}. Essentially, we again show that:
	\begin{equation*}
		\ellc{u}{ij} >\ellc{}{ij} \qquad\forall\ (i,j)\in\E, 
	\end{equation*}
	and, therefore, 
		\begin{alignat*}{8}
	\Pc{u}{ij} -\Pc{}{ij} &> 0, \quad 	\Qc{u}{ij} -\Qc{}{ij} > 0\qquad &&\forall\ (i,j)\in\E \\
	\nuc{u}{i} -\nuc{}{i} &< 0 \qquad&&\forall\ i\in\N.
	\end{alignat*}
	The proof is completed by taking the limit $\Delta \ptc{}{} \to 0$. 	
\end{proof}

\begin{proof}[Proof
	of \cref{prop:lossNPFisHigherThanlossLPF}]
	Let $(\third^\star, \xc{\star}{})$ be the optimal solution of the problem $\costMaxmin{}(\second)$. 
	For the fixed operator response $\third^\star$, the $\ptc{}{}$ and $\qtc{}{}$ vectors are uniquely determined. 
	Let $\xc{l}{}$ be the LPF solution for the $\ptc{}{}$ and $\qtc{}{}$ vectors. By applying \cref{prop:relatingQuantitiesToNetConsumption}, we can show that $\nuc{nom}{}\ge \nuc{l}{} \ge \nuc{}{}$. Therefore, we can claim that $(\third^\star,\xc{l}{})$ is a feasible solution for the problem $\costMaxmin{l}(\second)$. 
	
	Now, $\Lossc{}(\third^\star,\xc{\star}{}) - \Lossc{l}(\third^\star,\xc{l}{}) = \Clovr(\linfinityNorm{\nuc{nom}{} - \nuc{\star}{}} - \linfinityNorm{\nuc{nom}{} - \nuc{l}{}}) +  \Clineloss\ssum_{ij\in\E}\resistance{ij}\ellc{\star}{ij} \ge 0$, because both these terms are non-negative. 
%

	Let $\second^\star$ and $\hat{\second}^\star$ be the optimal attacker strategies to problems \eqref{eq:Mm-cascade} and \problemMaxminHat, respectively. Then, $\loss(\second^\star)\ge \loss(\hat{\second}^\star) \ge \losshat(\hat{\second}^\star)$, where the first inequality holds because of optimality of $\second^\star$, and the second inequality holds because of the first half of  \cref{prop:lossNPFisHigherThanlossLPF}. The proof completes by applying the definitions of $\resilienceMaxmin$ and $\resilienceMaxminhat$. 
\end{proof}

\tcbtext{
\begin{proof}[Proof of \cref{prop:monotonicity}]
For an attack $\second \in  \{0,1\}^{\N}$, the operator's subproblem involves minimization over the set $\Third(\second)$. If two attacks $\second', \second'' \in \{0,1\}^{\N}$ satisfy \eqref{eq:subsetAttack}, then the set of feasible operator strategies under $\second''$ is a subset of that under $\second'$, i.e., $\Third(\second'') \subseteq \Third(\second')$. Therefore, $\costMaxmin{}(\second')  \le \costMaxmin{}(\second'')$. 

Now, suppose that $\arcm',\arcm''$ are such that $0 \le \arcm' \le \arcm'' \le \N$. Furthemore, $\second'$ and $\second''$ are the optimal attacks  for attack cardinalities $\arcm'$ and $\arcm''$, respectively. We can construct an attack $\second''' \in \{0,1\}^{\N}$ such that $\abs{\second'''} = \arcm''$ and $\second'_i = 1 \implies \second'''_i = 1 \ \forall  \ i\in\N$. Then, 
\begin{equation*}
	\costMaxmin{}(\second') \le \costMaxmin{}(\second''') \le \costMaxmin{}(\second''), 
\end{equation*}
where the first inequality holds because $\second'$ and $\second'''$ satisfy \eqref{eq:subsetAttack}, and the second inequality holds because of the optimality of $\second''$ over attacks of cardinality $\arcm''$. The proof is completed by noting that
\begin{align*}\begin{aligned}
		\resilienceMaxmin^{\arcm'}  &= 100(1-\costMaxmin{}(\second')/\lcompleteShed)\\ 
		&\ge  100(1-\costMaxmin{}(\second'')/\lcompleteShed) = \resilienceMaxmin^{\arcm''}. 
	\end{aligned}\end{align*}
\end{proof}
}

\begin{proof}[Proof of \cref{prop:upstreamLoadConnectivityPreference}]
	Suppose for contradiction that $\third\in\Third$ is an optimal response such that $\kcc{}{i} = 1$, $\kcc{}{j} = 0$, and $\lcc{}{j} = a$ for some value $a\in [\lcc{min}{j},1]$. Then, we construct a response $\widecheck{\third}$ which is exactly the same as $\third$ except that $\kcc{u}{i} = 0$, $\kcc{u}{j} = 1$, and $\lcc{u}{i} = a$. Let $\xc{}{}$ and $\xc{u}{}$ be the corresponding network states. By \cref{prop:downstreamImpactIsHigher}, $\ellc{u}{} < \ellc{}{}$ and $\nuc{nom}{} > \nuc{u}{} > \nuc{}{} \ge \nucc{min}{}$. Therefore, $\xc{u}{}$ satisfies voltage bounds, and $\widecheck{\third}$ is a feasible operator strategy. Also, by \cref{prop:relatingQuantitiesToNetConsumption} and \cref{prop:downstreamImpactIsHigher}, the cost of voltage deviation and the cost of line loss is smaller because the increase in active and reactive load at $i$ (i.e. $a\pcc{max}{i}$ and $a\qcc{max}{i}$) is at most equal to the reduction in active and reactive load at $j$ ($a\pcc{max}{j}$ and $a\qcc{max}{j})$).  
	
	Now, the cost of load control and shedding in response $\widecheck{\third}$ is no worse than that in $\third$ (because $\Cload_j+ \Cload_i(1-a) \le \Cload_i+ \Cload_j(1-a)$). Moreover, the improved voltage profile may allow further reduction in cost of load control/shedding. Thus, $\third$ cannot be an optimal response. 
\end{proof}

\begin{proof}[Proof of \cref{prop:upstreamDGConnectivityPreference}]
	Suppose for contradiction that $\third\in\Third$ is an optimal response such that $\kgc{}{i} = 1$, and $\kgc{}{j} = 0$. Then, we can construct a response $\check{\third}$ which is exactly the same as $\third$ except that $\kgc{u}{i} = 0$, because DG $i$ was not disrupted by the attacker. Then, the cost of voltage deviations and line losses in response $\check{\third}$ is lesser than that in $\third$ by  \cref{prop:downstreamImpactIsHigher} and the fact that the decrease in active and reactive output of DG $i$ (i.e. $\pgc{max}{i}$ and $\etac{constant}{i}\pgc{max}{i}$) is smaller than the increase in active and reactive output of DG $j$ (i.e. $\pgc{max}{j}$ and $\etac{constant}{j}\pgc{max}{j}$). Thus, $\third$ cannot be an optimal response. 
\end{proof}

\clearpage




\section*{\Large \bf Technical Report}

This document provides additional discussions regarding our paper ``Evaluating Resilience of Electricity Distribution Networks via A Modification of Generalized Benders Decomposition Method", which we submitted to IEEE Transactions on Control of Network Systems. 

The outline for this report is as follows. In \cref{sec:introductionReport}, we describe the distinctions between different operations that the operator can exercise. In  \cref{sec:modelingDiscussion}, we provide a discussion on our attacker and operator modeling choices,  their technological feasibility, and their extensibility to other attacker and operator models. In~\cref{sec:problemFormulation}, we restate the formulations for the BiMISOCP, and provide an analogous formulation of the BiMILP problem defined using linear power flows. In~\cref{sec:technicalResults}, we provide additional details about the proof of Proposition 1 in the main manuscript. In \cref{sec:discusionOnBiMISOCP},  we describe the equivalence of the Min-cardinality and the Budget-k-max-loss problems, and provide additional details about our solution approach. Finally, in \cref{sec:computationalResults}, we add additional details about our computational study such as the network topology, and present a computational result, which we could not include in the main manuscript. 
\setcounter{subsection}{0}

\subsection{Introduction}
\label{sec:introductionReport}

The distinctions between operations (a) and response (b) and response (c) is summarized in \cref{tab:distinction}.

\begin{table}[h!]
	\resizebox{\textwidth}{!}{
		
		\begin{tabular}{|p{25mm}|x{25mm}|x{25mm}|x{25mm}|}
			\hline 
			Property & Operation (a) & Operation (b) & Operation (c)\\
			\hline
			Place of command initiation & Control center & DN node & Substation \\
			\hline
			Actions & Dispatch, fault/outage management & Disconnections & Load control, Disconnections, DG dispatch\\
			\hline
			Input & Node-level consumption, distributed generation, nodal voltages & Local nodal voltage & Node-level consumption, distributed generation, nodal voltages\\
			\hline
			Response time & 15 minutes or more & A few seconds to a few minutes & A few seconds \\
			\hline
			Coordinated & yes & no & yes \\
			\hline 
			Purpose & System-level optimization & Device protection & Prevention of network cascade\\
			\hline
			Attacked & \textcolor{red}{yes} & \textcolor{blue}{no} & \textcolor{blue}{no}\\
			\hline 
	\end{tabular}}
	\caption{Properties of operations  (a), (b), and (c).}\label{tab:distinction}
\end{table}

\subsection{Discussion on modeling choices}\label{sec:modelingDiscussion}


\subsubsection{Cyberphysical failure model} \label{subsec:cyberphysicalFailureModel}

Our model considers a generic cyber-physical failure model that captures the effects of DN-side component disruptions caused by security failures as well as effects of disturbances from the TN. Our model of TN-side disturbances is motivated by situations such as failure of a transmission line or a bulk generator, which result in low voltage conditions that last for a prolonged period (several minutes). We model its impact as a sudden drop in the substation node's voltage by $\vdc{}{0}$, which we assume to be exogenously given (and fixed). 
Indeed, $\vdc{}{0} = 0$ indicates no TN-side disturbance.\footnote{Note that a TN-side disturbance can also result in a change in frequency away from the nominal operating frequency of the network. In our future work, we extend our model to include frequency disturbances. }

On the other hand, our attack model is motivated by the security failure scenarios discussed in \cite{nescor}. Our attack model is relevant in the context of smart DNs, with a hierarchical control architecture; for further details we refer the reader to \cite{smartDistributionControl}. In this architecture, the main controller resides in the DN control center and performs the traditional tasks such as the optimization of DN operations and Volt-VAr control during nominal operations. Besides, it also provides flexibility to implement new functionalities such as DGMS. An attack on the DN control center server can affect one or more of these  functionalities.  
For the sake of concreteness, we limit our attention to a specific attack scenario in which the attacker targets the DGMS server, with the aim to simultaneously disrupt multiple DGs connected to the DN. However, our modeling approach is general in that it can also accommodate other important attack scenarios such as mass remote disconnects of loads or invalid load control commands~\cite{nescor}.\footnote{An attack on a DN control center can also be used to open circuit breakers. We consider this attack in the future work.}

Furthermore, our attack model considers that the control center functionalities such as DGMS are more viable targets for remote external attackers than local substation automation (SA) systems. Indeed, recent incidents \cite{ukraine} have confirmed that control center/DGMS servers can be targets of sophisticated phishing attacks (e.g. through a download of infected email attachments by the human operators who manage these servers). In contrast, a growing number of distribution utilities 
\ifLongVersion that implement the UVLS schemes
\fi
are regulated under NERC CIP standards which secure the substations against remote attacks via reperimetrisation of the substation cyber architecture~\cite{nerccip,Fuloria_theprotection}. In addition, SA is typically not prone to insecure actions by human insiders.

\begin{figure}[htbp!]
	\tikzstyle{lwstyle}=[line width=1.2pt]
	\tikzstyle{nwdth}=[draw,line width=1.2pt,align=left]
	\tikzstyle{nwdth1}=[nwdth, minimum width=3.4cm,inner sep=2]
	\tikzstyle{connector}=[->, lwstyle, >={latex[length=1pt]}]
	\begin{tikzpicture}[scale=0.8]
		\scalefont{0.7}
		\node[nwdth](ta) {Threat actors}; 
		\node[nwdth,right = 1cm of ta] (vu) {Vulnerabilities}; 
		\node[nwdth,right = 1cm of vu] (t) {Threats}; 
		\node[nwdth,fit = {(ta) (vu) (t)}, inner sep=5 ](avt) {};
		
		\node[nwdth,below = 0.2cm of avt, xshift=-3.5cm,align=center] (dc) {Disconnect\\commands}; 
		\node[nwdth,below = 0.4cm of dc,align=center,minimum height=0.5cm,minimum width=1.8cm] (dd) {DN-side \\disruptions}; 
		\node[nwdth,below = 0.3cm of dd,align=center,minimum height=0.5cm,minimum width=1.8cm] (td) {TN-side \\disturbances}; 
		
		\draw[connector] (avt.south) |- (dc.east); 
		
		\node[draw=red, right=0.35cm of dd, yshift=0.2cm, align=center,line width=1.5pt,minimum height=0.9cm,minimum width=2.2cm] (sdd) {Supply-demand\\disturbances\\at DN nodes};
		\node[draw=red, right=0.35cm of td, yshift=-0.2cm,align=center,line width=1.5pt,minimum height=0.9cm,minimum width=2.2cm] (svd) {Voltage \\disturbance at\\ substation node};
		
		\coordinate[below=0.15cm of dd, xshift=3.7cm] (pcli) {}; 
		
		\node[nwdth1,right=0.5cm of pcli,yshift=0.1cm] (lcs) {Operator response\\
			-- Load control, component \\ \quad disconnections\\
			(current version)\\
			-- microgrid islanding, \\ \quad DN restoration\\ (future work)}; 
		\node[nwdth1, above = 0.3cm of lcs](obv) {Contingencies}; 
		\node[nwdth1,below = 0.3cm of lcs] (pcl) {Post-contingency loss}; 
		\node[nwdth,fit = {(obv) (lcs) (pcl)}, inner sep=5] (impact) {};
		
		\foreach \from/\to in {ta/vu, vu/t,obv/lcs,lcs/pcl}
		\draw [connector, dashed] (\from) -> (\to);
		
		\foreach \from/\to in {dc/dd}
		\draw [connector] (\from) -> (\to);
		
		\foreach \from/\to in {dd/sdd, td/svd}
		\draw [connector] (\from.east) -> (\to.west);
		
		\foreach \from/\to in {sdd/pcli, svd/pcli}
		\draw [-,lwstyle] (\from.east) -| (\to);
		
		\draw [connector] (pcli) -- (impact);
		
	\end{tikzpicture}
	\caption[]{Framework for modeling impact of cyber-physical failures on DNs. The arrows in the top box indicate that threat actors exploit vulnerabilities of the system and pose risk of threats. One such threat is to issue disconnect commands leading to DN-side disruptions that cause supply-demand disturbances. Additionally, from the TN-side, certain events can cause voltage fluctuations which together with supply-demand disturbances can lead to contingencies. Then, the operator response to the adverse event determines the incurred post-contingency loss.}
	\label{fig:attackTrees}
\end{figure}
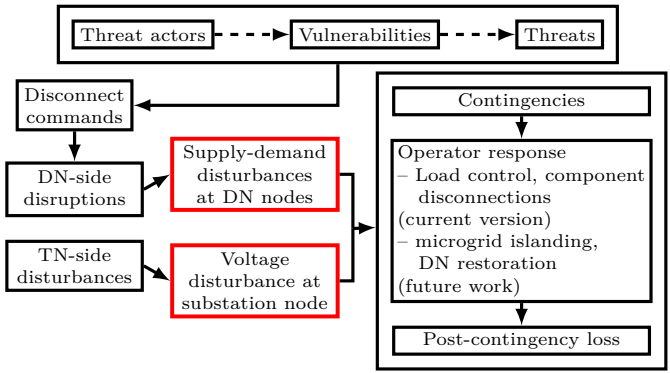

Our attack model is motivated by the security failure scenarios discussed in \cite{nescor}. These scenarios capture the capabilities of the following threat actors: (i) cyber-hackers of an enemy nation motivated to disrupt supply to critical facilities, 
(ii) a malicious adversary looking to extort ransom money from the utility, or (iii) a disgruntled internal employee motivated by revenge. In this paper, we are concerned with type (i) actors. Such actors can leverage existing vulnerabilities in DN cyber architecture such as 
non-confidentiality of control commands, lack of multi-factor authentication, and incorrect firewall rules that allow unauthorized access. Particularly, a threat actor can exploit these vulnerabilities to launch replay attacks~\cite{replayAttackZhuMartinez}, or a server-side attack at the control center, or hack operator credentials, any of which could allow him to perform malicious activities such as mass remote disconnect of components. We model the DN-side disruptions as nodal supply-demand disturbances. For example, mass disconnects of DGs (resp. loads) can cause loss of supply (resp. demand). Additionally, a threat actor could program his attack to be launched simultaneously with a TN-side disruption. A high-level framework for modeling impact of  cyber-physical disruptions to DN is illustrated in \cref{fig:attackTrees}. 

The disconnections of DGs and their inverters lead to a sudden drop in active as well as reactive power supply. Under heavy loading (high demand) conditions, reactive power often cannot be supplied from the bulk supply sources through the transmission lines. The reactive power shortfall may be exacerbated by a voltage dip resulting from a TN-side disturbance, as discussed below. This may result in sustained low-voltage conditions, e.g. a fault-induced delayed voltage recovery (FIDVR) event~\cite{voltageDisturbanceDuration,nercFIDVR} and/or result in voltage collapse. 

Now we model the impact of an attacker's actions on the DN state. If the attacker disrupts a DG at node~$i$, then that DG becomes non-operational, and is \emph{effectively disconnected} from the DN, i.e 
\begin{equation}
	\label{eq:dgConnectivityPostContingency1}
	\kgc{}{i} \ge \second_i \quad \forall\quad i\in \N. 
\end{equation}

\subsubsection{Discussion on Operator model} 

The emergency response capability (refer (c) in \cref{fig:attackerDefenderInteractions})  of modern SA systems is enabled by fine-grained data collection of node-level consumption, distributed generation, and nodal voltages. Many of the newer installations of smart meters are already equipped with data logging and communication capabilities. As a side note, the temporal frequency of data collected by low-voltage residential meters can vary from 15 minute to 24 hour intervals, depending on the desired control functionalities, customer privacy levels provided by the operator as well as the available communication bandwidth between DN nodes and the SA. In contrast, for the purpose of emergency response, meters installed at medium voltage to low voltage transformers at DN nodes can be utilized to provide aggregated node-level data from the customer meters in real-time (every second). With this capability, sudden changes in local DG output can also be detected by the SA, thereby enabling the operator to identify the attack vector $\second$. This level of monitoring does not involve individual customer meter readings, and hence, does not violate privacy regulations. 

Thus, the currently available capabilities of collection and processing of node-level data can be exploited by the operator to implement fast response strategies through SA. 

The best response by the DN operator can be computed sufficiently fast. Indeed, the window of opportunity which we allude to in Figure 2 of the manuscript can be of the order of thirty seconds or so for which the ride-through schemes are prescribed. In those 30 seconds, the operator may be able to detect the attacker's actions, by observing the sudden change in the active and reactive power consumption and generation via the Advanced Metering Infrastructures deployed at the DN medium-voltage nodes. Based on this knowledge, the operator can determine which nodes have been attacked. Then, the operator can solve the operator subproblem (MISOCP) within a few 100 milliseconds to determine an optimal response. 

\paragraph*{Traditional response to voltage regulation}
Indeed, other types of classical actions implemented through control of voltage regulators and capacitors as well as network reconfiguration can also form part of the operator response.  However, we chose load control and intentional disconnects due to timing requirements. 
The time-scale of the disturbance created by the attack can be very small (few seconds), and can trigger an immediate cascade of component disconnects due to operating bound violations. Typically, voltage regulators and capacitor banks require a longer response time; in fact, frequent activation of these devices is discouraged as they are subject to mechanical wear and tear~\cite{wearTearPVControl}. On the other hand, thanks to advances in SA and power electronics based control of loads/DGs, our response strategy can be implemented within a few milliseconds after the  information about the  timing and extent of the disruption is obtained by the SA. Our modeling approach can be extended to situations where appropriate changes in the settings of voltage regulators and capacitor banks are deemed to be  desirable aspects of operator response; these can be incorporated as integer decision variables in the inner problem of the considered bilevel formulation. 


We assume that the DN is connected to a \enquote{stiff} transmission network, barring the effect of transmission side disruption. As a result, one can argue that the transient effects arising due to attacker-operator interaction may not be as significant. 

\paragraph*{Loss function}
We have included the cost of load shedding, but not the cost of disconnection of customer-owned DGs because the customers are likely to face more inconvenience if there is load shedding, in comparison to DG disconnections during a contingency. However, we can easily account for the cost of DG disconnections in our formulation. 

Also, the cost of load control is modeled as an affine function as opposed to a quadratic function. Typically the incremental cost of load control should be larger when the distributed power decreases. Indeed, such a cost function can be handled provided the operator sub-problem is maintained as a SOCP.  However, we demonstrate that the performance of Benders decomposition approach in terms of its computational requirements significantly improves when the operator sub-problem is an MILP as opposed to when it is MiSOCP. Hence, to ensure that we can model the operator sub-problem as an MILP, we keep the cost of load control as an affine function.

\def \obj {L}
\def \lhsMatrix {A}
\def \rhsMatrix {B}


\subsection{Problem formulation}
\label{sec:problemFormulation}

\setcounter{subsubsection}{0}

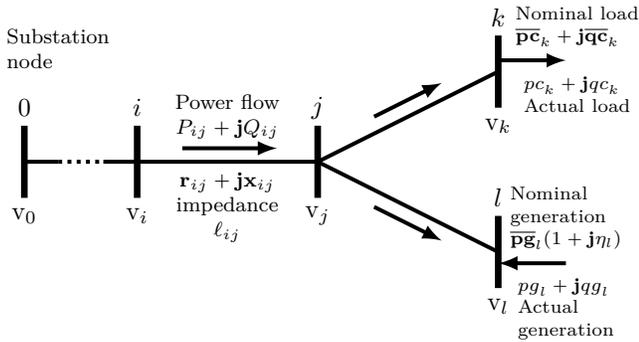
\begin{figure}[htbp!]
	\centering
	\tikzset{every node/.append ={font size=tiny}} 
	\tikzstyle{dnnode}=[draw,circle, minimum size=0.5pt, inner sep = 2]
	\tikzstyle{dnedge}=[-, line width=1pt]
	\tikzstyle{dernode}=[circle, fill=blue, minimum size=0.5pt, inner sep = 2]
	\tikzstyle{blackoutnode}=[circle, fill=black, minimum size=0.5pt, inner sep = 2]
	\tikzstyle{failededge}=[-, densely dotted]
	\def \drawgrid {\draw[step=1,gray, ultra thin, draw opacity = 0.5] (0,0) grid (3,4);}
	\def \drawSubstation {\draw[-, line width = 2pt] (0.8,4.05) -- (2.2,4.05)  node [midway,above] {};}
	\def \drawZero {\node (0) at (1.5,3.85) {};}
	\def \drawBus at (#1,#2); {\draw [line width=2.5pt,-] (#1,#2-0.8) -- (#1,#2+0.8);}
	\def \drawBuss at (#1); {\draw [line width=2.5pt,-] ([yshift=-0.7]#1) -- ([yshift=0.7]#1);}
	\def \drawLines at (#1,#2,#3,#4); {\draw [line width=1.5pt,-] (#1,#2) -- (#3,#4);}
	\def \drawArrows at (#1,#2,#3,#4); {\draw [line width=1.5pt,->, >=latex] (#1,#2) -- (#3,#4);}
	\scalefont{0.7}
	\begin{tikzpicture}[scale=0.6]
		
		\node (substation) at (0.5,0) {};
		
		\drawBus at (0.5,0); 	\drawBus at (3,0); \drawBus at (7,0); \drawBus at (11,2); \drawBus at (11,-2);
		
		\drawLines at (0.5,0,1.25,0); \drawLines at (2.25,0,7,0); \drawLines at (7,0,11,2); \drawLines at (7,0,11,-2); 
		\draw [line width=1.5pt,-,dotted] (1.25,0) -- (2.25,0);
		
		\drawArrows at (4,0.3,6.1,0.3); \drawArrows at (8.25,1,9.75,1+2/4*1.5); \drawArrows at (8.25,-1,9.75,-1-2/4*1.5); 
		
		\drawArrows at (11,2.25,12.5,2.25); 
		\node[align=left] at (12.75,3) {Nominal load\\ $\pcc{max}{k} + \j\qcc{max}{k}$};
		\node[align=left] at (12.75,1.5) {$\pcc{}{k} + \j\qcc{}{k}$\\Actual load};
		
		\drawArrows at (12.5,-2.25,11,-2.25); 
		\node[align=left] at (12.5,-3.25) { $\pgc{}{l} + \j\qgc{}{l}$\\
			Actual\\generation};
		\node[align=left] at (12.5,-1.25) {Nominal\\generation\\$\pgc{max}{l}(1 + \j\etac{constant}{l})$};
		
		\scalefont{1.4}
		\foreach \x/\y/\what  in {0.5/0/0, 3/0/i, 7/0/j, 11/2/k, 11/-2/l }
		{\node[] at (\x,\y+1.2) {$\what$}; \node[] at (\x,\y-1.2) {$\nuc{}{\what}$}; }
		
		\scalefont{0.8}
		\node[align=center] at (5,1) {Power flow\\$\Pc{}{ij} + \j\Qc{}{ij}$}; 
		\node[align=center] at (5,-1) {$\resistance{ij} + \j\reactance{ij}$\\impedance\\$\ellc{}{ij}$}; 
		\node[align=left] at (1.25,2.5) {Substation\\node}; 
		
		\scalefont{1}
		
	\end{tikzpicture}
	\caption{Illustration of a radial distribution network.}
	\label{fig:systemStateDistribution}
	
\end{figure}
\Cref{fig:systemStateDistribution} illustrates the topology and parameters of a radial distribution network pertaining to our problem.

\subsubsection{BiMILP formulation for $\lossMaxminhat$}
Solving a bilevel problem can be computationally difficult, especially when the inner subproblem is an MISOCP. To check whether using linear power flow (LPF) approximation provides any computational advantage, we propose an analogous Bilevel Mixed-Integer Linear Problem (BiMILP) based on LPF. Therefore, consider the classical \textit{LinDistFlow}  model~\cite{lindist}:
\begin{alignat}{2}
	\label{eq:conserveRealApprox} \Pc{}{ij} &= \ssum_{k:(j,k)
		\in\E} \Pc{}{jk} + \ptc{}{j} \qquad\qquad && \forall\ (i,j)\in\E 	\\
	\label{eq:conserveReactiveApprox} \Qc{}{ij} &= \ssum_{k:(j,k)
		\in\E}\Qc{}{jk} + \qtc{}{j} && \forall\ (i,j)\in\E	\\
	\label{eq:voltageApprox} \nuc{}{j} &=  \nuc{}{i} -  2\left(\resistance{ij}\Pc{}{ij} +  \reactance{ij}\Qc{}{ij}\right)  \qquad&& \forall\ (i,j)\in\E,
\end{alignat}
where  \crefrange{eq:conserveRealApprox}{eq:conserveReactiveApprox} are the approximate power conservation equations and \eqref{eq:voltageApprox} is the  voltage drop equation. 

We approximate the loss function in \eqref{eq:costGCregime} as the sum of following costs: (i) cost due to loss of voltage regulation, (ii) cost of load control, and (iii) cost of load shedding:
\begin{equation}\label{eq:costLinearApprox}\fontsize{9}{8}\selectfont
	\begin{split}
		\hspace{-0.3cm}\Lossc{l}(\third,\xc{}{}) = \Clovr\linfinityNorm{\nuc{nom}{} - \nuc{}{}} 
		+  \ssum_{i\in\N}\ \Cload _i\left(\unity-\lcc{}{i}\right)  \pcc{max}{i}\\ + \ssum_{i\in\N}\ \left(\Cshed_i-\Cload_i\right)\kcc{}{i}\pcc{max}{i},
	\end{split} 
\end{equation}
where we omit the line loss term for the sake of BiMILP formulation. 

Let $\Xc{l}{}$ denote the set of post-contingency states $\xc{}{}$ that satisfy the constraints \eqref{eq:postContingencyVoltage1}, \eqref{eq:loadControlParameterConsumptionConstraint}-\eqref{eq:totalPowerConsumption}, \eqref{eq:currentTrue}, and \eqref{eq:conserveRealApprox}-\eqref{eq:voltageApprox}. 
Again, we can denote the attacker-operator interaction under LPF constraints as follows:
\begin{align}\tag{$\widehat{\text{Mm}}$}\label{eq:Mm-cascade-linear}
	\begin{aligned}
		\hspace{-0.3cm}\lossMaxminhat\; \coloneqq\;  \max_{\second \in\Second}&  \quad \costMaxmin{l}(\second) \\
		\text{s.t.} & \quad \costMaxmin{l}(\second) \coloneqq \min_{\third\in\Third(\second), \xc{}{} \in \Xc{l}{}\left(\third\right)} \;  \Lossc{l}\left(\third,\xc{}{}\right).
	\end{aligned}
\end{align}
Note that current-magnitude-squared variables $\ellc{}{}$ do not affect the loss function $\Lossc{l}$, and do not impact the choice of other decision variables in  \eqref{eq:Mm-cascade-linear} as $\ellc{}{}$ only appear in \eqref{eq:currentTrue}. Hence, the problem \eqref{eq:Mm-cascade-linear} is still effectively a BiMILP despite having a non-linear equation~\eqref{eq:currentTrue}. 

\subsubsection{Features of our bilevel formulation}
Some features of Problem~\eqref{eq:Mm-cascade} are as  follows. We model the TN-side disruption as a sudden drop in substation voltage by $\vdc{}{0}$, which we assume to be exogenously given (and fixed). Indeed, $\vdc{}{0} = 0$ indicates no TN-side disturbance. In the attack model, we only consider disruption of DGs at nodes. However, our model is extensible to include attacks on loads. In our model, DGs may disconnect due to voltage bound violations. However, the DGs may also disconnect for other reasons such as frequency bound violations, which we will consider in our future work. The DG model is chosen such that there is no tradeoff between active and reactive power output of the DG. Our loss model can also be extended to include cost of DG disconnections, which we have not considered only for the sake of simplicity. The operator model can also be extended to include the traditional response mechanisms such as voltage regulators and capacitors. However, we do not consider them due to timing requirements.

\subsection{Technical results}\label{sec:technicalResults}
\setcounter{subsubsection}{0}
\subsubsection{Efficient computation of LPF and NPF solutions}
For fixed $\ptc{}{}$ and $\qtc{}{}$, let $\Pc{l}{}, \Qc{l}{}, \nuc{l}{}$ and $\ellc{l}{}$ be the LPF solutions of \eqref{eq:postContingencyVoltage1}, \eqref{eq:currentTrue} and  \crefrange{eq:conserveRealApprox}{eq:voltageApprox}. Since $\Pc{l}{}, \Qc{l}{}, \nuc{l}{}$ do not depend on $\ellc{l}{}$ and are linear functions of $\ptc{}{}$ and $\qtc{}{}$, $\Pc{l}{}, \Qc{l}{}, \nuc{l}{}$ and $\ellc{l}{}$ can be solved for in $\O(\abs{\N})$ time. 


Again, for fixed $\ptc{}{}$ and $\qtc{}{}$, let $(\Pc{}{},\Qc{}{},\nuc{}{},\ellc{}{})$ be the solution of the problem:
\begin{align}\label{eq:computeNPFquantitiesForFixedpq}
	\begin{aligned}
		\mmin_{\Pc{}{},\Qc{}{},\nuc{}{},\ellc{}{}} & \quad&&\Clovr\linfinityNorm{\nuc{nom}{} - \nuc{}{}} 
		+ \ssum_{(i,j)\in\E}\resistance{ij}\ellc{}{ij} \\
		\text{s.t.} &&& \eqref{eq:postContingencyVoltage1}, \eqref{eq:conserveRealTrue}-\eqref{eq:voltageTrue},  \eqref{eq:currentApproxConvex}.
	\end{aligned}
\end{align} Note that problem \eqref{eq:computeNPFquantitiesForFixedpq} is the same as the optimal power flow problem~\cite{exactConvexRelaxation} such that the lower and upper bounds for the net nodal consumption at each node are equal to $\ptc{}{i}$ and $\qtc{}{i}$. Furthermore, problem \eqref{eq:computeNPFquantitiesForFixedpq} is a SOCP and has a cost function that is strictly increasing in $\ellc{}{}$. Therefore, under NRPF, it has a unique solution~\cite{exactConvexRelaxation}.

Now, the objective in problem  \eqref{eq:computeNPFquantitiesForFixedpq} is strictly increasing in $(\ptc{}{},\qtc{}{},\ellc{}{})$ and $\ptc{}{i}$ and $\qtc{}{i}$ is fixed $\forall\ i\in\N$. Furthermore, we have assumed the NRPF condition. Hence, as shown in \cite{exactConvexRelaxation}, the solution of problem  \eqref{eq:computeNPFquantitiesForFixedpq} is unique. Thus, even $\ellc{}{}$ can be considered a function of $(\ptc{}{},\qtc{}{})$. 

The following lemma states the conditions under which the partial derivatives of the flow and voltage quantities can be defined. 

\begin{lemma}\label{lem:relatingPartialConsumptionGeneration} 
	Let $c\in \consumptionSet$ be a net nodal consumption quantity. Let $\hat{f} \in \Flowc{l}{}$ denote a flow quantity and $\hat{v}\in\Voltc{l}{}$ a voltage quantity computed using LPF. Let $f \in \Flowc{}{}$ and $v \in \Voltc{}{}$ be corresponding NPF quantities. The partial derivatives $\frac{\partial \hat{f}}{\partial c}$ and $\frac{\partial \hat{v}}{\partial c}$ exist with or without NRPF. Furthermore, under NRPF, the partial derivatives $\frac{\partial f}{\partial c}$ and $\frac{\partial v}{\partial c}$ also exist. 
	
	Consequently, the following hold:
	\begin{alignat*}{8}
		&\frac{\partial e}{\partial a} &&= \frac{\partial e}{\partial b} &&= -\frac{\partial e}{\partial c} \quad &&\forall \ e \in \Flowc{}{}\cup\Flowc{l}{}\cup\Voltc{}{}\cup\Voltc{l}{}, (a,b,c)\in\setConsumptionGenerationTuple.
	\end{alignat*}	
\end{lemma}


Henceforth, with a slight abuse of notation, we use the notation $\Third$ to denote the projection of the set $\{\third \in \R^{5\NN}  \text{ such that } \qgc{}{i} = \etac{}{i}\pgc{}{i} = \etac{}{i}\pgc{max}{i}(1-\kgc{}{i}) \quad \forall \ i \in \N \text{ and } \eqref{eq:integralityConstraints}-\eqref{eq:loadControlSheddingConstraint} \text{ hold}\}$ onto the space of $(\lcc{}{},\kcc{}{},\kgc{}{})-$variables. Then, an operator response can be denoted by $\third = (\lcc{}{},\kcc{}{},\kgc{}{}) \in\Third$.


\subsubsection{Detailed proof of \cref{prop:relatingQuantitiesToNetConsumption}}

For $i\in\N$, let $\pathFromRoot{i}\subseteq \E$ denote the subset of DN edges on the path from the substation node $0$ to node $i$. For $i,j\in\N$, let $\commonResistance{ij}$ (resp. $\commonReactance{ij}$) denote the sum of resistances (resp. reactances) of the edges common to $\pathFromRoot{i}$ and $\pathFromRoot{j}$, i.e., 
\begin{align*}
	\begin{aligned}
		\commonResistance{ij} &\coloneqq \ssum_{\left((k,l)\in\pathFromRoot{i}\cap\pathFromRoot{j}\right)} \resistance{kl}, \qquad \\
		\commonReactance{ij} &\coloneqq \ssum_{\left((k,l)\in\pathFromRoot{i}\cap\pathFromRoot{j}\right)} \reactance{kl}. 
	\end{aligned}
\end{align*}

Let $\subtree{i}\subseteq \N$ be the subset of nodes that form the subtree rooted at node $i$, which includes node $i$, and let $\subtreeEdges{i}\subseteq\E$ be the subset of edges that form the subtree $\subtree{i}$. Then, the following equations can be derived using recursion on the radial tree topology. 
\begin{alignat}{8}
	\label{eq:realFlowLinear}	\Pc{l}{ij} &= \ssum_{k\in\subtree{j}}\ptc{}{k}  && \forall\ (i,j)\in\E \\
	\label{eq:reacFlowLinear}	\Qc{l}{ij} &= \ssum_{k\in\subtree{j}}\qtc{}{k} && \forall\ (i,j)\in\E\\
	\label{eq:voltageLinear}	\nuc{l}{j} &= \nuc{nom}{} - 2\ssum_{k}\left(\commonResistance{jk}\ptc{}{k} +\commonReactance{jk}\qtc{}{k}\right)  &\quad& \forall\ j\in\N \\
	\label{eq:realFlowFunction}	\Pc{}{ij} &= \Pc{l}{ij} +  \ssum_{(k,l)\in\subtreeEdges{i}}\resistance{kl}\ellc{}{kl}  && \forall\ (i,j)\in\E \\
	\label{eq:reacFlowFunction}	\Qc{}{ij} &= \Qc{l}{ij} +  \ssum_{(k,l)\in\subtreeEdges{i}}\reactance{kl}\ellc{}{kl}  && \forall\ (i,j)\in\E 
\end{alignat}
\begin{equation}\label{eq:voltageFunction}
	\begin{split}
		\nuc{}{j} = \nuc{l}{j} - 2\ssum_{(k,l)\in\E}(\commonResistance{jl}\resistance{kl} + \commonReactance{jl}\reactance{kl})\ellc{}{kl} + \\ \ssum_{(k,l)\in\pathFromRoot{j}}(\resistance{kl}^2 +\reactance{kl}^2)\ellc{}{kl}
	\end{split}
\end{equation}
Thus, we can write $(\Pc{l}{},\Qc{l}{},\nuc{l}{})$ as functions of $(\ptc{}{},\qtc{}{})$ and  $(\Pc{}{},\Qc{}{},\nuc{}{})$ as functions of $(\ptc{}{},\qtc{}{},\ellc{}{})$.

Consider the iterative Backward-Forward Sweep (BFS) algorithm~\cite{backwardForwardSweep} used to compute the NPF values, which we modify to consider the TN-side voltage disturbance. Let  $(\Pc{\ts}{},\Qc{\ts}{},\nuc{\ts}{},\ellc{\ts}{})$ be the values computed in $\ts^{th}$ iteration of the FBS algorithm. 

\textbf{Initialization:} 
\begin{alignat}{8}
	\nuc{0}{i} &= \nuc{nom}{}-\vdc{}{0}\quad &&\forall\ i\in\N\backslash\{0\}\\
	\nuc{\ts}{0} &= \nuc{nom}{}-\vdc{}{0}\quad &&\forall\ \ts\in[1\isep \nperiod]\\
	\ellc{0}{ij} &= 0, \Pc{0}{ij} = \Pc{l}{ij}, \Qc{0}{ij} = \Qc{l}{ij} \qquad &&\forall\ (i,j)\in\E. 
\end{alignat}

\textbf{Backward Sweep:} Starting from the leaf nodes to the substation node, compute:
\begin{alignat}{8}
	\label{eq:currentSweep}		\ellc{\ts}{ij} &= \left({(\Pc{\ts-1}{ij})}^2+{(\Qc{\ts-1}{ij})}^2\right)/\nuc{\ts-1}{i} \qquad &&\forall\ (i,j)\in\E\\ 
	\label{eq:realFlowSweep}		\Pc{\ts}{ij} &= \ptc{}{j}+\resistance{ij}\ellc{\ts}{ij}+\ssum_{k:(j,k)\in\E}\Pc{\ts}{jk} \qquad &&\forall\ (i,j)\in\E\\
	\label{eq:reacFlowSweep}		\Qc{\ts}{ij} &= \qtc{}{j}+\reactance{ij}\ellc{\ts}{ij}+\ssum_{k:(j,k)\in\E}\Qc{\ts}{jk}\qquad &&\forall\ (i,j)\in\E.
\end{alignat}

\textbf{Forward Sweep:} Starting from the children nodes of the substation node to the leaf nodes, compute $\forall\ (i,j)\in\E$:
\begin{equation}\label{eq:voltageSweep}
	\nuc{\ts}{j} = \nuc{\ts}{i} - 2\left(\resistance{ij}\Pc{\ts}{ij}+\reactance{ij}\Qc{\ts}{ij}\right)+(\resistance{ij}^2+\reactance{ij}^2)\ellc{\ts}{ij}. 
\end{equation}

The BFS algorithm is bound to converge under mild assumptions of power flows in the DNs, for e.g., small line losses, small line impedances; see~\cite{shelarAminTCNS} for technical definitions of these assumptions.

\begin{proof}[Proof of \cref{prop:relatingQuantitiesToNetConsumption}]
	Let $\{(\Pc{\ts}{},\Qc{\ts}{},\nuc{\ts}{},\ellc{\ts}{})\}_{\ts=1}^\nperiod$ be the values computed by the BFS algorithm in iteration $\ts=[1 \isep \nperiod]$ where $\nperiod$ is a fixed large number of iterations. Now, suppose that $\ptc{}{k}$ increases  marginally to $\ptc{}{k}+\Delta\ptc{}{k}$, while all other consumption values remain constant. Let $\{(\Pc{u\ts}{},\Qc{u\ts}{},\nuc{u\ts}{},\ellc{u\ts}{})\}_{\ts=1}^\nperiod$ be the new values computed by the BFS algorithm. 
	
	From \eqref{eq:realFlowSweep} and \eqref{eq:reacFlowSweep}, we get:
	%
	%
	%
	%
	%
	\begin{subequations}\label{eq:firstStep}
		\begin{alignat}{12}
			\Pc{u0}{ij} &= \Pc{0}{ij} + \Delta \ptc{}{k}\quad &&\forall\ (i,j) \in \pathFromRoot{k}\\
			\Pc{u0}{ij} &= \Pc{0}{ij} \quad &&\forall\ (i,j) \in \E\backslash\pathFromRoot{k}\\
			\Qc{u0}{ij} &= \Qc{0}{ij} &&\forall\ (i,j) \in \E. 
		\end{alignat}
	\end{subequations}
	By applying \eqref{eq:currentSweep} and \eqref{eq:firstStep}, we get
	\begin{subequations}\label{eq:secondStep}
		\begin{alignat}{12}
			\ellc{u1}{ij} &> \ellc{1}{ij}\qquad &&\forall\ (i,j)\in \pathFromRoot{k}\\
			\ellc{u1}{ij} &= \ellc{1}{ij}\quad &&\forall\ (i,j)\in\E\backslash \pathFromRoot{k}. 
		\end{alignat}
	\end{subequations}
	Next, from \eqref{eq:voltageSweep} and \eqref{eq:secondStep}, we get:
	\begin{equation*}
		\nuc{u1}{i} < \nuc{1}{i} \qquad \forall\ i\in\N, 
	\end{equation*}
	which, in turn, implies
	\begin{equation*}
		\ellc{u2}{ij} > \ellc{2}{ij} \qquad \forall\ (i,j)\in\E. 
	\end{equation*} 
	
	Now, by making an inductive argument based on \eqref{eq:currentSweep}-\eqref{eq:voltageSweep}, we can show that 
	\begin{equation*}
		\ellc{u\ts}{ij} > \ellc{\ts}{ij} \quad\forall\ (i,j)\in\E,\ts\ge 2. 
	\end{equation*}
	Furthermore, we can also show that 
	\begin{equation}\label{eq:motonicalDifference}
		\ellc{u\ts}{ij} - \ellc{u\ts-1}{ij} > \ellc{\ts}{ij} - \ellc{\ts-1}{ij} \quad\forall\ (i,j)\in\E,\ts\ge 2. 
	\end{equation}
	The detailed argument for the previous inequality \eqref{eq:motonicalDifference} is as follows. 
	
	Let $f_j = \frac{{\Pc{}{j}}^2+{\Qc{}{j}}^2}{\partial \nuc{}{i}}$. Also, let the shorthand for $\frac{\partial a}{\partial ch\ellc{}{k}}$ be written as $\partial_{\ellc{}{k}} a$. Then,	
	\begin{align*}
		\begin{aligned}
			\partial_{\ellc{}{k}}\Pc{}{j}  &= \resistance{k}\unity\{k\in\subtree{j}\} && \forall\ j \in\N\\
			\partial_{\ellc{}{k}}\Qc{}{j}  &= \reactance{k}\unity\{k\in\subtree{j}\} && \forall\ j \in\N\\
			\partial_{\ellc{}{k}}\nuc{}{i}  &= -\left(\commonResistance{ik}\resistance{k}+\commonReactance{ik}\reactance{k}\right)(2-\unity\{k\in\P_i\}) && \forall\ i \in\N,
		\end{aligned}
	\end{align*}
	where $\unity\{\}$ is an indicator function. 
	
	Therefore, 
	\begin{align*}
		\begin{aligned}
			\partial_{\ellc{}{k}} f_j &= \frac{2\left(\Pc{}{j}\partial_{\ellc{}{k}}\Pc{}{j} + \Qc{}{j}\partial_{\ellc{}{k}}\Qc{}{j}\right)}{\nuc{}{i}} - \frac{\left({\Pc{}{j}}^2+{\Qc{}{j}}^2\right)\partial_{\ellc{}{k}}\nuc{}{i}}{{\nuc{}{i}}^2}\\
			&=  \frac{2\left(\resistance{k}\Pc{}{j} + \reactance{k} \Qc{}{j}\right)\unity\{k\in\subtree{j}\}}{\nuc{}{i}} \\
			&\quad +  \frac{\left({\Pc{}{j}}^2+{\Qc{}{j}}^2\right)\left[\left(\commonResistance{ik}\resistance{k}+\commonReactance{ik}\reactance{k}\right)(2-\unity\{k\in\P_i\})\right]}{{\nuc{}{i}}^2}
		\end{aligned}
	\end{align*}
	
	Now, suppose $\ellc{u}{}$ and $\ellc{l}{}$ are such that $\ellc{u}{k} > \ellc{l}{k}\quad \forall\ k \in \N$, then $\Pc{u}{j} > \Pc{l}{j}$ and $\Qc{u}{j} > \Qc{l}{j} \quad \forall \ j \in \N$ as well as $\nuc{u}{i} < \nuc{l}{i} \quad \forall \ i \in \N$. Since the resistances and reactances are positive,
	\begin{align*}
		\begin{aligned}
			2\left(\resistance{k}\Pc{u}{j} + \reactance{k} \Qc{u}{j}\right) &> 2\left(\resistance{k}\Pc{l}{j} + \reactance{k} \Qc{l}{j}\right),\\
			\left[\left(\commonResistance{ik}\resistance{k}+\commonReactance{ik}\reactance{k}\right)(2-\unity\{k\in\P_i\})\right] &> 0, \\
			\text{and}\qquad		{\Pc{u}{j}}^2+{\Qc{u}{j}}^2 &> {\Pc{l}{j}}^2+{\Qc{l}{j}}^2.
		\end{aligned}
	\end{align*}
	Therefore, 
	\begin{equation*}
		\partial_{\ellc{}{k}} f_j \bigg\rvert_{\ellc{}{} = \ellc{u}{}} \quad > \quad  \partial_{\ellc{}{k}} f_j \bigg\rvert_{\ellc{}{} = \ellc{l}{}} \qquad 
		\forall\ j,k \in \N, 
	\end{equation*}
	i.e., the partial derivative of $f_j$ with respect to $\ellc{}{k}$ when evaluated at $\ellc{}{} = \ellc{u}{}$ is greater than the partial derivative of $f_j$ with respect to $\ellc{}{k}$ when evaluated at $\ellc{}{} = \ellc{l}{}$. 
	
	(Here, the detailed argument for the proof of  \eqref{eq:motonicalDifference} ends, and we return to the rest of the proof of \cref{prop:relatingQuantitiesToNetConsumption}. )
	
	Thus, $\ellc{u\ts}{ij}$ and $\ellc{\ts}{ij}$ are the $\ts^{th}$ terms of two monotonically increasing and converging sequences such that the difference between consecutive terms of the former sequence are strictly greater than the corresponding difference of the latter. Therefore, the relative ordering also remains true for the converged values in the final iteration, i.e., 
	\begin{equation*}
		\ellc{u}{ij} >\ellc{}{ij} \quad\forall\ (i,j)\in\E. 
	\end{equation*}
	Then, by applying   \eqref{eq:realFlowFunction}-\eqref{eq:voltageFunction}, we can show that 
	\begin{alignat*}{8}
		\Pc{u}{ij} -\Pc{}{ij} &> \Pc{u0}{ij} - \Pc{0}{ij} &&\ge 0\qquad &&\forall\ (i,j)\in\E \\
		\Qc{u}{ij} -\Qc{}{ij} &> \Qc{u0}{ij} - \Qc{0}{ij} &&= 0\qquad &&\forall\ (i,j)\in\E \\
		\nuc{u}{i} -\nuc{}{i} &< \nuc{u0}{i} - \nuc{0}{i} &&< 0 \qquad&&\forall\ i\in\N.
	\end{alignat*}
	Then, taking the limit ${\Delta \ptc{}{k} \to 0}$, 
	\begin{equation*}
		\frac{\partial\Pc{}{ij}}{\partial\ptc{}{k}} > \frac{\partial\Pc{l}{ij}}{\partial\ptc{}{k}} \ge 0 > \frac{\partial\nuc{l}{l}}{\partial\ptc{}{k}} > \frac{\partial\nuc{}{l}}{\partial\ptc{}{k}}\quad  \forall\ (i,j)\in\E, l\in\N.
	\end{equation*}
	We conclude the proof by noting that a similar argument can be made had $\qtc{}{k}$ been increased instead of $\ptc{}{k}$.  
\end{proof}

\subsection{Discussion on BiMISOCP formulation}
\label{sec:discusionOnBiMISOCP}

\setcounter{subsubsection}{0}
\subsubsection{Example for problem \eqref{eq:innerConvexSP}}

\def \obj {L}
\def \lhsMatrix {A}
\def \lhssocpMatrix {E}
\def \lhssocpVector {f}
\def \rhssocpCoeff {g}
\def \rhssocpScalar {h}
\def \rhsMatrix {B}

Note that problem \eqref{eq:innerConvexSP} with parameters ($\second^\star$, $\configurationc{}{}^\star$) can be simplified and rewritten as the following problem:
\begin{align}
	\nonumber\mmin_{w} \  &\transpose{c}w &&\\
	\text{s.t. } &\lhsMatrix w &&\ge b + \rhsMatrix\second^\star && \tag{O-SOCP2}\label{eq:OSP2}\\
	& \nonumber\normSquared{\lhssocpMatrix^iw+\lhssocpVector^i} &&\le \transpose{\rhssocpCoeff^i}w + \rhssocpScalar^i && \ \forall\ i \in [1\isep\NN],
\end{align}
where $\normSquared{\cdot}$ is the L-squared norm; $w$ is the primal decision vector variable; $\lhsMatrix$, $\rhsMatrix$, and $\lhssocpMatrix^i$ for $i\in [1\isep\NN]$ are matrices; $b$,  $\lhssocpVector^i$ and $\rhssocpCoeff^i$ for $i\in [1\isep\NN]$ are vectors of appropriate dimensions; and ${\rhssocpScalar^i}s$ are scalars. The $\NN$ second-order cone constraints correspond to \eqref{eq:currentApproxConvex}. 

The dual of problem \eqref{eq:OSP2} is as follows:
\begin{align}\tag{D-SOCP2}\label{eq:dspSP2}\fontsize{9}{10}\selectfont
	\hspace{-0.5cm}\begin{aligned}
		\max_{\lambda\ge\zero}   &&& \transpose{\left(b+\rhsMatrix\second^\star\right)}\lambda+\sum_{i=1}^\NN \big( \transpose{\lhssocpVector^i}\alpha^i-\beta^i\rhssocpScalar^i\big)\\
		\text{s.t. }&&& c - \transpose{A}\lambda+\sum_{i=1}^\NN \big(\transpose{\lhssocpMatrix^i}\alpha^i - \beta^i \rhssocpCoeff^i\big) = \zero&& \\
		&&& \normSquared{\alpha^i} \le \beta^i \qquad\forall\ i \in [1\isep\NN]
	\end{aligned}
\end{align}
Here $w$ and $\lambda$ are the primal and dual decision vector variables; $\lhsMatrix = \transpose{[\transpose{\lhsMatrix_{eq}} \transpose{\lhsMatrix_{in}}]}$, $\rhsMatrix = 
\transpose{[\transpose{\rhsMatrix_{eq}} \transpose{\rhsMatrix_{in}}]}
$ are matrices and $b= 
\transpose{[\transpose{b_{eq}} \transpose{b_{in}}]}$ is a vector of appropriate dimensions.


Recall the primal problem in \eqref{eq:OSP2}. With the help of an example, we show how to instantiate the primal problem. Consider a DN $\G$ with nodes $\{0,1\}$ and line $(0,1)$. Then the variable $w$ is given as: 
\begin{equation*}
	w = (\lcc{}{1}, \pgc{}{1},\qgc{}{1}, \ptc{}{1},\qtc{}{1},\Pc{}{1}, \Qc{}{1}, \nuc{}{0}, \nuc{}{1}, \ellc{}{1}, t), 
\end{equation*} 
where $t$ is an auxiliary variable. The corresponding cost vector $c$ is given as:
\begin{equation}
	c = (-\Cload, 0, 0, 0, 0, \Csupply, 0, 0, 0, 0, \Clineloss, \Clovr). 
\end{equation}
Furthermore, we are given the parameters $\second$ and $\configurationc{}{} = (\kcc{}{1},\kgc{}{1})$. Then, the constraints of the problem \eqref{eq:dspSP2} are given as follows:  
\ifOneColumn
\def \ineqheight {3.2cm}
\else
\def \ineqheight {2.0cm}
\fi 
\setcounter{MaxMatrixCols}{20}
\begin{alignat*}{11}\setstretch{1}
	\overbrace{\resizebox{!}{\ineqheight}{
			$\begin{bmatrix}
				0 & 1 & 0 & 0 & 0 & 0 & 0 & 0 & 0 & 0 & 0\\
				0 & -1 & 0 & 0 & 0 & 0 & 0 & 0 & 0 & 0 & 0\\
				0 & \eta_1 & 1 & 0 & 0 & 0 & 0 & 0 & 0 & 0 & 0\\
				0 & -\eta_1 & -1 & 0 & 0 & 0 & 0 & 0 & 0 & 0 & 0\\
				0 & 1 & 0 & 0 & 0 & 0 & 0 & 0 & 0 & 0 & 0\\
				0 & -1 & 0 & 0 & 0 & 0 & 0 & 0 & 0 & 0 & 0\\
				0 & 0 & 0 & 0 & 0 & 0 & 0 & 0 & 1 & 0 & 0\\
				0 & 0 & 0 & 0 & 0 & 0 & 0 & 0 & -1 & 0 & 0\\
				0 & 0 & 0 & 0 & 0 & 0 & 0 & 0 & 1 & 0 & 0\\
				0 & 0 & 0 & 0 & 0 & 0 & 0 & 0 & -1 & 0 & 0\\
				1 & 0 & 0 & 0 & 0 & 0 & 0 & 0 & 0 & 0 & 0\\
				-1 & 0 & 0 & 0 & 0 & 0 & 0 & 0 & 0 & 0 & 0\\
				0 & 0 & 0 & 0 & 0 & 0 & 0 & 0 & 1 & 0 & 1\\
				0 & 0 & 0 & 0 & 0 & 0 & 0 & 0 & -1 & 0 & 1\\
				0 & 0 & 0 & 0 & 0 & 0 & 0 & 0 & 0 & 0 & 1\\
			\end{bmatrix}$
	}}^{A_{in}}w \ge 
	\overbrace{\resizebox{!}{\ineqheight}{$\begin{bmatrix}
				0\\
				-\pgc{max}{1}\\
				0\\
				0\\
				(1-\kgc{}{1})\pgc{max}{1}\\
				-(1-\kgc{}{1})\pgc{max}{1}\\
				\nugc{min}{1} - \kgc{}{1}\\
				-\nugc{max}{1} + \kgc{}{1}\\
				\nucc{min}{1} - \kgc{}{1}\\
				-\nucc{max}{1} + \kgc{}{1}\\
				(1-\kcc{}{1})\lcc{min}{1}\\
				-(1-\kcc{}{1})\\
				\nuc{nom}{}\\
				-\nuc{nom}{}\\
				0
			\end{bmatrix}$}}^{b_{in}} + 
	\overbrace{\resizebox{!}{\ineqheight}{$\begin{bmatrix}
				0\\
				\pgc{max}{1}\\
				0\\
				0\\
				0\\
				0\\
				0\\
				0\\
				0\\
				0\\
				0\\
				0\\
				0\\
				0\\
				0
			\end{bmatrix}$}}^{\rhsMatrix_{in}} \second
\end{alignat*}

\ifOneColumn
\def \ineqheight {1.2cm}
\else
\def \ineqheight {0.6cm}
\fi 
\begin{alignat*}{10}\setstretch{1}
	\overbrace{\resizebox{!}{\ineqheight}{$\begin{bmatrix}
				0 & 0 & 0 & -1 & 0 & 1 & 0 & 0 & 0 & \resistance{01} & 0\\
				0 & 0 & 0 & 0 & -1 & 0 & 1 & 0 & 0 & \reactance{01} & 0\\
				0 & 0 & 0 & 0 & 0 & 2\resistance{01} & 2\reactance{01} & -1 & 1 & -(\resistance{01}^2+\reactance{01}^2) & 0\\
				-\pcc{max}{1} & 1 & 0 & 1 & 0 & 0 & 0 & 0 & 0 & 0 & 0\\
				-\qcc{max}{1} & 0 & 1 & 0 & 0 & 0 & 0 & 0 & 0 & 0 & 0\\
			\end{bmatrix}$}}^{A_{eq}}w =
	\overbrace{\resizebox{!}{\ineqheight}{$\begin{bmatrix}
				0\\
				0\\
				0\\
				0\\
				0
			\end{bmatrix}$}}^{b_{eq}} + 
	\overbrace{\resizebox{!}{\ineqheight}{$\begin{bmatrix}
				0\\
				0\\
				0\\
				0\\
				0
			\end{bmatrix}$}}^{\rhsMatrix_{eq}} \second
\end{alignat*}
\ifOneColumn
\def \ineqheight {1.2cm}
\else
\def \ineqheight {1.6cm}
\fi 
\begin{alignat}{10}\setstretch{1}
	\normSquared{\overbrace{\resizebox{!}{\ineqheight}{$\begin{bmatrix}
					0 \\
					0\\
					0\\
					0\\
					0\\
					\sqrt{2}\\
					\sqrt{2}\\
					1\\
					0\\
					1\\
					0
				\end{bmatrix}$}}^{\transpose{E^1}}w 
		+ 
		\overbrace{\resizebox{!}{\ineqheight}{$\begin{bmatrix}
					0\\
					0\\
					0\\
					0\\
					0\\
					0\\
					0\\
					0\\
					0\\
					0\\
					0
				\end{bmatrix}$}}^{\transpose{f^1}} } \le 
	\overbrace{\resizebox{!}{\ineqheight}{$\begin{bmatrix}
				0\\
				0\\
				0\\
				0\\
				0\\
				0\\
				0\\
				1\\
				0\\
				1\\
				0
			\end{bmatrix}$}}^{\transpose{g^1}}w + 
	\overbrace{\resizebox{!}{\ineqheight}{$\begin{bmatrix}
				0\\
				0\\
				0\\
				0\\
				0\\
				0\\
				0\\
				0\\
				0\\
				0\\
				0
			\end{bmatrix}$}}^{h^1} \label{eq:normcurrentequation}
\end{alignat}
Note that the previous equation \eqref{eq:normcurrentequation} is equivalent to 
\begin{align*}
	\begin{aligned}
		& \normSquared{(\sqrt{2}P_1,\sqrt{2}Q_1,\nuc{}{0},\ellc{}{1})} &&\le (\nuc{}{0}+\ellc{}{1})\\
		\iff & 2P_1^2+2Q_1^2 + {\nuc{}{0}}^2 + {\ellc{}{1}}^2 &&\le {\nuc{}{0}}^2 + {\ellc{}{1}}^2 + 2{\nuc{}{0}}{\ellc{}{1}}\\
		\iff & \ellc{}{1} &&\ge \frac{P_1^2+Q_1^2}{\nuc{}{0}}
	\end{aligned}
\end{align*}

Finally, $\lhsMatrix = \begin{bmatrix}
	\lhsMatrix_{in}\\\lhsMatrix_{eq}
\end{bmatrix}$, $\rhsMatrix = \begin{bmatrix}
	\rhsMatrix_{in}\\\rhsMatrix_{eq}
\end{bmatrix}$, and $b=\begin{bmatrix}
	b_{in}\\b_{eq}
\end{bmatrix}$.

\subsubsection{Equivalence of min-cardinality problem and Budget-k-max-loss problem}
\label{subsec:equivalence}

The min-cardinality problem is equivalent to \eqref{eq:Mm-cascade} in the following sense. The loss $\lossMaxmin$ in \eqref{eq:Mm-cascade} is non-decreasing with $\arcm$ (due to the inequality constraint $\sum_{i\in\N}\second_i \le \arcm$). Therefore, if the parameter $\ltarget$ is gradually increased then the minimum attack cardinality computed by the min-cardinality problem will be non-decreasing in $\ltarget$. Thus, for a fixed budget $\arcm$, the smallest $\ltarget$ value at which the minimum attack cardinality changes from $\arcm$ to $\arcm+1$ will be the optimal value of problem \eqref{eq:Mm-cascade}. By implementing a binary search on the parameter $100\ltarget/\lcompleteShed$ between $0-100\%$, we can determine the smallest $\ltarget$ at which the minimum attack cardinality changes from $\arcm$ to $\arcm+1$. Conversely, if we can solve \eqref{eq:Mm-cascade}, then by implementing a binary search on the parameter $\arcm$ between $0$ and $\NN$, we can determine the minimum attack cardinality whose optimal loss exceeds $\ltarget$. However, \eqref{eq:mincardinalityProblem} is a relatively  easier problem to solve using GBD because the master problem of \eqref{eq:mincardinalityProblem} has fewer variables and constraints than the corresponding master problem of \eqref{eq:Mm-cascade}. 

The quantity $100\ltarget/\lcompleteShed$ is also relevant from the viewpoint of DN resilience. For example, if we want to evaluate whether or not a DN is 80\% resilient to a $\arcm$ cardinality attack, we can set $\ltarget = 0.2\lcompleteShed$, and then check if the optimal value of the min-cardinality problem is smaller than or equal $\arcm$. 

\subsubsection{Generalized Benders Cut}
Note that problem \eqref{eq:innerConvexSP} with parameters ($\second^\star$, $\configurationc{}{}^\star$) can be simplified and rewritten as the following problem:
\begin{align}
	\nonumber\min_{w} \  &\transpose{c}w &&\\
	\text{s.t. }  &\lhsMatrix w &&\ge b + \rhsMatrix\second^\star && : (\lambda) \tag{O-SOCP2}\label{eq:OSP2}\\
	& \nonumber\normSquared{\lhssocpMatrix^jw+\lhssocpVector^j} &&\le \transpose{\rhssocpCoeff^j}w + \rhssocpScalar^j && : (\alpha^j,\beta^j) &&\ \forall\ j \in \N,
\end{align}
where $\normSquared{\cdot}$ is the L-squared norm; $w$ is the primal decision vector variable; $\lhsMatrix$, $\rhsMatrix$, and $\lhssocpMatrix^j$s are matrices; $b$,  $\lhssocpVector^j$s and $\rhssocpCoeff^j$s are vectors of appropriate dimensions; and ${\rhssocpScalar^j}s$ are scalars. The $\NN$ second-order cone constraints correspond to \eqref{eq:currentApproxConvex}. Also, $\lambda$ and ($\alpha^j, \beta^j$)  for $j \in \N$ are the dual variables corresponding to the linear and SOCP inequalities, respectively. 

The dual of problem \eqref{eq:OSP2} is as follows:
\begin{align}\tag{D-SOCP2}\label{eq:dspSP2}\fontsize{9}{10}\selectfont
	\hspace{-0.45cm}\begin{aligned}
		\max_{\substack{\lambda\ge\zero,\\\alpha^j,\beta^j}}   &&& \transpose{\left(b+\rhsMatrix\second^\star\right)}\lambda \\
		\text{s.t. }&&& c - \transpose{A}\lambda+\sum_{j\in\N} \big(\transpose{\lhssocpMatrix^j}\alpha^j - \beta^j \rhssocpCoeff^j\big) = \zero&& \\
		&&& \normSquared{\alpha^j} \le \beta^j \qquad\forall\ j \in \N
	\end{aligned}
\end{align}
We solve the dual problem (thanks to strong duality, the optimal values are the same) in  \eqref{eq:dspSP2} to compute $\ploss\left({\second}^\star,{\configurationc{}{}}^\star\right)$ and an optimal dual solution $(\lambda^\star,\alpha^{j\star},\beta^{j\star})$. This furnishes a \emph{generalized Benders cut}, which is added to the master problem in the next iteration. 
In particular, if the dual problem in \eqref{eq:dspSP2} has an optimal solution $(\lambda^\star,\alpha^{j\star},\beta^{j\star})$, and its optimal value is $\obj^\star$, then 
\begin{equation}\label{eq:generalizedBendersCut}
	\transpose{\left(b+\rhsMatrix\second\right)}\lambda^\star \ge \obj^\star + \epsilon 
\end{equation} 
is the desired generalized Benders cut  where $\epsilon$ is a non-negative number. In a classical generalized Benders cut the value of $\epsilon$ is 0. If the inner subproblem of \eqref{eq:Mm-cascade} were convex, such a cut would indeed be useful in eliminating sub-optimal attacker strategies~\cite{generalizedBenders}. However, this cut is not useful in the presence of discrete inner variables, i.e. it does not eliminate any attack vector. 

Hereafter, we refer to the generalized Benders cut as simply the Benders cut. Note that $\second^\star$ does not satisfy the Benders cut constraint when $\epsilon>0$ because $\transpose{\left(b+\rhsMatrix\second^\star\right)}\lambda^\star+\sum_{j\in\N} \big( \transpose{\lhssocpVector^j}\alpha^j-\beta^j\rhssocpScalar^j\big) = \ploss\left(\second^\star, {\configurationc{}{}}^\star\right) = \obj^\star < \obj^\star + \epsilon$, where the first equality holds because of strong duality in second-order cone programs. Thus, choosing $\epsilon >0$ is a modification to the Benders cut which will help eliminate $\second^\star$ from attacker's set of feasible strategies. 
However, due to numerical issues,  an off-the-shelf solver \enquote{stalls} sometimes, and is unable to generate dual vector values required for the Benders cut. To address this issue, we add the following cut: 
\begin{equation}
	\ssum_{\left(i\in\N:\second^\star_i = 1\right)}\second_i + \ssum_{\left(i\in\N:\second^\star_i = 1\right)} (1-\second_i) \le \NN-1, 
\end{equation} 
which definitely eliminates $\second^\star$.

\subsection{Exact expression of Benders cut}
\def \l#1{\lambda^{(\text{#1})}}

\setcounter{subsubsection}{0}
An equality constraint $a = b$ can be reformulated as $a \ge b$ and $-a \ge -b$. Now, suppose that the equality constraints in  \eqref{eq:postContingencyVoltage1}-\eqref{eq:voltageTrue} are similarly reformulated. Then, let the dual variables corresponding to the constraints \eqref{eq:postContingencyVoltage1}-\eqref{eq:voltageTrue} are as follows: \\
$\l{2a} , \l{2b} , \l{4} , \l{5} , \l{6} , \l{7a}, 
\l{7b} , \l{9a} $, $\l{9b} , \l{9c} , $\\
$\l{9d} , \l{10a} , \l{10b} , \l{10c} , \l{10d} , \l{11a} , \l{11b} ,\l{12a} ,  $\\ 
$\l{12b} ,\l{13a} , \l{13b} , \l{13c} ,\l{13d} , \l{14a} , \l{14b} \l{15a} ,$ \\
$\l{15b} , \l{16a} , \l{16b} $.\\
Here, for example, $\l{2a} $ and $\l{2b} $ are the dual variables correspond to the constraint \eqref{eq:postContingencyVoltage1}. Also, the two sets of equalities in \eqref{eq:totalPowerConsumption} would correspond to four sets of inequalities. Hence, the dual variables $\l{10a} , \l{10b} , \l{10c} $, and $\l{10d} $.

The vector $b$ can be written as follows: \\
$b = [\nuc{nom}{}- \Delta \nuc{}{0}, -\nuc{nom}{} + \Delta \nuc{}{0}, -\pgc{max}{}, -\pgc{max}{}, \zero_{3\NN}, \lcc{min}{}$, $-\unity_{\NN}, \zero_{4\NN}, \nugc{min}{} - \kgc{}{}, -\nugc{max}{} - \kgc{}{}, \nucc{min}{} - \kcc{}{}, -\nucc{max}{} - \kcc{}{}, \zero_{10\NN}]$. 

The vector $(B\second) = [0, 0, \pgc{max}{}\emult\second, \zero_{24\NN}]$. 

Then, the Benders cut \eqref{eq:generalizedBendersCut} can be written as follows: \\
$(\l{2a} - \l{2b} )(\nuc{nom}{}- \Delta \nuc{}{0})  +\sum_{i\in\N} \big(\lambda^{(4)}_i (-\pgc{max}{i} +  \pgc{max}{i}\second_i) - \lambda^{(5)}_i\pgc{max}{i} + \lambda^{(9a)}_i\lcc{min}{i} - \lambda^{(9b)}_i + \lambda^{(11a)}_i (\nugc{min}{i} - \kgc{}{i}) +  \lambda^{(11b)}_i (-\nugc{max}{i} - \kgc{}{i}) + \lambda^{(12a)}_i (\nucc{min}{i} - \kcc{}{i}) +  \lambda^{(12b)}_i (-\nucc{max}{i} - \kcc{}{i})\big) \ge \obj^\star  \epsilon$

\subsubsection{Comparative remarks about solution approach}
We now offer some comparative remarks about our solution approach to \eqref{eq:Mm-cascade} 	which -- as mentioned earlier -- is a BiMISOCP with conflicting objectives in the inner (operator) and outer (attacker) problems. In general, one can reformulate a BiMISOCP as a single level MISOCP (for example, via a high-point relaxation (HPR) problem~\cite{bard1990,xuWangBnBBiMIP}), and use an advanced branch-and-bound algorithm to solve the problem. Note, however, the  HPR is a weak relaxation of the original BiMISOCP due to directly conflicting objectives~\cite{kevinwood,baldickBowenHua}. More recent work has developed intersection cuts~\cite{fischetti2016new,fischetti2018new} and disjunction cuts~\cite{disjuctionCutsWaterMelon,disjunctionCutsValueFunction} -- these approaches introduce stronger cuts for the HPR problem. However, these approaches are only suitable for BiMISOCPs in which the inner problem has integer coefficients in the constraints. In contrast, our problem \eqref{eq:Mm-cascade} has fractional coefficients. A recent paper by Hua et. al~\cite{baldickBowenHua} addresses this issue by applying a Generalized Benders decomposition method but without the min-cardinality reformulation; as a result, the master problem in their approach needs to handle a relatively larger number of variables and constraints. 
Since in our solution approach we apply the Min-cardinality reformulation, the resulting master problem has fewer variables and constraints. Another approach by Zeng and An~\cite{zeng2014solving} uses a Column Constraint Generation (CCG) method, whose iterations progressively add variables and constraints (particularly the disjuntive constraints resulting from the KKT conditions for the inner problem with fixed binary variables). While these approaches are certainly of interest in solving \eqref{eq:Mm-cascade}, we find that our proposed approach achieves desirable computational performance as discussed in the case study in the main paper. 

\subsection{Evaluating $\resilienceNoResponse$ - A Two-Step Approach}\label{sec:autonomousDisconnects}

\setcounter{subsubsection}{0}
We restate the intermediate and final problems for our two-step approach here. The intermediate problem is as follows:
\begin{align}\tag{P-IN}\label{eq:noresponseIntermediate}
	\begin{aligned}
		&& \mmin_{\third_\intermediateState,\xc{}{\intermediateState}} &\  &&\Lossc{}\left(\third_\intermediateState,\xc{}{\intermediateState}\right) \\ 
		& &\text{s.t. }&&&
		\third_\intermediateState \in \Third(\second), \quad && \xc{}{\intermediateState}\in\Zc{}{}\left(\third_\intermediateState\right)\\
		&&&&& \lcc{in}{i} = 1 \quad &&  \forall \ i\in\N,
	\end{aligned}
\end{align} 
The problem to compute the final state under the autonomous disconnections is as follows:  
\begin{align}\tag{P-FN}\label{eq:noresponseFinal} 
	\begin{aligned}
		\hspace{-0.3cm}
		\min_{\third_\text{nr},\xc{}{\text{nr}}} &&& \Lossc{}\left(\third_\noresponseSmall,\xc{}{\noresponseSmall}\right)\\ 
		\text{s.t. } &&& \third_\noresponseSmall \in \Third, \quad && \xc{}{\noresponseSmall}\in\Xc{}{}\left(\third_\noresponseSmall\right) \\
		&&&\lcc{nr}{i} = \kcc{nr}{i} &&\forall \ i\in\N\\
		&&& \kgc{nr}{i} \ge \kgc{in\star}{i}(\second) &&\forall \ i\in\N\\
		&&& \kcc{nr}{i} \ge \nuc{}{i}-\nuc{in\star}{i}(\second) \;  &&\forall \ i\in\N
	\end{aligned}
\end{align} 
The optimal solution of the above problem will provide us $(\third^{\star}_\noresponseSmall,\xc{\star}{\noresponseSmall})$, i.e the final autonomous disconnect action and the post-contingency state. 
\vspace{1cm}
\begin{algorithm}[tbp!]
	\caption{Uncontrolled cascade under response (b)}
	\begin{algorithmic}[1]
		\Require attacker action $\second$ (initial contingency)
		\State $\third^\star_\noresponseSmall,\xc{\star}{\noresponseSmall}\gets \Call{GetCascadeFinalState($\second$)}{}$
		\Function{GetCascadeFinalState}{$\second$}
		\State Compute  $\third^\star_\intermediateState(\second),\xc{\star}{\intermediateState}(\second)$ by solving \eqref{eq:noresponseIntermediate}
		\State Extract parameters $(\kgc{in\star}{},\nuc{in\star}{})$ from $(\third^\star_\intermediateState,\xc{\star}{\intermediateState})$
		\State Instantiate \eqref{eq:noresponseFinal} with parameters $(\kgc{in\star}{},\nuc{in\star}{})$
		\State Solve \eqref{eq:noresponseFinal} to compute the final state  $\third^\star_\noresponseSmall,\xc{\star}{\noresponseSmall}$
		\State \Return $\third^\star_\noresponseSmall,\xc{\star}{\noresponseSmall}$ 
		\EndFunction
	\end{algorithmic}
	\label{algo:noDefenderResponse}
\end{algorithm}

\paragraph*{Remarks about  \cref{algo:noDefenderResponse}}
In \cref{algo:noDefenderResponse}, DGs disconnect before the loads disconnect. This can be justified by considering a sudden voltage drop. Such voltage behavior can be indicative of a fault within the DN, and therefore DGs supplying power to a fault can be potentially dangerous. Therefore, according to \cite{ieee1547}, when voltage bound violations occur, the DGs are supposed to disconnect within two seconds or less, depending on the extent of the voltage drop. On the other hand, the loads can continue to operate even a minute after mild or moderate voltage bound violations occur. Indeed, we can infer this from the fact that the response time of voltage regulators along DN feeders is typically at least 15 to 30 seconds~\cite{voltageRegulatorResponseTime,voltageDisturbanceDuration}. However, the disconnect actions of loads happen due to activation of protection devices which operate based on local measurements, i.e. they operate independent of each other. Therefore, in the worst case all loads experiencing voltage bound violations may disconnect together. Hence, our choice to consider the disconnection of all the DGs followed by the disconnection of all loads is reasonable. 

\subsection{Computational studies} \label{sec:computationalResults}

\paragraph*{Topologies of modified IEEE test networks}

\begin{figure}[h!]
	\tikzstyle{ndc} = [draw,circle,inner sep=0,minimum width=0.55cm]
	\tikzstyle{ndcSyn} = [ndc,fill,gray,opacity=0.2]
	\tikzstyle{ndcVsi} = [ndc,fill,pattern=vertical lines]
	\tikzstyle{ndcattack} = [ndc,fill,line width=2pt,red,opacity=0.3]
	\tikzstyle{mgl} = [draw]
	\subfloat[24 node network.]{
		\resizebox{8.7cm}{!}{
			\begin{tikzpicture}
				\def \nodePar {0/1,1/2,2/3,3/4,4/5,5/6,6/7,7/8,8/9,9/10,10/11,11/12,12/13,13/14,14/15,15/16,16/17,17/18,18/19,19/20,20/21,21/22,22/23,23/24,24/25,25/26,26/27,27/28,28/29,29/30,30/31,31/32,32/33,33/34,34/35,35/36}
				
				\def \ny {0}		
				\node[ndc] (0) at (0,\ny) {0};		
				\foreach \x in {1,...,12}
				{
					\node[ndc](\x) at (\x,\ny) {\x};
				}
				
				\def \ny {1}		
				\foreach \x in {13,...,16}
				{
					\pgfmathtruncatemacro{\nx}{\x-10}
					\node[ndc](\x) at (\nx,\ny) {\x};
				}
				\def \ny {-1}		
				\foreach \x in {17,...,24}
				{
					\pgfmathtruncatemacro{\nx}{\x-13}
					\node[ndc](\x) at (\nx,\ny) {\x};
				}
				
				\def \synNodes {{1,3,5,7,9,11,13,15,17,19,21,23}}
				\foreach \x in {1,3,5,7,9,11,13,15,17,19,21,23}
				\node[ndcSyn] at  (\x) {};
				
				\foreach \x/\y in {0/1,1/2,2/3,3/4,4/5,5/6,6/7,7/8,8/9,9/10,10/11,11/12,2/13,13/14,14/15,15/16,3/17,17/18,18/19,6/20,20/21,21/22,22/23,23/24}
				\draw (\x) -- (\y);
				\foreach \x/\y in {0/1,6/7,2/13,3/17,6/20}
				\draw[mgl] (\x) -- (\y);
	\end{tikzpicture}}}\qquad\\
	\subfloat[36 node network.]{\resizebox{6cm}{!}{
			\begin{tikzpicture}
				\foreach \z/\x/\y in {0/0/0,1/1/0,2/2/0,16/3/0,21/4/0,22/5/0,23/6/0,24/7/0,6/2/1,7/2/2,8/3/2,9/4/2,14/1/2,15/0/3,13/1/3,10/2/3,11/3/3,12/4/3,4/1/-1,3/2/-1,5/2/-2,17/3/-1,18/3/-2,19/4/-2,20/3/-3,25/5/-1,26/5/-2,27/6/-1,28/7/-1,29/7/-2,31/6/-2,30/8/-2,32/8/-1,33/8/0,34/8/1,35/8/2,36/7/1}
				\node[ndc] (\z) at (\x,0.8*\y) {\z};
				\foreach \x in {2, 6, 7, 11, 12, 13, 16, 18, 19, 21, 23, 25, 28, 29, 32, 34, 35, 36}
				\node[ndcSyn] at  (\x) {};
				\foreach \x/\y in {0/1,1/2,2/3,3/4,3/5,2/6,6/7,7/8,8/9,7/10,10/11,11/12,10/13,13/14,13/15,2/16,16/17,17/18,18/19,18/20,16/21,21/22,22/23,23/24,22/25,25/26,25/27,27/28,28/29,29/30,29/31,28/32,32/33,33/34,34/35,34/36}
				\draw (\x) -- (\y);
				\foreach \x/\y in {0/1,2/6,7/10,2/16,16/21,25/27,28/32}
				\draw[mgl] (\x) -- (\y);
			\end{tikzpicture}	
	}}\\
	\subfloat[118 node network.]{
		\resizebox{8.7cm}{!}{
			\begin{tikzpicture}
				
				\def \lateral(#1,#2,#3,#4) {
					\foreach \x in {#1,...,#2}
					{
						\pgfmathtruncatemacro{\nx}{\x-#3}
						\node[ndc](\x) at (\nx,0.8*#4) {\x};
					}
					\pgfmathtruncatemacro{\ns}{#1+1}
					\foreach \x in {\ns,...,#2}
					{
						\pgfmathtruncatemacro{\ef}{\x-1}
						\draw (\x) -- (\ef);				
					}
				}
				
				
				\lateral(18,27,16,7)
				\lateral(10,17,9,6)
				\lateral(4,9,2,5)
				\node[ndc] (2) at (1,5*0.8) {2};
				\node[ndc] (3) at (0,5*0.8) {3};
				\lateral(38,46,35,4)
				\lateral(28,35,26,3)
				\lateral(47,54,42,2)
				\lateral(36,37,33,2)
				\lateral(55,62,53,1)
				\lateral(96,99,92,0)
				\lateral(0,1,0,0)
				\lateral(89,95,87,-1)
				\lateral(63,69,62,-2)
				\lateral(70,77,66,-3)
				\lateral(78,80,77,-3)
				\lateral(81,85,77,-4)
				\lateral(86,88,85,-4)
				\lateral(100,106,100,-5)
				\lateral(107,113,102,-6)
				\lateral(114,118,114,-6)
				
				\foreach \x in {1,3,5,7,9,11,13,15,17,19,21,23,25,27,29,31,33,35,37,39,41,43,45,47,49,51,53,55,57,59,61,63,65,67,69,71,73,75,77,79,81,83,85,87,89,91,93,95,97,99,101,103,105,107,109,111,113,115,117}
				\node[ndcSyn] at  (\x) {};
				\foreach \x/\y in {1/2,2/3,2/4,4/28,30/36,80/81,100/114,64/78,79/86,65/89,91/96}
				\draw (\x) -- (\y);
				\foreach \x/\y in {0/1,2/10,11/18,29/30,29/38,35/47,29/55,1/63,65/89,69/70,79/80,106/107,1/100,100/101}
				\draw[mgl] (\x) -- (\y);
	\end{tikzpicture}}}
	\caption{Modified IEEE test networks. DG nodes are  indicated by gray color. }\label{fig:testNetworks}
\end{figure}
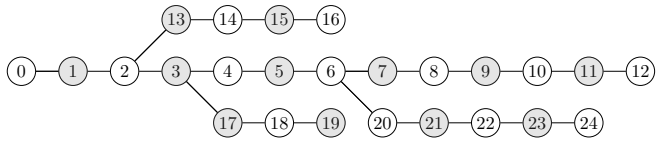
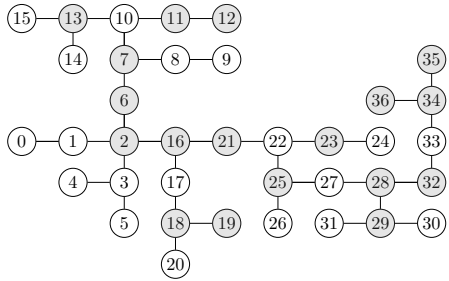
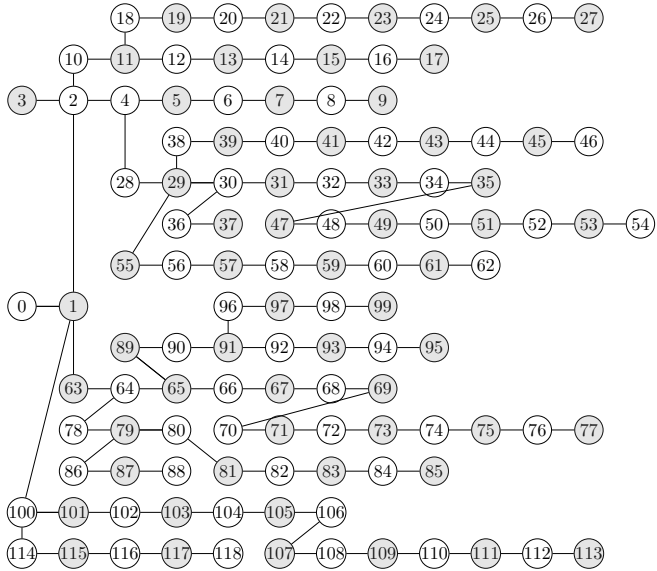

The modified IEEE test networks are shown in \cref{fig:testNetworks}. 

\normalsize

\end{document}